\documentclass[AMA,Times1COL]{WileyNJDv5} 

\articletype{Research article}%

\received{Date Month Year}
\revised{Date Month Year}
\accepted{Date Month Year}
\journal{Journal}
\volume{00}
\copyyear{2025}
\startpage{1}

\usepackage{bm}
\usepackage{bbm}
\usepackage[version=4]{mhchem}
\usepackage{siunitx}
\usepackage{longtable,tabularx}
\usepackage{xspace}

\usepackage{booktabs}
\usepackage{subcaption} 

\usepackage{graphicx}

\usepackage{float}
\floatstyle{plaintop}
\restylefloat{table}

\usepackage{tikz}

\usepackage{subcaption}

\usepackage{cancel}

\usepackage{algorithm,algpseudocode}

\usepackage{dsfont}

\usepackage{amsthm}

\usepackage{stackengine} 

\usepackage{lipsum}

\newcommand{\Table}[1]{\mbox{Table~\ref{#1}}}

\newcommand{\Fig}[1]{\mbox{Fig.~\ref{#1}}}
\newcommand{\Figs}[2]{\mbox{Figs.~\ref{#1}} and \ref{#2}}
\newcommand{\Figss}[3]{\mbox{Figs.~\ref{#1}}, \ref{#2}, and \ref{#3}}

\newcommand{\Sec}[1]{\mbox{Section~\ref{#1}}}
\newcommand{\Secs}[2]{\mbox{Sections~\ref{#1}} and \ref{#2}}

\newcommand{\Eq}[1]{\mbox{Eq.~(\ref{#1})}}
\newcommand{\Eqs}[2]{\mbox{Eqs.~(\ref{#1})} and (\ref{#2})}
\newcommand{\Eqss}[3]{\mbox{Eqs.~(\ref{#1})}, (\ref{#2}), and (\ref{#3})}

\newcommand{\bsym}[1]{\boldsymbol{#1}}

\newcommand{\mat}[1]{\ensuremath{\mathsf{#1}}}

\newcommand{\Idty}[0]{\mat{I}} 

\definecolor{Orange}{RGB}{255,69,0}
\definecolor{Red}{RGB}{235,40,25}

\definecolor{Green}{RGB}{34,139,34}

\newcommand\p[2]{\frac{\partial #1}{\partial #2}}
\newcommand\pp[3]{\frac{\partial^2 #1}{\partial #2 \partial #3}}

\newcommand{\order}[1]{\mathcal{O}\left(#1\right)}

\DeclareMathOperator{\cov}{cov}

\DeclareMathOperator{\mydiag}{diag}

\newcommand\lb[1]{\ensuremath{\underline{#1}}}
\newcommand\ub[1]{\ensuremath{\overline{#1}}}

\newcommand{\argmin}[0]{\ensuremath{\operatornamewithlimits{argmin}}}

\newcommand{\med}[0]{\operatorname{med}}

\newcommand{\condmax}[0]{\kappa_{\max}}

\newcommand{\eg}[0]{e.g.\ }
\newcommand{\ie}[0]{i.e.\ }

\usepackage{amssymb}
\usepackage{pifont}
\newcommand{\ignore}[1]{} 

\newcommand\yesnumber{\addtocounter{equation}{1}\tag{\theequation}}

\newcommand{\RNum}[1]{\uppercase\expandafter{\romannumeral #1\relax}}

\newcommand{\zero}[0]{\ensuremath{\textbf{0}}}

\newcommand{\one}[0]{\ensuremath{\textbf{1}}}

\newcommand{\avec}[0]{\ensuremath{\bsym{a}}}

\newcommand{\fvec}[0]{\ensuremath{\bsym{f}}}

\newcommand{\kvec}[0]{\ensuremath{\bsym{k}}}
\newcommand{\lvec}[0]{\ensuremath{\bsym{l}}}
\newcommand{\mvec}[0]{\ensuremath{\bsym{m}}}

\newcommand{\rvec}[0]{\ensuremath{\bsym{r}}}
\newcommand{\svec}[0]{\ensuremath{\bsym{s}}}

\newcommand{\uvec}[0]{\ensuremath{\bsym{u}}}
\newcommand{\vvec}[0]{\ensuremath{\bsym{v}}}

\newcommand{\xvec}[0]{\ensuremath{\bsym{x}}}
\newcommand{\yvec}[0]{\ensuremath{\bsym{y}}}

\newcommand{\betavec}[0]{\ensuremath{\bsym{\beta}}}

\newcommand{\phivec}[0]{\ensuremath{\bsym{\phi}}}

\newcommand{\gammavec}[0]{\bsym{\gamma}}

\newcommand{\A}[0]{\mat{A}}

\newcommand{\E}[0]{\mat{E}}

\newcommand{\K}[0]{\mat{K}}

\newcommand{\W}[0]{\mat{W}}
\newcommand{\X}[0]{\mat{X}}

\newcommand{\noptz}[0]{\ensuremath{n_{\text{optz}}}}
\newcommand{\nLhs}[0]{\ensuremath{n_{\text{LHS}}}}
\newcommand{\nmed}[0]{\ensuremath{n_{\text{med}}}}

\newcommand{\rdot}[0]{\ensuremath{\dot{r}}}

\newcommand{\rvecdot}[0]{\ensuremath{\dot{\rvec}}}

\newcommand{\fgrad}{\ensuremath{\bsym{f}_{\nabla}}}

\newcommand{\mgrad}{\ensuremath{\bsym{m}_{\nabla}}}

\newcommand{\onemod}[0]{\ensuremath{\check{\one}}}

\newcommand{\kG}[0]{\ensuremath{k_{\text{G}}}}
\newcommand{\kM}[0]{\ensuremath{k_{\text{M}{\frac{5}{2}}}}}
\newcommand{\kRq}[0]{\ensuremath{k_{\text{rq}}}}

\newcommand{\kvecgrad}{\ensuremath{\bsym{k}_\nabla}}

\newcommand{\Kg}[0]{\ensuremath{\mat{K}_{\nabla}}}

\newcommand{\Kgdot}[0]{\ensuremath{\dot{\mat{K}}_{\nabla}}}

\newcommand{\etaKg}{\ensuremath{\eta_{\Kg}}}
\newcommand{\etaKgdot}{\ensuremath{\eta}_{\Kgdot}}

\newcommand{\Vg}[0]{\ensuremath{\mat{V}_{\nabla}}}

\newcommand{\Lmat}[0]{\mat{L}}

\newcommand{\Sigmag}[0]{\ensuremath{\Sigma_{\nabla}}}
\newcommand{\Sigmagdot}[0]{\ensuremath{\dot{\Sigma}_{\nabla}}}

\newcommand{\Pmat}[0]{\mat{P}}
\newcommand{\Pinv}[0]{\mat{P}^{-1}}

\newcommand{\sigK}{\ensuremath{\hat{\sigma}_{\K}}}
\newcommand{\sigKf}{\ensuremath{\hat{\sigma}_{\K,f}}}

\newcommand{\hpstdfval}{\ensuremath{\hat{\sigma}_f}}
\newcommand{\hpstdfgrad}{\ensuremath{\hat{\sigma}_{\nabla f}}}

\newcommand\muGP{\ensuremath{\tilde{\mu}}}
\newcommand\sigmaGP{\ensuremath{\tilde{\sigma}}}

\newcommand\Xall[0]{\ensuremath{\X_{\text{all}}}}

\newcommand\Xdata[0]{\ensuremath{\X_{\text{data}}}}
\newcommand\Xacq[0]{\ensuremath{\X_{\text{acq}}}}

\newcommand\Jdata[0]{\ensuremath{\bsym{J}_{\text{data}}}}

\newcommand\noisy[1]{\ensuremath{\widetilde{#1}}} 

\newcommand\pdf[0]{\ensuremath{\theta_\text{pdf}}}
\newcommand\cdf[0]{\ensuremath{\theta_\text{cdf}}}

\newcommand\nclose[0]{\ensuremath{n_{x,\text{close}}}}
\newcommand\lclose[0]{\ensuremath{\ell_{\text{close}}}}

\newcommand\nlast[0]{\ensuremath{n_{x,\text{last}}}}
\newcommand\llast[0]{\ensuremath{\ell_{\text{last}}}}

\newcommand\ndata[0]{\ensuremath{n_{\text{data}}}} 

\newcommand\npara[0]{\ensuremath{n_{\text{para}}}} 

\newcommand\xnext[0]{\ensuremath{\xvec_{\text{next}}}}

\newcommand\xsol[1]{\ensuremath{\xvec_{\text{sol}}^{(#1)} }}

\newcommand\xbest[0]{\ensuremath{\xvec_{\text{best}}}}
\newcommand\Jbest[0]{\ensuremath{J_{\text{best}}}}

\newcommand\ldataregion[0]{\ensuremath{l_{\text{data}}}}

\newcommand\rhoInc[0]{\ensuremath{\rho_{\text{inc}}}}
\newcommand\rhoDec[0]{\ensuremath{\rho_{\text{dec}}}}

\newcommand\gtrc[0]{\ensuremath{g_{\text{trc}}}}
\newcommand\ubtrc[1]{\ensuremath{\ub{g}_{\text{trc}}^{#1}}}

\newcommand\gtrsig[0]{\ensuremath{g_{\text{tr}\sigmaGP_f}}}
\newcommand\ubtrsig[1]{\ensuremath{\ub{g}_{\text{tr}\sigmaGP_f}^{#1}}}

\newcommand\ubtrsigMin[0]{\ensuremath{\ub{g}_{\text{tr}\sigmaGP_f, \min}}}
\newcommand\ubtrsigMax[0]{\ensuremath{\ub{g}_{\text{tr}\sigmaGP_f, \max}}}

\newcommand\qEI[0]{\ensuremath{q_{\text{EI}}}}

\newcommand\fbest[0]{\ensuremath{f_{\text{best}}}}

\newcommand\gradfNoisy[0]{\ensuremath{\widetilde{\nabla f}}}

\usepackage{soul}

\definecolor{verylightblue}{rgb}{0.7,0.7,0.85}

\begin{document}

\title{Efficient Gradient-Enhanced Bayesian Optimizer with Comparisons to Conjugate-Gradient and Quasi-Newton Optimizers for Unconstrained Local Optimization}

\author[1,2]{Andr\'e L.\ Marchildon}

\author[1]{David W.\ Zingg}

\authormark{Marchildon and Zingg}
\titlemark{An Efficient Local Optimization Framework for Gradient-Enhanced Bayesian Optimizers}

\address[1]{\orgdiv{Institute for Aerospace Studies}, \orgname{University of Toronto}, \orgaddress{\state{Ontario}, \country{Canada}}}

\address[2]{This work was completed prior to the author starting his employment at Amazon}

\corres{Corresponding author Andr\'e L.\ Marchildon: \email{andre.marchildon@alumni.utoronto.ca}}


\fundingInfo{Natural Sciences and Engineering Research Council of Canada and the Ontario Graduate Scholarship Program.}

\abstract[Abstract]{The probabilistic surrogates used by Bayesian optimizers make them popular methods when function evaluations are noisy or expensive to evaluate. While Bayesian optimizers are traditionally used for global optimization, their benefits are also valuable for local optimization. In this paper, a framework for gradient-enhanced unconstrained local Bayesian optimization is presented. It involves selecting a subset of the evaluation points to construct the surrogate and using a probabilistic trust region for the minimization of the acquisition function. The Bayesian optimizer is compared to conjugate-gradient and quasi-Newton optimizers from MATLAB and SciPy for unimodal problems with 2 to 40 dimensions. The Bayesian optimizer converges the optimality as deeply as the optimizers used for comparison and often does so using significantly fewer function evaluations. For the minimization of the 40-dimensional Rosenbrock function for example, the Bayesian optimizer requires half as many function evaluations as the MATLAB and SciPy optimizers to reduce the optimality by 10 orders of magnitude. For test cases with noisy gradients, the probabilistic surrogate of the Bayesian optimizer enables it to converge the optimality several additional orders of magnitude relative to the conjugate-gradient and quasi-Newton optimizers. The final test case involves the chaotic Lorenz 63 model and inaccurate gradients. For this problem, the Bayesian optimizer achieves a lower final objective evaluation than the SciPy quasi-Newton optimizer for all initial starting solutions. The results demonstrate that a Bayesian optimizer can be competitive with quasi-Newton and conjugate-gradient optimizers when accurate gradients are available, and significantly outperforms them when the gradients are innacurate\footnote{This work comes in large part from the following thesis: Marchildon, A. L. Aug. 2024. “The Development of a Versatile and Efficient Gradient-Enhanced Bayesian Optimizer for Nonlinearly Constrained Optimization with Application to Aerodynamic Shape Optimization”. PhD Thesis. Toronto, Canada: University of Toronto. URL: https://utoronto.scholaris.ca/items/2037ae09-d3f7-4781-a8b6-71c2743d0ef8.}.}


 
\maketitle

\renewcommand\thefootnote{}

\renewcommand\thefootnote{\fnsymbol{footnote}}
\setcounter{footnote}{1}


\section{Introduction} \label{Sec_Intro}

Numerical optimization is used extensively for various applications to select parameters that provide more efficient solutions and designs \cite{nocedal_numerical_2006}. For example, the shape of an aircraft can be parameterized and then a numerical optimizer can be used to find an aircraft shape that maximizes fuel efficiency subject to various constraints, including a constraint associated with the numerical solution of the flow over the aircraft \cite{jameson_optimum_1998}. Various algorithms have been developed that are well suited to a specific type of numerical optimization, such as local or global, and gradient-free or gradient-enhanced optimization \cite{nocedal_numerical_2006}. One effective algorithm is Bayesian optimization, which is typically used without gradients and for global optimization \cite{rasmussen_gaussian_2006}. Bayesian optimizers are popular since they can quantify the uncertainty in noisy data and they use function evaluations effectively \cite{shahriari_taking_2016}, which is important when these are computationally expensive. In this paper, a framework is presented that enables Bayesian optimizers to be efficiently applied to local optimization problems.

Bayesian optimizers require two components: a probabilistic surrogate and an acquisition function \cite{shahriari_taking_2016}. The probabilistic surrogate, which is commonly a Gaussian process (GP), approximates the function of interest. Since the GP is probabilistic, the uncertainty of its posterior can be quantified \cite{rasmussen_gaussian_2006}. Another benefit of using a GP is that using a surrogate to approximate the function of interest enables the minimum of a function of interest to be found with fewer function evaluations since their locations in the parameter space can be selected effectively. This is particularly useful when function evaluations are computationally expensive, such as when flow simulations are involved \cite{paul-dubois-taine_sensitivity-based_2013}, since forming a surrogate can itself be computationally expensive. The second ingredient that Bayesian optimizers require are acquisition functions, which are formed from the surrogate and are minimized to select the next point in the parameter space where the function of interest will be evaluated \cite{brochu_tutorial_2010,zhan_expected_2020}. The acquisition function is constructed using the mean and variance of the probabilistic surrogate to balance exploration and exploitation.

Many efficient local optimization algorithms, such as conjugate-gradient and quasi-Newton methods, utilize gradients \cite{davidon_variable_1991}. Gradients are particularly useful when the parameter space is high-dimensional. While Bayesian optimizers have typically been used without gradients, they can leverage gradients when they are available \cite{morris_bayesian_1993,wu_exploiting_2018,cheng_gradient-enhanced_2023,marchildon_gradient-enhanced_2024}. Some of the challenges in doing so include a higher computational cost to train and evaluate the probabilistic surrogate, and the ill-conditioning of the covariance matrix \cite{dalbey_efficient_2013,marchildon_non-intrusive_2023}. Some of the methods that have been used to address this ill-conditioning include constraining how close function evaluations are to each other in the parameter space \cite{osborne_gaussian_2009}, or limiting the number of function evaluations that are included in the covariance matrix \cite{march_gradient-based_2011,dalbey_efficient_2013,shende_systematic_2022}. Unfortunately, these methods do not ensure that the condition number of the covariance matrix remains bounded. Furthermore, some of these methods, such as constraining how close the evaluation points of function evaluations can get to each other in the parameter space, are impractical when performing local optimization where these evaluation points naturally get closer to each other as the optimizer converges to a minimum. To address this, a previously developed preconditioning and regularization method is used that guarantees that the condition number of the gradient-enhanced covariance matrix stays below a user-set threshold \cite{marchildon_solution_2024}.

Application of Bayesian optimizers to local optimization problems has been limited since they are typically used for global optimization \cite{shahriari_taking_2016}. Mortished $\etal$ applied a gradient-enhanced Bayesian optimizer to a local optimization problem and found it to be more computationally efficient than using a genetic algorithm \cite{mortished_aircraft_2016}. March $\etal$ used a gradient-enhanced Bayesian optimizer to perform local optimization of structural and aerodynamic problems \cite{march_gradient-based_2011}. In both cases, the local Bayesian optimization frameworks that were used encountered ill-conditioned gradient-enhanced covariance matrices.

This paper presents a framework that enables gradient-enhanced Bayesian optimizer to be competitive with conjugate-gradient and quasi-Newton optimizers for local unconstrained optimization. The optimizers are compared based on the total number of function evaluations that are needed to reduce the optimality, \ie the $L_2$ norm of the gradient, by 10 orders of magnitude. This metric is used to identify which optimizer is more efficient for problems where the computational cost of the function evaluations is significantly higher than that of the computational cost of the optimization algorithm. 

The objective of this paper is to present a versatile optimization framework that expands the effective use of Bayesian optimizers from global to local optimization problems. Specifically, the goal of this paper is to demonstrate that Bayesian optimizers require no more function evaluations than conjugate-gradient and quasi-Newton optimizers to achieve deep convergence for local optimization problems with accurate gradients. Furthermore, for problems with inaccurate gradients, the goal is for the Bayesian optimizer to more consistently achieve deeper convergence with fewer function evaluations relative to conjugate-gradient and quasi-Newton optimizers. The number of function evaluations is used as the metric to compare optimizers since this is what drives the overall computational cost for many practical engineering problems. For example, aerodynamic shape optimization requires expensive computational fluid dynamics simulations that can take several minutes or even hours to perform, which is several orders of magnitude longer than the time required for an optimizer to select the next point in the design space to evaluate.

Many optimization problems of interest contain nonlinear constraints \cite{nocedal_numerical_2006,martins_engineering_2021}. While this paper is focused on local unconstrained optimization, its framework can be generalized for application to nonlinearly constrained local optimization problems \cite{marchildon_framework_2025}. Furthermore, the framework from this paper has also been used to perform nonlinearly constrained aerodynamic shape optimization \cite{marchildon_gradient-enhanced_2024}. 

The notation used in this paper is introduced in \Sec{Sec_Notation}. An overview of Gaussian processes and acquisition functions is provided in \Secs{Sec_Gp}{Sec_Bo_Acq}, respectively. The various components of the optimization framework are detailed in \Sec{Sec_LocalOptz_Framework}. Test cases are presented in \Sec{Sec_LocalOptz_TestCases} and these are used in \Sec{Sec_LocalOptz_Studies} to select parameters for the optimization framework. The Bayesian optimizer is then compared to conjugate-gradient and quasi-Newton optimizers in \Secs{Sec_LocalOptz_UnconQN}{Sec_LocalOptz_NoisyGrad} with accurate and noisy gradients, respectively. In \Sec{Sec_LocalOptz_LorenzOptz} the optimizers are compared for the optimization of a chaotic system. Finally, the conclusions for the paper can be found in \Sec{Sec_LocalOptz_Summary}.

\section{Notation} \label{Sec_Notation}

Non-bold lowercase Greek and sans-serif Latin letters are used for scalars. Vectors are denoted by bold lowercase symbols. Matrices are denoted by uppercase Greek letters or sans-serif Latin letters, \eg $\Sigma$ and $\X$. The vectors $\xvec_{i:}$ and $\xvec_{:j}$ denote the $i$-th row and $j$-th column of $\X$, respectively, and $x_{ij}$ is the entry at the $i-$th row and $j$-th column of $\X$. The vectors $\zero$ and $\one$ have all of their entries equal to zero and one, respectively, and $\Idty$ is the identity matrix. Integer quantities are denoted with the letter $n$ along with a subscript, \eg $n_x$ and $n_d$ denote the number of evaluation points and the number of dimensions, respectively. To get the diagonal of a matrix, or to form a diagonal matrix from a vector, we use $\avec = \mydiag(\A)$ and $\A = \mydiag(\avec)$, respectively.

Elementwise lower and upper bounds on variables such as the vector $\xvec$ are denoted by $\lb{\xvec}$ and $\ub{\xvec}$, respectively. Hyperparameters that approximate unknown parameters are denoted with the hat accent $\hat{\cdot}$. The symbol $\nabla$ before a scalar denotes its gradient, \eg $\nabla f$ is the gradient of the function $f$, while a $\nabla$ in the subscript is used to indicate both function and gradient evaluations are contained, \eg $f_{\nabla} = [f, \p{f}{x_1}, \ldots, \p{f}{x_{n_d}} ]$.


\section{Gaussian processes} \label{Sec_Gp}

\subsection{Mean and covariance function} \label{Sec_Gp_MeanCovFun}

A GP requires a mean and a covariance function to be fully defined. The mean function $m(\xvec)$ is often simply a constant, \ie $m(\xvec) = \beta$. This constant is set as a hyperparameter by maximizing the marginal log likelihood \cite{zhang_exploiting_2005,svensson_marginalizing_2015,ollar_gradient_2017,chen_exploiting_2022}, which is introduced in \Sec{Sec_Bo_Gp_Likelihood}. Various kernels can be used for the covariance function, such as the Gaussian, Mat\'ern $\frac{5}{2}$, and rational quadratic kernels:
\begin{alignat}{3}
	\kG(\xvec, \yvec; \gammavec) 
	&= \kG(\rvecdot)
	&&= e^{-\frac{1}{2} \| \rvecdot \|^2} \label{Eq_kern_Gaussian} \\
	\kM(\xvec, \yvec; \gammavec) 
	&= \kM(\rvecdot) 
	&&= \left(1 + \sqrt{3} \| \rvecdot \| + \| \rvecdot \|^2 \right) e^{- \sqrt{3} \| \rvecdot \|} \label{Eq_kern_Mat5f2}\\
	\kRq(\xvec, \yvec; \gammavec, \alpha) 
	&= \kRq(\rvecdot; \alpha) 
	&&= \left(1 + \frac{ \| \rvecdot \|^2 }{2 \alpha} \right)^{-\alpha}, \label{Eq_kern_RatQd}
\end{alignat}
where $\gammavec \in \mathbb{R}_+^{n_d}$, $\rdot_i = \gamma_i (x_i - y_i)$, and $\alpha >0$ is a hyperparameter for the rational quadratic kernel. The parameter $\rvecdot$ denotes a nondimensional radius. These three kernels are stationary since they depend only on $\rvecdot$, \ie the relative location of two evaluation points from one another, rather than their specific locations in the parameter space. These three kernels are all at least twice continuously differentiable and thus they can be used for gradient-free and gradient-enhanced GPs \cite{chen_optimization_2020}. Also, the limit of the rational quadratic kernel as $\alpha \rightarrow \infty$ is the Gaussian kernel \cite{rasmussen_gaussian_2006}. The function $k(\cdot,\cdot)$ is used to denote an arbitrary kernel, which is not limited to those presented in this section.

\subsection{Gradient-enhanced Gaussian processes} \label{Sec_Gp_GradEnhanced}

When gradient evaluations of the function of interest are available, a gradient-enhanced GP can be used, which is generally more accurate than a gradient-free GP \cite{morris_bayesian_1993,han_improving_2013,wu_bayesian_2017,chen_optimization_2020}. The function and derivative evaluations, which are evaluated at the rows of the matrix $\X \in \mathbb{R}^{n_x \times n_d}$, can be noisy:
\begin{alignat}{3}
	\noisy{f}_i 
		&= f(\xvec_{i:}) + \left( \epsilon_f \right)_i \quad 
		&& \forall \, i \in \{1, \ldots, n_x \} \label{Eq_noise_model_fun} \\
	\left(\noisy{\p{f}{x_d}} \right)_{\xvec_{i:}} 
		&= \p{f}{x_d} + \left( \epsilon_{\nabla f} \right)_{i,d} \quad 
		&& \forall \, i \in \{1, \ldots, n_x\}, d \in \{1, \ldots, n_d \}, \label{Eq_noise_model_grad}
\end{alignat}
where the noise is assumed to be additive, Gaussian, independent and identically distributed (IID) for each function and derivative evaluation, and to have variance $\sigma_f^2$ and $ \sigma_{\nabla f}^2$ for the function and gradient evaluations, respectively:
\begin{alignat}{3}
	\left( \epsilon_f \right)_i
		&\sim \mathcal{N}(0, \sigma_f^2) \quad
		&& \forall \, i \in \{1, \ldots, n_x \} \\
	\left( \epsilon_{\nabla f} \right)_{i,d}
		&\sim \mathcal{N} \left(0, \sigma_{\nabla f}^2 \right) \quad 
		&& \forall \, i \in \{1, \ldots, n_x\}, d \in \{1, \ldots, n_d \}, 
\end{alignat}
which is a homoscedastic noise model since the variance is the same for all derivative evaluations. The priors for the noisy function and derivative evaluations are
\begin{align}
	\cov \left(\noisy{\fvec}(\X) \right) 
		&= \Sigma(\X; \sigK, \gammavec, \sigma_f) \label{Eq_Sigma_common} \\
		&= \sigK^2 \K(\X; \gammavec) + \sigma_f^2 \Idty \\
	\cov \left(\noisy{\p{\fvec}{x_d}} \right)_{\X}
		&= \sigK^2 \left( \p{^2\K}{x_d^2} \right)_{\X} + \sigma_{\nabla f}^2 \Idty \quad \forall \, d \in \{1, \ldots, n_d\}, \label{Eq_prior_noise_model_grad}
\end{align}
where $\Sigma$ is the gradient-free covariance matrix and $\K$ is the gradient-free $n_x \times n_x$ kernel matrix:
\begin{equation} \label{Eq_K}
	\K(\X; \gammavec) =
	\begin{bmatrix}
		k(\xvec_{1:}, \xvec_{1:}; \gammavec) 	& k(\xvec_{1:}, \xvec_{2:}; \gammavec) 		& \ldots 	& k(\xvec_{1:}, \xvec_{n_x:}; \gammavec) \\
		k(\xvec_{2:}, \xvec_{1:}; \gammavec) 	& k(\xvec_{2:}, \xvec_{2:}; \gammavec) 		& \ldots 	& k(\xvec_{2:}, \xvec_{n_x:}; \gammavec) \\
		\vdots 	&	\vdots						& \ddots 	& \vdots \\
		k(\xvec_{n_x:}, \xvec_{1:}; \gammavec) 	& k(\xvec_{n_x:}, \xvec_{2:}; \gammavec) 	& \ldots 	& k(\xvec_{n_x:}, \xvec_{n_x:}; \gammavec)
	\end{bmatrix}.
\end{equation}
The notation $\left( \p{^2\K}{x_d^2} \right)_{\X}$ from \Eq{Eq_prior_noise_model_grad} indicates that the entries of $\K$ from \Eq{Eq_K} are each differentiated twice with respect to $x_d$ and then evaluated with the input $\X$. The joint distribution for noisy function and gradient evaluations is given by
\begin{equation}
	\begin{bmatrix}
		\noisy{\fvec}(\X) \\
		\left(\noisy{\p{\fvec}{x_1}} \right)_{\X} \\
		\vdots \\
		\left(\noisy{\p{\fvec}{x_{n_d}}} \right)_{\X} 
	\end{bmatrix}
	\sim \mathcal{N} \left(
	\begin{bmatrix}
		\one_{n_x} \beta \\
		\zero_{n_x n_d}
	\end{bmatrix},
	\Sigmag
	\right),
\end{equation}
where $\Sigmag$ is the gradient-enhanced covariance matrix. 

The gradient-enhanced covariance matrix can be significantly ill-conditioned, as will be discussed further in \Sec{Sec_LocalOptz_Framework_Precon}. Consequently, it is common to include a nugget $\etaKg >0$. The gradient-enhanced covariance matrix with the addition of a nugget is given by:
\begin{equation} \label{Eq_Sigmag}
	\Sigmag(\X; \sigK, \gammavec, \etaKg, \W, \hpstdfval, \hpstdfgrad) 
	= \sigK^2 \left(\Kg(\X;\gammavec) + \etaKg \W \right) + \Vg(\hpstdfval, \hpstdfgrad),
\end{equation}
where $\W$ is a diagonal matrix with nonnegative entries, $\Kg$ is the gradient-enhanced kernel matrix, and $\Vg$ is given by
\begin{equation} \label{Eq_Vg}
	\Vg(\hpstdfval, \hpstdfgrad) 
	= \mydiag \left( \hpstdfval^2 \one_{n_x}^\top, 
	\hpstdfgrad^2 \one_{n_x n_d}^\top \right),
\end{equation}
where $\hpstdfval$ and $\hpstdfgrad$ are hyperparameters that estimate the true noise variance for the function and gradient evaluations, respectively. The addition of a nugget is mathematically equivalent to having noisy data. While $\W$ could be any diagonal matrix with nonnegative entries, we are specifically interested in diagonal matrices of the following form:
\begin{equation}
	\mydiag(\W) = [1, w_1, \ldots, w_{n_d}] \otimes \Idty_{n_x},
\end{equation}
where $w_d \geq 0 \, \forall \, d \in \{1, \ldots, n_d\}$. The prior for the noisy function and derivative evaluations along with the nugget that is added is:
\begin{align}
	\cov \left(\noisy{\fvec}(\X) \right) 
		&= \sigK^2 \K(\X; \gammavec) + \left( \sigK^2 \etaKg + \sigma_f^2 \right) \Idty \\
	\cov \left(\noisy{\p{\fvec}{x_d}} \right)_{\X}
		&= \sigK^2 \left( \p{^2\K}{x_d^2} \right)_{\X} 
		+\left( \sigK^2 \etaKg w_d + \hpstdfgrad^2 \right) \Idty \quad \forall \, d \in \{1, \ldots, n_d\}, \label{Eq_noise_model_grad_w_nugget}
\end{align}
where the terms on the right-hand side with the nugget $\etaKg$ are usually small since the nugget is on the order of $10^{-9}$ when $\condmax = 10^{10}$ \cite{marchildon_solution_2024}. With this structure, the inclusion of the nugget $\etaKg$ results in homoscedastic noise on the evaluation of the function. For the gradient evaluations, the noise is homoscedastic when the derivative with respect to each variable $x_i \, \forall \, i \{1, \ldots, n_x \}$ is considered individually. However, if all of the derivatives are considered together, and $\W$ is not a multiple of the identity matrix, then the noise is heteroscedastic since its variance is not the same for all of the derivatives.

The gradient-enhanced kernel matrix can be constructed with either the indirect or the direct method \cite{zimmermann_maximum_2013}. The former uses a kernel matrix of the same form as the gradient-free kernel matrix and adds additional evaluation points to approximate the gradients. This method is well suited if the gradients are calculated with finite differences. In contrast, the direct method is generally used if the gradients are calculated analytically \cite{han_improving_2013,dalbey_efficient_2013,wu_exploiting_2018,laurent_overview_2019}. 
The gradients in this paper are calculated analytically and therefore, the direct method is used. The gradient-enhanced kernel matrix $\Kg$ from \Eq{Eq_Sigmag}, which uses the direct method, is given by:
\begin{equation} \label{Eq_Kg} 
	\Kg(\X; \gammavec) = 
	\begin{bmatrix}
		\K \left( \X \right) 				& \left(\p{\K}{y_1} \right)_\X			& \ldots 	& \left( \p{\K}{y_{n_d}} \right)_\X \\
		 \left( \p{\K}{x_1} \right)_\X 	& \left( \pp{\K}{x_1}{y_1} \right)_\X 	& \ldots 	& \left( \pp{\K}{x_1}{y_{n_d}} \right)_\X \\ 
		\vdots 			& 	\vdots				& \ddots 	& \vdots \\
		 \left( \p{\K}{x_{n_d}} \right)_\X 	& \left( \pp{\K}{x_{n_d}}{y_1} \right)_\X 	& \ldots 	& \left( \pp{\K}{x_{n_d}}{y_{n_d}} \right)_\X
	\end{bmatrix}.
\end{equation} 

The joint distribution of the gradient-enhanced GP for the function and gradient evaluations at the rows $\X$ and the function evaluation at a point $\xvec'$ is given by
%
\begin{equation}
	\begin{bmatrix}
		\noisy{\fgrad} \\
		f(\xvec')
	\end{bmatrix}
	\sim \mathcal{N} \left(
	\begin{bmatrix}
		\mgrad(\X) \\
		\mvec(\xvec')
	\end{bmatrix},
	\begin{bmatrix}
		\Sigmag & \sigK^2 \kvecgrad(\X; \xvec') \\
		\sigK^2 \kvecgrad(\xvec', \X) & \sigK^2 k(\xvec', \xvec')
	\end{bmatrix}
	\right),
\end{equation}
where $\left( \kvecgrad(\xvec', \X) \right)^\top = \kvecgrad(\X, \xvec')$ and
\begin{align} 
	\mgrad(\X) =
	\begin{bmatrix}
		\mvec(\X) \\
		\p{\mvec(\X)}{x_1} \\
		\vdots \\
		\p{\mvec(\X)}{x_{n_d}}
	\end{bmatrix}, \quad
	\kvecgrad(\X; \xvec') =
	\begin{bmatrix}
		\kvec(\X, \xvec') \\
		\left(\p{\kvec}{x_1} \right)_{(\X, \xvec')} \\
		\vdots \\
		\left(\p{\kvec}{x_{n_d}} \right)_{(\X, \xvec')}
	\end{bmatrix}, \quad
	\noisy{\fgrad}(\X) = 
	\begin{bmatrix}
		\noisy{\fvec}(\X) \\
		\left(\noisy{\p{\fvec}{x_1}} \right)_{\X} \\
		\vdots \\
		\left(\noisy{\p{\fvec}{x_{n_d}}} \right)_{\X}
	\end{bmatrix}. \label{Eq_def_kvecgrad_fgrad}
\end{align}

The mean and variance of the posterior for the gradient-enhanced GP are formed by conditioning the prior of the GP on the observations, \ie the function and gradient evaluations $\noisy{\fgrad}(\X)$:
\begin{align}
	\muGP_f(\xvec') 
	&= m(\xvec') + \sigK^2 \, \kvecgrad(\xvec', \X) \Sigmag^{-1} \left(\noisy{\fgrad}(\X) - \mgrad(\X) \right) \label{Eq_muGp_w_grad} \\
	\sigmaGP_f^2(\xvec') 
	&= \sigK^2 \left( k(\xvec', \xvec') - \sigK^{-2} \, \kvecgrad(\xvec', \X) \Sigmag^{-1} \kvecgrad(\X, \xvec') \right). \label{Eq_sigmaGp_w_grad} 
\end{align}
The matrix $\Sigmag$ is symmetric positive definite for $\etaKg >0$ and positive diagonal values for $\W$. Therefore, the Cholesky decomposition of $\Sigmag$ can be calculated once with a cost that scales as $\order{n_d^3 (n_x + 1)^3}$. This enables $\muGP_f(\xvec')$ and $\sigmaGP_f(\xvec')$ to be evaluated at a cost that scales as $\order{n_d^2 (n_x + 1)^2}$ for different points $\xvec'$ in the parameter space. The structure of the gradient-enhanced covariance matrix $\Sigmag$ can be leveraged to evaluate $\muGP_f$ and $\sigmaGP_f$ efficiently without needing to first calculate the Cholesky decomposition \cite{de_roos_high-dimensional_2021}. However, using the Cholesky decomposition is efficient if $\muGP_f$ and $\sigmaGP_f$ are evaluated at several points in the parameter space.

\subsection{Marginal log-likelihood} \label{Sec_Bo_Gp_Likelihood}

A GP has several hyperparameters, including $\beta$, $\sigK$, $\gammavec$, and $\alpha$ if the rational quadratic kernel from \Eq{Eq_kern_RatQd} is used. There are also the hyperparameters $\hpstdfval$ and $\hpstdfgrad$ to consider if the function and gradient evaluations are noisy, respectively. The most popular method of selecting values for these hyperparameters is by maximizing the marginal likelihood function \cite{toal_kriging_2008,toal_development_2011,snoek_practical_2012}, which is given by
\begin{equation} \label{Eq_marginal_lkd_w_integral}
	p( \noisy{\fvec} | \X ) = \int p( \noisy{\fvec} | \fvec, \X ) p \left( \fvec | \X \right) d \fvec,
\end{equation}
where, for example, $p(a|b,c)$ denotes the probability of $a$ given $b$ and $c$, the integral is the marginalization over all possible function evaluations $\fvec$ that could have been observed at the rows of $\X$, $p( \noisy{\fvec} | \fvec, \X )$ is the likelihood, and $p \left( \fvec | \X \right)$ is the prior. The likelihood for a GP model is a factorized Gaussian, $\ie$ $\noisy{\fvec} | \fvec, \X \sim \mathcal{N}(\fvec, \sigma_f^2 \Idty)$, and its prior is Gaussian: $\fvec | \X \sim \mathcal{N}(\zero, \Sigma)$ \cite{rasmussen_gaussian_2006}. The integral in \Eq{Eq_marginal_lkd_w_integral} thus involves the product of two Gaussian distributions and it evaluates to
\begin{equation} \label{Eq_lkd}
	L(\gammavec, \beta, \sigK, \hpstdfval, \hpstdfgrad; \X, \fgrad, \etaKg, \W) 
	= \frac{e^{-\frac{ \left(\fgrad - \onemod \beta \right)^\top \Sigmag^{-1} \left( \fgrad - \onemod \beta \right) }{2}} }{ \left(2 \pi \right)^{\frac{n_x(n_d+1)}{2}} \sqrt{\det \left(\Sigmag \right) }},
\end{equation}
where $\onemod = [\one_{n_x}^\top, \zero_{n_x n_d}^\top]^\top$. The hyperparameters $\gammavec$, $\beta$, $\sigK$, $\hpstdfval$, and $\hpstdfgrad$ can have a significant impact on the mean and variance of the posterior for the GP. The marginal likelihood is either maximized to select the hyperparameters or marginalized by integrating over all potential values of the hyperparameters. This latter method cannot be done analytically and thus requires the use of a Monte Carlo method, for example \cite{rasmussen_gaussian_2006}. The marginal likelihood has often been found to be orders or magnitude larger at its global maximum than at other local maxima \cite{rasmussen_gaussian_2006}. In such cases, the hyperparameters that are selected by maximizing the marginal likelihood provide a GP with a posterior that is a good approximation to the function of interest and it is thus the method that is used in this paper.


To simplify the calculations, the log of the marginal likelihood function is usually used and constant terms are dropped:
\begin{equation} \label{Eq_ln_lkd_noisy}
	\ln(L) = -\frac{1}{2} \ln \left( \det \left( \Sigmag \right) \right) - \frac{1}{2} \left( \fgrad(\X) - \onemod \beta \right)^\top \Sigmag^{-1} \left( \fgrad(\X) - \onemod \beta \right).
\end{equation}
Since $\Sigmag$ is symmetric positive definite for $\etaKg >0$ and nonzero diagonal entries of $\W$, a Cholesky decomposition $\Lmat \Lmat^\top = \Sigmag$ can be used. The logarithm of the determinant of $\Sigmag$ can then be efficiently calculated with $\ln \left( \det \left( \Sigmag \right) \right) = 2 \sum_i \ln(\Lmat_{ii})$. There is a closed-form solution for the value of $\beta$ that maximizes \Eq{Eq_ln_lkd_noisy}, which can be found in Appendix \ref{Apdx_NoiseFreeLkd}. The same appendix also provides the closed form solution for the value of $\sigK$ that maximizes \Eq{Eq_ln_lkd_noisy} when the function and gradient evaluations are both noise-free.

The optimization of hyperparameters is the most expensive step in using GPs and its cost grows as the dimensionality of the problem increases since there are more variables. Several methods have been developed to reduce the number of hyperparameters in order to reduce this computational cost. Some of the methods include partial least squares \cite{bouhlel_improved_2016,amine_bouhlel_efficient_2018}, maximum information coefficient \cite{zhao_efficient_2020}, active subspaces \cite{chen_exploiting_2022,cheng_gradient-enhanced_2023}, and sliced gradient-enhanced Kriging \cite{cheng_sliced_2023}. Another method that has been used to reduce the computational cost is to maximize the marginal log likelihood by constraining all of the entries in $\gammavec$ to be equal \cite{chung_using_2002}. This last method significantly reduces the computational cost but it also limits the flexibility of the kernel, making it simply a radial basis function.

\section{Acquisition functions} \label{Sec_Bo_Acq}

The next point in the design space where the function of interest and its gradient are evaluated is denoted by $\xnext$ and it is selected by solving
\begin{equation} \label{Eq_acq_uncon}
	\xnext = \argmin_{\xvec} q(\xvec),
\end{equation}
where $q(\cdot)$ is a user-selected acquisition function. Acquisition functions generally depend on both the mean $\muGP_f$ and the variance $\sigmaGP_f^2$ of the posterior of the GP that is approximating the function of interest. The acquisition function balances exploration and exploitation in order to identify the next point in the parameter space where the function of interest and its gradient should be evaluated. One such example is the upper confidence acquisition function, which is given by
\begin{equation} \label{Eq_acq_UC}
	q_{\text{UC}}(\xvec; \omega) = \muGP_f(\xvec) - \omega \sigmaGP_f(\xvec),
\end{equation}
where $\omega \geq 0$ is a parameter that is set by the user. Large values of $\omega$ promote exploration, which is useful for avoiding getting stuck in a local minimum, while a small value or a value of zero promotes exploitation. Another popular acquisition function for global optimization is expected improvement, which was introduced by Jones \etal \cite{jones_efficient_1998}:
\begin{align*}
	q_{\text{EI}}(\xvec; \fbest) 
	&= \int_{-\infty}^{\fbest} (\fbest - f) \pdf \left(\frac{f - \muGP_f(\xvec)}{\sigmaGP_f(\xvec)} \right) df \\
	&= \left( \fbest - \muGP_f(\xvec) \right) \cdf \left(\frac{\fbest - \muGP_f(\xvec)}{\sigmaGP_f(\xvec)} \right) + \sigmaGP_f(\xvec) \pdf \left(\frac{\fbest - \muGP_f(\xvec)}{\sigmaGP_f(\xvec)} \right), \yesnumber \label{Eq_acq_EI}
\end{align*}
where $\fbest$ is the minimum evaluation of the function of interest, while $\pdf(\cdot)$ and $\cdf(\cdot)$ are the standard normal probability and cumulative density functions, respectively. The expected improvement acquisition function is extensively used and probably the most popular acquisition function \cite{zhan_expected_2020}. When the function evaluations are noisy, then $\fbest$ can be replaced in \Eq{Eq_acq_EI} by the lowest evaluation of the mean of the posterior $\muGP_f$ at rows of $\X$ \cite{bernardo_optimization_2011}.


\section{Bayesian optimization framework for local minimization} \label{Sec_LocalOptz_Framework}

\subsection{Solution to the ill-conditioning of $\Sigmag$} \label{Sec_LocalOptz_Framework_Precon}

The evaluation of the mean and variance of the posterior of the GP with \Eqs{Eq_muGp_w_grad}{Eq_sigmaGp_w_grad}, respectively, or of the marginal log-likelihood from \Eq{Eq_ln_lkd_noisy} involves taking the product of the inverse of the gradient-enhanced covariance matrix $\Sigmag$ with different vectors. The covariance matrix is symmetric positive definite when a positive nugget is used and thus a Cholesky decomposition can be used. However, since $\Sigmag$ can be significantly ill-conditioned, the Cholesky decomposition may fail. To address this the preconditioning method from Marchildon and Zingg \cite{marchildon_solution_2024} is used. A preconditioned gradient-enhanced covariance matrix is formed with
\begin{equation}
	\Sigmagdot = \sigK^2 \left( \Kgdot + \etaKgdot \Idty \right),
\end{equation}
where 
\begin{align} 
	\Kgdot 
		&= \Pinv \Sigmag(\etaKg = 0) \Pinv \\
	\etaKgdot 
		&= \frac{ \max_{i} \sum_{j=1}^{n_d (n_d + 1)} \left| \Kgdot \right|_{ij} }{\condmax - 1}, \label{Eq_etaKgdot}
\end{align}
and the preconditioning matrix is
\begin{equation}
	\Pmat = \mydiag \left( \mydiag \sqrt{ \sigK^{-2} \mydiag \left( \Sigmag(\etaKg = 0) \right) } \right).
\end{equation}
The Cholesky decomposition $\dot{\Lmat} \dot{\Lmat}^\top = \Sigmagdot$ is performed since the condition number of $\Sigmagdot$ is bounded from above by $\condmax$, which is selected by the user, for all combinations of evaluation points (even if they are collocated) and all positive values of the hyperparameters, \ie $\kappa(\Sigmagdot(\X;\sigK, \gammavec, \etaKgdot, \W, \hpstdfval, \hpstdfgrad) \Pinv) \leq \condmax \, \forall \, \X \in \mathbb{R}^{n_x \times n_d}, \gammavec \in \mathbb{R}_+^{n_d}, \sigK, \hpstdfval, \hpstdfgrad \geq 0$ \cite{marchildon_solution_2024}. 
The Cholesky decomposition of $\Sigmag$ is then recovered with $\Lmat = \Pmat \dot{\Lmat}$, where $\Lmat \Lmat^\top = \Sigmag$ with $\W = \Pmat \Pmat$ and $\etaKg$ coming from \Eq{Eq_etaKgdot}. The nugget value $\etaKgdot$ from \Eq{Eq_etaKgdot} scales as $\order{n_x \sqrt{n_d}}$ when the Gaussian kernel from \Eq{Eq_kern_Gaussian} is used \cite{marchildon_solution_2024}.

\subsection{Data region} \label{Sec_LocalOptz_Framework_DataRegion}

All of the function and gradient evaluations can be used when the marginal log likelihood is maximized to select the hyperparameters and to evaluate the posterior of the GP. However, as will be demonstrated in the numerical results in \Sec{Sec_LocalOptz_Studies_nclose}, it is advantageous for local minimization to use only the function and gradient evaluations near $\xbest$. This is the evaluation point with the lowest evaluation of the merit function, \ie the objective function for unconstrained problems. If all of the evaluation points are used, then the hyperparameters that maximize the marginal log likelihood will provide a surrogate that is relatively accurate near all evaluation points, while being less accurate around $\xbest$ than a surrogate constructed using only the evaluation points in its vicinity. A secondary advantage of using only a subset of the evaluation points is that this reduces the computational cost to select the hyperparameters of the GP and to evaluate its posterior. Mortished $\etal$ for example used only the function and gradient evaluations in a hyperrectangle that was set slightly larger than the hyperrectangle for their rectangular trust region \cite{mortished_aircraft_2016}.

The data region is used to denote the region where all of the function and gradient evaluations are used to select the hyperparameters of the GP and to evaluate its posterior. Algorithm~\ref{Alg_DataRegion} details how the data region is selected. There are two requirements for the data region: to have at least the $\min(\nclose, n_x)$ closest evaluation points to $\xbest$, and the $\min(\nlast, n_x)$ most recent evaluation points. The default values that are used are $\nclose = 20$ and $\nlast = 3$. The evaluation points are selected for the data region regardless of whether they provide a reduction in the merit function or not. This is in contrast to quasi-Newton methods where the function and gradient evaluations are not used to update the Hessian approximation if progress is not made \cite{nocedal_numerical_2006}.

\begin{algorithm}[t!]
	\caption{Selecting the evaluation points $\Xdata$ for the data region}
	\label{Alg_DataRegion}
	\begin{algorithmic}[1]
		\Statex{\textbf{Required:} All $n_x$ previous evaluation points in the matrix $\Xall$ and $\xbest$, \ie the evaluation point with the lowest merit function evaluation.}
		\Statex{\textbf{Select:} $\nlast (3)$ and $\nclose(20)$, with default values indicated in parentheses.}
		\State{Calculate the distance to $\xbest$: $l_i = \|\xvec_i - \xbest \|_2 \, \forall \, i \in \{1, \ldots n_x \}$}
		\If{$n_x \leq \nclose$}
		\State{$\Xdata = \Xall$ and $\ldataregion = \max(\lvec)$}
		\Else
		\State{$\llast = \max(l_{m}, \ldots, l_{n_x})$, where $m = \max(n_x - n_{\text{recent}} +1, 1)$}
		\State{Identify the $\nclose$-th smallest $l_i$: $\lclose = \left(\text{sort}(\lvec) \right)_{\nclose}$}
		\State{$\ldataregion = \max(\llast, \lclose)$}
		\State{Append $x_i$ as a row for $\Xdata \, \forall i \in \{1, \ldots, n_x\}$ if $l_i \leq \ldataregion$}
		\EndIf
		\State{\textbf{Return:} Data region radius $\ldataregion$ and the matrix of evaluation points $\Xdata$ of size $\ndata \times n_d$}
	\end{algorithmic}
\end{algorithm}

\subsection{Selecting the hyperparameters} \label{Sec_LocalOptz_Framework_HpOptz}

The hyperparameters of the GP are selected by maximizing the marginal log likelihood, as presented in \Sec{Sec_Bo_Gp_Likelihood}. Only the function and gradient evaluations from the data region from Algorithm~\ref{Alg_DataRegion} are used to evaluate the marginal log likelihood. In the noise-free case, \ie $\hpstdfval = \hpstdfgrad = 0$, the only hyperparameters that are selected by numerically minimizing \Eq{Eq_ln_lkd_noisefree_final} are $\gammavec$, and $\alpha$ if the rational quadratic kernel from \Eq{Eq_kern_RatQd} is used. The hyperparameters $\beta$ and $\sigK$ are calculated with \Eqs{Eq_lkd_beta}{Eq_lkd_sigK2_noisefree}, respectively. 

In the case when there are noisy function or gradient evaluations, the hyperparameters that must be selected with numerical optimization are $\gammavec$, $\sigK$, $\hpstdfval$ if the function evaluations are noisy, $\hpstdfgrad$ if the gradient evaluations are noisy, and $\alpha$ if the rational quadratic kernel is used. The marginal log likelihood is evaluated with \Eq{Eq_ln_lkd_noisy} and the hyperparameter $\beta$ once again comes from \Eq{Eq_lkd_beta}.

The marginal log-likelihood function is often multimodal, and thus the selection of hyperparameters from a local maximum of the marginal log likelihood can result in a surrogate that poorly approximates the function of interest \cite{toal_development_2011}. However, the optimization of the marginal log likelihood is the most expensive step in using a Bayesian optimizer since it requires the covariance matrix to be constructed and its Cholesky decomposition to be calculated each time the hyperparameters are changed. The size of the gradient-enhanced covariance matrix $\Sigmag$ is $n_x (n_d+1) \times n_x (n_d+1)$ and thus the cost of the Cholesky decomposition scales as $n_x^3 (n_d+1)^3$. It is therefore important to optimize the hyperparameters efficiently to get an accurate surrogate. The SciPy SLSQP (sequential least squares programming) numerical optimizer is used \cite{virtanen_scipy_2020}, which is a quasi-Newton optimizer using the BFGS updating formula for the approximation of the Hessian \cite{broyden_convergence_1970,fletcher_new_1970,goldfarb_family_1970,shanno_conditioning_1970}. This optimizer was selected since it was found to efficiently perform the optimization using the analytical derivatives of the marginal log likelihood. Since the hyperparameter values can span several orders of magnitude, the logarithmic values of the hyperparameters are used as the variables for the maximization of the marginal log likelihood.

\begin{algorithm}[t!]
	\caption{Latin hypercube sampling for the $\npara$ hyperparameters of the GP} 
	\label{Alg_GpHyperparameters}
	\begin{algorithmic}[1] 
		\Statex{\textbf{Select:} Initial hyperparameters $\gamma_{\text{init}} (10^{-2})$, along with $\hat{\sigma}_{f,\text{init}} (10^{-5})$, $\hat{\sigma}_{\nabla f,\text{init}} (10^{-5})$, and $\hat{\sigma}_{\K,\text{init}} (1)$ if the function or gradient evaluations are noisy, also select $\nmed (5)$, $n_{\log} (3)$, and $\nLhs (50)$, with default values indicated in parentheses.}
		\Statex{\textbf{Required:} $\Xdata$ from Algorithm~\ref{Alg_DataRegion} and vectors of the hyperparameters at the $j$-th optimization iteration of the Bayesian optimizer $\phivec^{j} \, \forall \, j \in \{m, \ldots, \noptz \}$, where \newline $m = \max(1, \noptz - \nmed + 1)$.}
		\If{$\noptz = 0$}
		\State{Set $\phivec_{\text{med}}$ to the initial hyperparameters}
		\Else
		\State{$\phi_{\text{med,i}} = \med \left( \phivec_i^{(m)}, \ldots, \phivec_i^{(\noptz)} \right) \forall \, i \{1, \ldots, \npara \}$, where $m = \max(1, \noptz - \nmed +1)$}
		\Statex{\hspace{1.5em}and $\npara$ is the number of hyperparameters being selected numerically}
		\EndIf
		\State{$\ub{\phivec}_{\text{LHS}} = \log(\phivec_{\text{med}}) + n_{\log}$}
		\State{$\lb{\phivec}_{\text{LHS}} = \log(\phivec_{\text{med}}) - n_{\log}$}
		\State{Calculate the Latin hypercube sampling: $\Phi = \text{LHS}(\nLhs, \lb{\phivec}_{\text{LHS}}, \ub{\phivec}_{\text{LHS}})$}
		\State{\textbf{Return:} The matrix $\Phi$, where each row contains the logarithmic values of the hyperparameters}
	\end{algorithmic}
\end{algorithm}

Algorithm~\ref{Alg_GpHyperparameters} uses a Latin hypercube sampling to select several points in the hyperparameter space. The bounds for Latin hypercube sampling are selected to be $\pm n_{\log}$ of the median of the previous $\nmed$ values of the hyperparameters, where the default values for $n_{\log}$ and $\nmed$ are 3 and 5, respectively. The median is used since it is less sensitive than the mean to outlier values. The Latin hypercube sampling returns the logarithmic values of the hyperparameters since these are the variables of the numerical optimizer.

Two techniques were used to avoid selecting hyperparameters from a poor local maximum of the marginal log likelihood. The first technique involves using Algorithm~\ref{Alg_GpHyperparameters} with $\nLhs = 5$, which returns 5 points in the hyperparameter space. An independent local optimization using each of these $\nLhs$ starting points is performed using the SciPy SLSQP optimizer. The second method also uses Algorithm~\ref{Alg_GpHyperparameters} but with $\nLhs = 50$. The marginal log likelihood is evaluated at each of these $\nLhs$ initial hyperparameter values. However, only a single local optimization is performed starting with the hyperparameter values that provide the highest marginal log likelihood evaluation. Generally, the numerical optimizer was found to take about 50 function and gradient evaluations for each local minimization. The first method thus takes about 250 function and gradient evaluations. The second method only requires about 100 function evaluations and 50 gradient evaluations, making it significantly less expensive than the first method. Both methods were found to provide similar performance and thus, the latter method was used since it is less computationally expensive. Using 50 starting points to sample the hyperparameter space was found to be sufficient for problems with up to 40 design variables. However, for higher-dimensional problems the number of starting points will likely need to be increased.

\subsection{Trust regions} \label{Sec_LocalOptz_Framework_TrustRegion}

Trust regions are often used to constrain how far the next evaluation point is from previous evaluation points. The trust region for quasi-Newton optimizers is commonly a hypersphere around $\xbest$ \cite{nocedal_numerical_2006,martins_engineering_2021}. Bayesian optimizers are most often used for global optimization and have thus not typically used trust regions. However, some Bayesian optimizers have used trust regions, such as a hyperrectangle \cite{mortished_aircraft_2016}.

Two trust regions are used for the present Bayesian optimizer; a circular trust region and a probabilistic trust region that leverages the uncertainty quantification of the GP's posterior. The circular trust region is simply a hypersphere around $\xbest$:
\begin{align}
	\gtrc(\xvec; \xbest) 
	&= \| \xvec - \xbest \|_2^2 \label{Eq_tr_circle_val} \\
	&\leq \ubtrc{j}, \nonumber
\end{align}
where $\ubtrc{j}$ is the maximum $\ell_2$ squared distance allowed around $\xbest$ at the $j$-th optimization iteration. \Eq{Eq_tr_circle_val} involves the squared distance between $\xvec$ and $\xbest$ since this simplifies its gradient calculations and avoids dividing by zero when $\xvec = \xbest$. Algorithm~\ref{Alg_TrustRegion_circle} indicates how the upper bound $\ubtrc{j}$ is selected. The upper bound for the circular trust region $\ubtrc{j}$ is only increased if progress is made after the latest function evaluation, \ie $J^{i} < \Jbest^{i-1}$, where $i$ indicates the function evaluation iteration. If progress was not made during the last two iterations, then $\ubtrc{j}$ is decreased. Otherwise, $\ubtrc{j}$ is kept the same as the previous optimization iteration.

\begin{algorithm}[t!]
	\caption{Selecting the upper bound for the circular trust region: $\ubtrc{j}$} 
	\label{Alg_TrustRegion_circle}
	\begin{algorithmic}[1]
		\Statex{\textbf{Required:} The merit function evaluations $J^{i}$, $J^{i-1}$, and $\Jbest^{i-1}$, where $i$ is the function evaluation iteration, the last circular trust region evaluation and upper bound $\gtrc^{j-1}$ and $\ubtrc{j-1}$, respectively, where $j$ is the optimization iteration, and $\ndata$ and $\ldataregion$ come from Algorithm~\ref{Alg_DataRegion}.}
		\Statex{\textbf{Select:} $0 < \ubtrc{0} (1)$, $0 < \rhoDec (0.5) < 1 < \rhoInc (2)$, and $0 < \rho_{\text{data}} (0.9)$ with default values indicated in parentheses.}
		\If{$\ndata = 1$}
		\State $\ubtrc{j} = \ubtrc{0}$
		\Else
		\If{$J^{i} < \Jbest^{i-1}$}
		\State{$\ubtrc{j} = \max \left( \rhoInc \cdot \gtrc^{j-1}, \ubtrc{j-1} \right)$}
		\ElsIf{$J^{i-1} \leq \Jbest^{i-1}$}
		\State{$\ubtrc{j} = \ubtrc{j-1}$}
		\Else 
		\State{$\ubtrc{j} = \rhoDec \cdot \ubtrc{j-1}$}
		\EndIf
		\EndIf
		\If{$\ndata \geq 5$}
		\State{$\ubtrc{j} = \min \left(\ubtrc{j}, \rho_{\text{data}} \cdot \ldataregion) \right)$}
		\EndIf
		\State{\textbf{Return:} $\ubtrc{j}$, \ie the upper bound of the circle trust region}
	\end{algorithmic}
\end{algorithm}

The second trust region that is used leverages the probabilistic component of the GP:
\begin{align}
	\gtrsig(\xvec; \gammavec, \sigKf, \hpstdfval, \hpstdfgrad)
	&= \frac{\sigmaGP_f^2(\xvec; \gammavec, \sigKf, \hpstdfval, \hpstdfgrad)}{\sigKf^2} \label{Eq_tr_sigma_val} \\
	&\leq \ubtrsig{j}, \nonumber
\end{align}
where $\ubtrsig{j}$ is the upper-bound of the trust region at the $j$-th optimization iteration, and the $f$s in the subscripts of $\sigmaGP_f$ and $\sigKf$ indicate these are for the posterior of the GP that approximates the objective function $f(\xvec)$. This distinction is made since separate GPs can be used to approximate nonlinear constraints when they are present. For kernels such as the Gaussian, Mat\'ern, and rational quadratic kernels from \Eqss{Eq_kern_Gaussian}{Eq_kern_Mat5f2}{Eq_kern_RatQd}, respectively, which all have $k(\xvec,\xvec) = 1$, we have the following relation:
\begin{equation} \label{Eq_frac_sigf_sigK}
	\frac{\sigmaGP_f^2}{\sigKf^2} 
	= \left(1 - \sigKf^{-2} \, \kvecgrad(\X, \xvec)^\top \Sigmag^{-1} \kvecgrad(\X, \xvec) \right),
\end{equation}
which is always between 0 and 1 since $\Sigmag$ is symmetric positive definite. Therefore, the $\sigma$ trust region requires that $\ubtrsig{j} > 0 \, \forall \, j > 0$ and is only active if $\ubtrsig{j} < 1 \, \forall \, j > 0$. The parameter $\ubtrsig{j}$ is selected using Algorithm~\ref{Alg_TrustRegion_sigma}. The $\sigma$ trust region is only active when $\ndata \geq 10$ since the surrogate is not accurate if there are too few evaluation points. Just like the circular trust region from Algorithm~\ref{Alg_TrustRegion_circle}, $\ubtrsig{j}$ is increased if progress is made, decreased if progress is not made during the last two consecutive function evaluations, and kept constant otherwise. Lower and upper bounds $\ubtrsigMin$ and $\ubtrsigMax$, respectively, are used to ensure that $\ubtrsig{j}$ does not get too small or too large.

\begin{algorithm}[t!]
	\caption{Selecting the upper bounds for the $\sigma$ trust region: $\ubtrsig{j}$}
	\label{Alg_TrustRegion_sigma}
	\begin{algorithmic}[1]
		\Statex{\textbf{Required:} The number of data points $\ndata$ in the data region from Algorithm~\ref{Alg_DataRegion}, the merit function evaluations $J^{i}$, $J^{i-1}$, and $\Jbest^{i-1}$, where $i$ is the function evaluation iteration, and the previous $\sigma$ trust region evaluation and upper bound $\gtrsig^{j-1}$ and $\ubtrsig{j-1}$, respectively, where $j$ is the optimization iteration.}
		\Statex{\textbf{Select:} $0 < \ubtrsigMin (0.05^2) \leq \ubtrsig{0} (0.2^2)< \ubtrsigMax (0.4^2) \leq 1$, and $0 < \rhoDec (0.5) < 1 < \rhoInc (2)$, with default values indicated in parentheses.}
		\If{$\ndata < 10$}
		\State $\ubtrsig{j} = \infty$
		\ElsIf{$\ndata = 10$}
		\State $\ubtrsig{j} = \ubtrc{0}$
		\Else
		\If{$J^{i} < \Jbest^{i-1}$}
		\State{$\ubtrsig{j} = \max \left(\min \left( \rhoInc \cdot \gtrsig^{j-1}, \ubtrsigMax \right), \ubtrsig{j-1} \right)$}
		\ElsIf{$J^{i-1} \leq \Jbest^{i-1}$}
		\State{$\ubtrsig{j} = \ubtrsig{j-1}$}
		\Else 
		\State{$\ubtrsig{j} = \max \left( \rhoDec \cdot \ubtrc{j-1}, \ubtrsigMin \right)$}
		\EndIf
		\EndIf
		\State{\textbf{Return:} $\ubtrsig{j}$, \ie the upper bound of the $\sigma$ trust region}
	\end{algorithmic}
\end{algorithm}

\begin{figure}[t!]
	\centering
	\begin{subfigure}[t]{0.45\textwidth}
		\includegraphics[width=\textwidth]{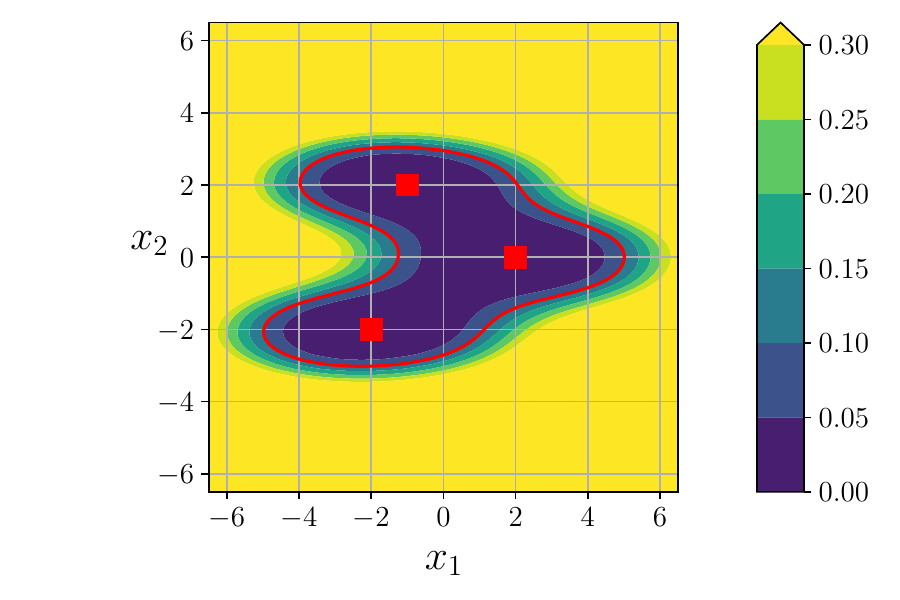}
		\caption{$\gammavec = [0.25, 0.75]^\top$}
		\label{Fig_TrustRegion_a}
	\end{subfigure}	
	\hspace{5pt}
	\begin{subfigure}[t]{0.45\textwidth}
		\includegraphics[width=\textwidth]{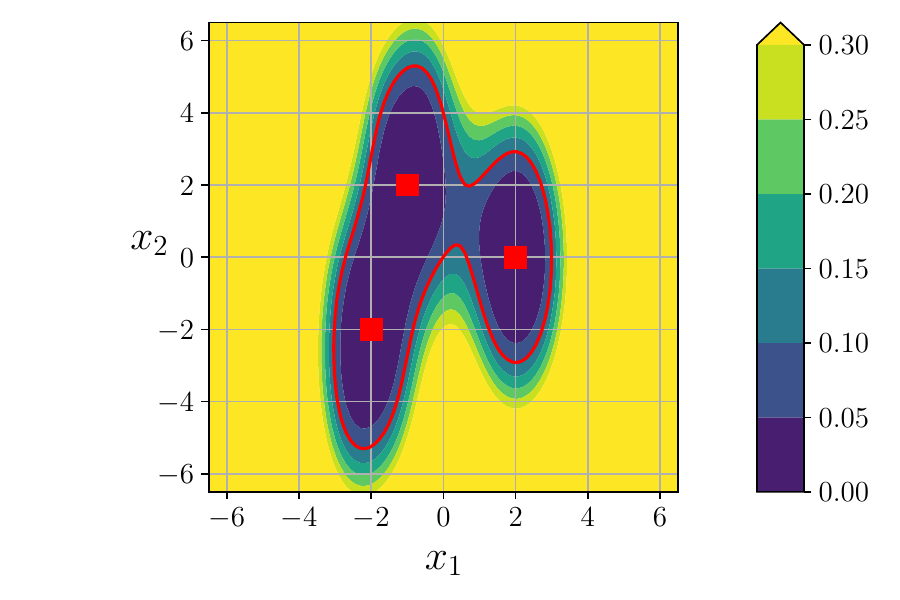}
		\caption{$\gammavec = [0.75, 0.25]^\top$}
		\label{Fig_TrustRegion_b}
	\end{subfigure}	
	%
	\caption[Examples of probabilistic trust regions.]{Trust region $\gtrsig(\xvec)$ from \Eq{Eq_tr_sigma_val} with the contour for $\frac{\sigmaGP_f^2(\xvec)}{\sigKf^2}$ using the Gaussian kernel with $\hpstdfval = \hpstdfgrad = 0$ and with red squares indicating the evaluation points. The region within the red line is where the constraint is satisfied for $\ubtrsig{} = 0.1$.}
	\label{Fig_TrustRegion}
\end{figure}

An example of the $\sigma$ trust region is shown in \Fig{Fig_TrustRegion} with three evaluation points. The contours represent $\frac{\sigmaGP_f^2(\xvec)}{\sigKf^2}$ and the red line is for $\ubtrsig{} = 0.1$. The trust region does not depend on the function and gradient evaluations directly, as is evident from \Eq{Eq_frac_sigf_sigK}, but it is sensitive to the values of the hyperparameters, as is evident from \Figs{Fig_TrustRegion_a}{Fig_TrustRegion_b}, which use different values for $\gammavec$. 

\subsection{Acquisition function minimization} \label{Sec_LocalOptz_Framework_AcqOptz}

The evaluation of the acquisition function is significantly less expensive than the evaluation of the marginal log likelihood. The most computationally expensive calculations involve $\Sigmag^{-1}$ from the evaluation of $\muGP_f$ and $\sigmaGP_f$ from \Eqs{Eq_muGp_w_grad}{Eq_sigmaGp_w_grad}, respectively. The covariance matrix $\Sigmag$ depends only on $\Xdata$ from Algorithm~\ref{Alg_DataRegion} and the hyperparameters such as $\gammavec$; it is independent of where the posterior of the GP is evaluated in the parameter space. Therefore, the Cholesky decomposition of the covariance matrix only needs to be calculated once for the acquisition function to be evaluated at several points in the parameter space. The evaluation of the GP only requires matrix-vector products, which are relatively inexpensive for modest $\ndata$ and $n_d$. Since the acquisition function can be multimodal, gradient-based multistart optimization is used. The minimization of the acquisition function for an unconstrained optimization problem is given by
\begin{align*}
	\xvec^{i+1} = \argmin_{\xvec} q(\xvec) \quad \text{s.t.} \quad \quad 
	\gtrc(\xvec) 	&\leq \ubtrc{j} \yesnumber \label{Eq_acq_optz_uncon} \\
	\gtrsig(\xvec) 	&\leq \ubtrsig{j},
\end{align*}
where $\xvec^{i+1}$ is the next point in the parameter space where the function and gradient will be evaluated, $q(\xvec)$ is the acquisition function, $\gtrc(\xvec)$ and $\gtrsig(\xvec)$ are the circular and $\sigma$ trust regions from \Eqs{Eq_tr_circle_val}{Eq_tr_sigma_val}, respectively, and $\ubtrc{j}$ and $\ubtrsig{j}$ are the upper bounds for the circular and $\sigma$ trust regions for the $j$-th optimization iteration, respectively. The circular trust region is active during all iterations, while the $\sigma$ trust region is only active once $\ndata \geq 10$, as indicated in Algorithm~\ref{Alg_TrustRegion_sigma}. The expected improvement acquisition function from \Eq{Eq_acq_EI} is used by default for $q(\xvec)$. Other acquisition functions were tested and the results are in \Sec{Sec_LocalOptz_Studies_AcqFun}.

\Eq{Eq_acq_optz_uncon} is solved with Algorithm~\ref{Alg_AcqSol}. Several initial points are selected to start the local minimization of \Eq{Eq_acq_optz_uncon}. Half of the starting points come from a Latin hypercube sampling centred at $\xbest$, and the other half are from the evaluation points in the data region with the lowest evaluations of the merit function. For each of the starting points, the gradient-based optimizer SciPy SLSQP is used to minimize \Eq{Eq_acq_optz_uncon}. Finally, the solution to \Eq{Eq_acq_optz_uncon} with the lowest acquisition function evaluation is the one that is returned by the Bayesian optimizer. This is the point in the parameter space where the function and gradient evaluations will be evaluated next.

\begin{algorithm}[t!]
	\caption{Acquisition function minimization} 
	\label{Alg_AcqSol}
	\begin{algorithmic}[1]
		\Statex{\textbf{Required:} The evaluation points $\Xdata$ from the data region from Algorithm~\ref{Alg_DataRegion} along with their merit function evaluations $\Jdata$ and the upper bound for the circular trust region $\ubtrc{j}$ from Algorithm~\ref{Alg_TrustRegion_circle}, where $j$ indicates the optimization iteration}
		\Statex{\textbf{Select:} $\nLhs (5)$, $n_{\text{best}}(5)$, and $q (\qEI)$, with default values indicated in parentheses}
		\State{Identify $\xbest$, \ie $\xbest = \left( \Xdata \right)_{i^*:}$, where $i^* = \argmin_i \left(\Jdata \right)_i$}
		\State{$\ub{\xvec}_{\text{LHS}} = \xbest + \ubtrc{j}$}
		\State{$\lb{\xvec}_{\text{LHS}} = \xbest - \ubtrc{j}$}
		\State{$\Xacq = \text{LHS}(\nLhs, \lb{\xvec}_{\text{LHS}}, \ub{\xvec}_{\text{LHS}})$}
		\State{For $i \in \{1, \ldots, \ndata \}$ append the $i$-th row of $\Xdata$ as a row to $\Xacq$ if $\left(\Jdata \right)_i$ is one of the $n_{\text{best}}$-th lowest entries in $\Jdata$}
		\For{$i \in \{1, \ldots, (\nLhs + n_{\text{best}}) \}$}
		\State{Solve \Eq{Eq_acq_optz_uncon} with SciPy trust-constr for $\xsol{i}$ using the $i$-th row of $\Xacq$ as the initial solution}
		\EndFor
		\State{\textbf{Return:} The $\xsol{i}$ with the lowest acquisition function evaluation and that satisfies the constraints in \Eq{Eq_acq_optz_uncon}}
	\end{algorithmic}
\end{algorithm}

\section{Unconstrained test cases} \label{Sec_LocalOptz_TestCases}

The test cases used to benchmark the Bayesian optimizer and compare it with conjugate-gradient and quasi-Newton optimizers are the following functions:
\begin{alignat}{3} 
	\text{Quadratic function:} \quad 
		& f(\xvec) 
		&&= \frac{1}{2} \left(\xvec - \one \right)^\top \A \left(\xvec - \one \right) \label{Eq_Quadratic_fun} \\
	\text{Bowl function:} \quad 
		& f(\xvec) 
		&&= 1 - e^{-\frac{1}{2} \left(\xvec - \one \right)^\top \A \left(\xvec - \one \right)} 
	+ \frac{ \| \xvec - \one \|_2^2}{100} 
	+ \frac{ \| \xvec - \one \|_4^4}{1000} \label{Eq_Bowl} \\
	\text{Rosenbrock function:} \quad 
		& f(\xvec) 
		&&= \sum_{i=1}^{n_d-1} \left[ a \left(x_{i+1} - x_i^2 \right)^2 + \left( 1 - x_i \right)^2 \right], \label{Eq_Rosenbrock}
\end{alignat}
where $\A$ is a symmetric matrix with entries given by $a_{ij} = \frac{1}{10} e^{-\frac{1}{2}(i-j)^2} \, \forall \, i,j \in \{1, \ldots, n_d \}$, and $a > 0$ is a constant for the Rosenbrock function. The exponential term in the bowl function ensures it is non-polynomial while the two polynomial terms are included to avoid having vanishing gradients far from the minimum. The Rosenbrock, quadratic, and bowl functions are all unimodal with their respective minima evaluating to zero at $\xvec = \one$. These functions can be used with an arbitrary number of dimensions, which enables the impact of dimensionality to be studied for the Bayesian optimizer.

\begin{figure}[t!]
	\centering
	\begin{subfigure}[t]{0.328\textwidth}
		\includegraphics[width=\textwidth]{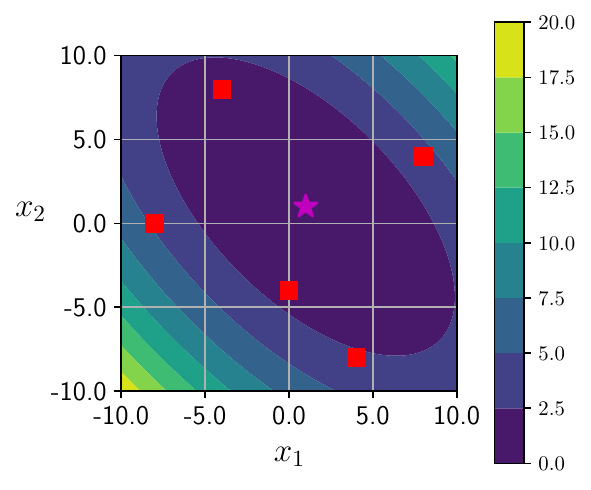}
		\caption{Quadratic function}
		\label{Fig_BoTestCases_Quad_zoomF}
	\end{subfigure}	
	\begin{subfigure}[t]{0.328\textwidth}
		\includegraphics[width=\textwidth]{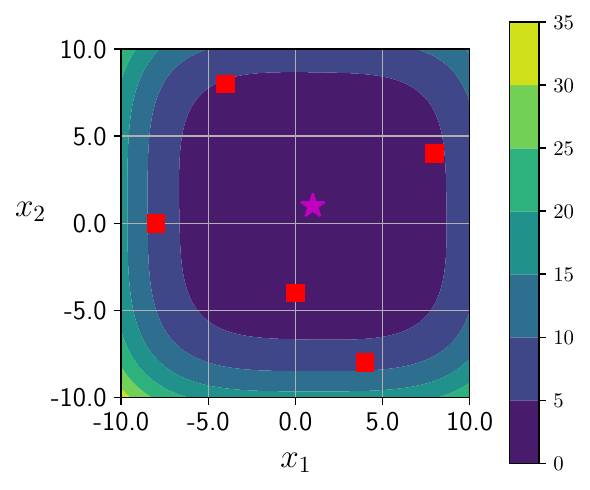}
		\caption{Bowl function} 
		\label{Fig_BoTestCases_Bowl_zoomF}
	\end{subfigure}	
	\begin{subfigure}[t]{0.328\textwidth}
		\includegraphics[width=\textwidth]{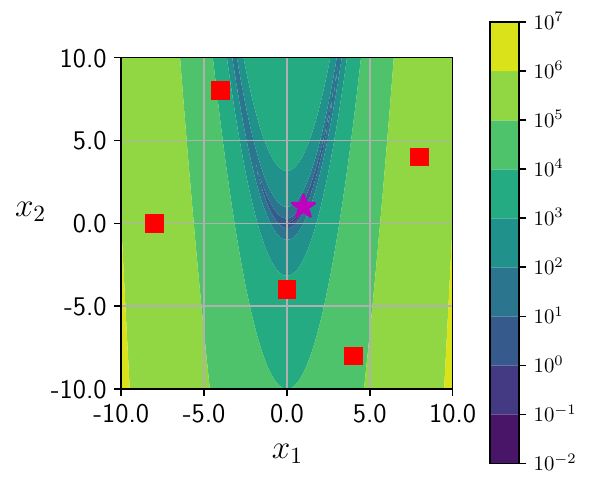}
		\caption{Rosenbrock $a=100$}
		\label{Fig_BoTestCases_Rosen_a100_zoomF}
	\end{subfigure}	
	\begin{subfigure}[t]{0.328\textwidth}
		\includegraphics[width=\textwidth]{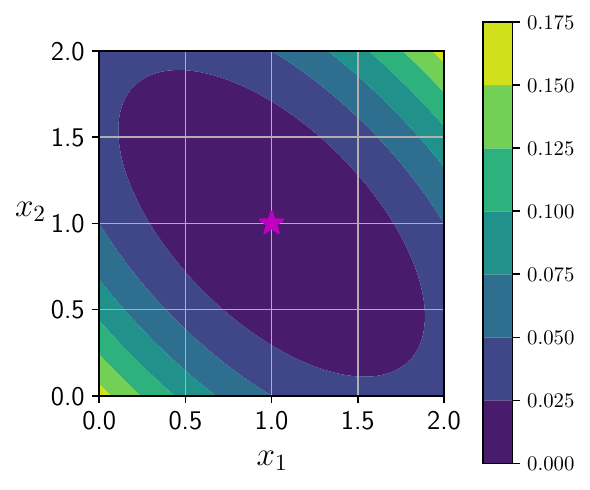}
		\caption{Quadratic function}
		\label{Fig_BoTestCases_Quad_zoomT}
	\end{subfigure}	
	\begin{subfigure}[t]{0.328\textwidth}
		\includegraphics[width=\textwidth]{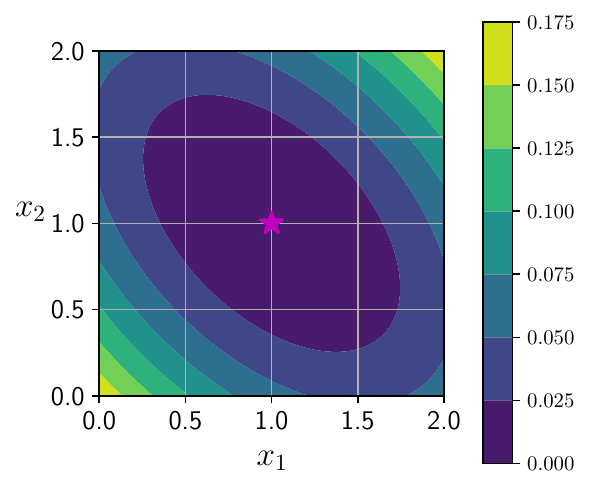}
		\caption{Bowl function} 
		\label{Fig_BoTestCases_Bowl_zoomT}
	\end{subfigure}	
	\begin{subfigure}[t]{0.328\textwidth}
		\includegraphics[width=\textwidth]{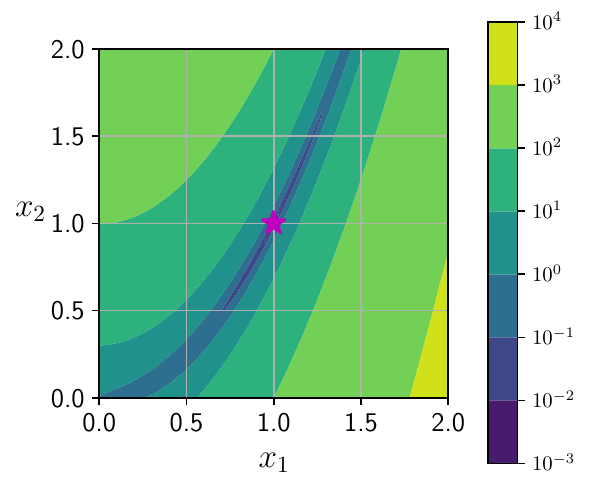}
		\caption{Rosenbrock $a=100$}
		\label{Fig_BoTestCases_Rosen_a100_zoomT}
	\end{subfigure}	
	\caption[Unconstrained test cases with the optimization starting points from a Latin hypercube sampling.]{Plots for the two-dimensional quadratic, bowl, and $a=100$ Rosenbrock functions from \Eqss{Eq_Quadratic_fun}{Eq_Bowl}{Eq_Rosenbrock}, respectively. The red squares in the subfigures of the top row indicate the starting points for the optimizer that were selected with a Latin hypercube sampling and the minimum of each function is labelled with a magenta star. The subfigures in the bottom row are centred at the minimum of the test cases.}
	\label{Fig_BoTestCases}
\end{figure}

The two-dimensional quadratic, bowl, and $a=100$ Rosenbrock functions for $-10 \leq x_1, x_2 \leq 10$ can be seen in \Figss{Fig_BoTestCases_Quad_zoomF}{Fig_BoTestCases_Bowl_zoomF}{Fig_BoTestCases_Rosen_a100_zoomF}, respectively. The quadratic function is straightforward to minimize, particularly for the quasi-Newton optimizers that approximate the function of interest with their own quadratic function. This will thus enable the Bayesian optimizer to be compared to a test case that is ideally suited for quasi-Newton optimizers. The bowl function should also be straightforward to minimize but unlike the quadratic function, it is not a polynomial. It is also clear from \Figs{Fig_BoTestCases_Quad_zoomT}{Fig_BoTestCases_Bowl_zoomT} that the quadratic and bowl functions are very similar near their respective minima. Finally, the Rosenbrock function is the most challenging to minimize and is a common optimization test case. It has a valley that is generally quickly found by all of the optimizers. However, traversing this valley to the minimum at $\xvec = \one$ is challenging to do efficiently, \ie with few function evaluations. Larger values of the parameter $a$ make the walls of the valley steeper, which makes the Rosenbrock function more challenging to minimize.

For each of the test cases five independent optimization runs are performed to avoid having outlier results. Each optimization run is initiated with one evaluation point. The starting points for the independent optimizations are selected from a Latin hypercube sampling and can be seen as the red squares in \Fig{Fig_BoTestCases}. The Bayesian optimizer could be started with an arbitrary number of initial evaluation points. However, it was initiated with only one starting point since it is compared in \Sec{Sec_LocalOptz_UnconQN} to deterministic optimizers, which can only be started with a single starting point.

The Bayesian optimizer is being developed to be used for the optimization of problems with expensive function evaluations, such as aerodynamic shape optimization that involves computationally intensive flow evaluations. The cost of the function and gradient evaluations will thus be far greater than the computational cost of the optimizers themselves for these types of problems. Therefore, the optimizers are compared based on the number of function and gradient evaluations. One iteration represents one function and one gradient evaluation. For the objective plots, the lowest evaluated objective up to the given iteration is shown. This is done, rather than comparing the function evaluation at each iteration, to make it easier to compare the progress of the optimizers. The normalized optimality will also be plotted, which is simply the $\ell_2$ norm of the gradient for an unconstrained test case normalized by the optimality of the starting evaluation point. Just like the plot of the objective, the normalized optimality of the evaluation point with the lowest objective evaluation is shown at each iteration.

\section{Unconstrained studies for the Bayesian optimizer} \label{Sec_LocalOptz_Studies}

In this section several settings for the Bayesian optimizer are varied to determine their impact on the optimization results. Each study focuses on one setting and the results from the default settings are shown in green. The default settings involve using the data region from Algorithm~\ref{Alg_DataRegion} with $\ndata = 20$. The starting points for the maximization of the marginal log-likelihood come from Algorithm~\ref{Alg_GpHyperparameters} and the preconditioning method detailed in \Sec{Sec_LocalOptz_Framework_Precon} with $\condmax = 10^{10}$. The upper bounds for the circular and $\sigma$ trust regions are selected with Algorithms~\ref{Alg_TrustRegion_circle} and \ref{Alg_TrustRegion_sigma}, respectively, and the expected improvement acquisition function $\qEI(\xvec)$ is minimized with Algorithm~\ref{Alg_AcqSol}. 

\subsection{Selecting $\nclose$ for the data region} \label{Sec_LocalOptz_Studies_nclose}

\begin{figure}[t!]
	\centering
	\begin{subfigure}[t]{0.328\textwidth}
		\includegraphics[width=\textwidth]{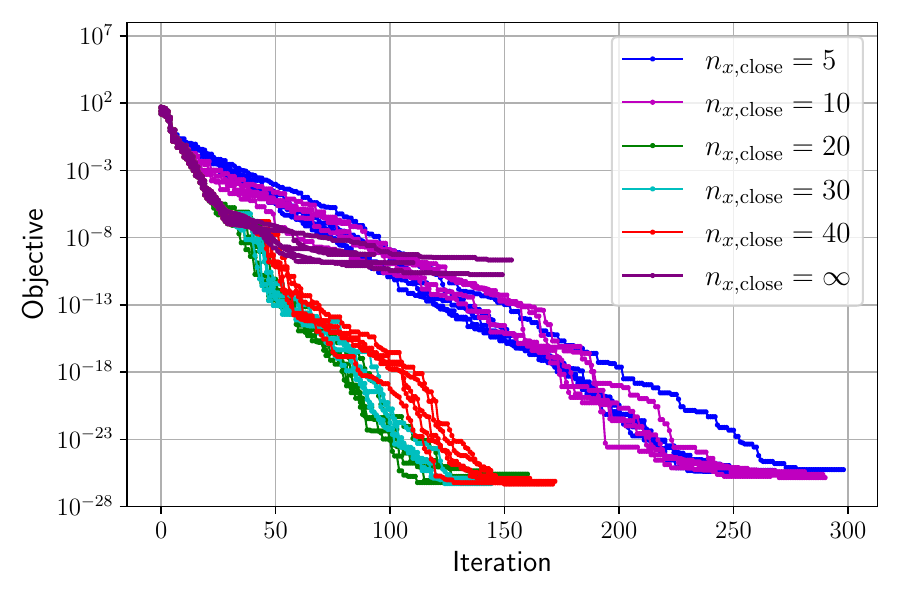}
		\caption{Objective: quadratic function}
		\label{Fig_Study_nclose_Obj_Quad_d20}
	\end{subfigure}	
	\begin{subfigure}[t]{0.328\textwidth}
		\includegraphics[width=\textwidth]{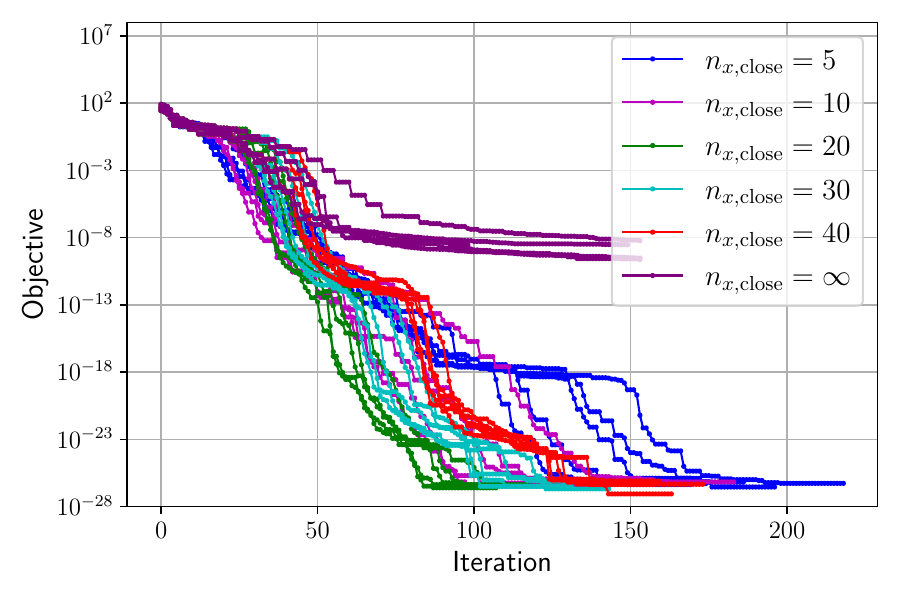}
		\caption{Objective: bowl function}
		\label{Fig_Study_nclose_Obj_Bowl_d20}
	\end{subfigure}	
	\begin{subfigure}[t]{0.328\textwidth}
		\includegraphics[width=\textwidth]{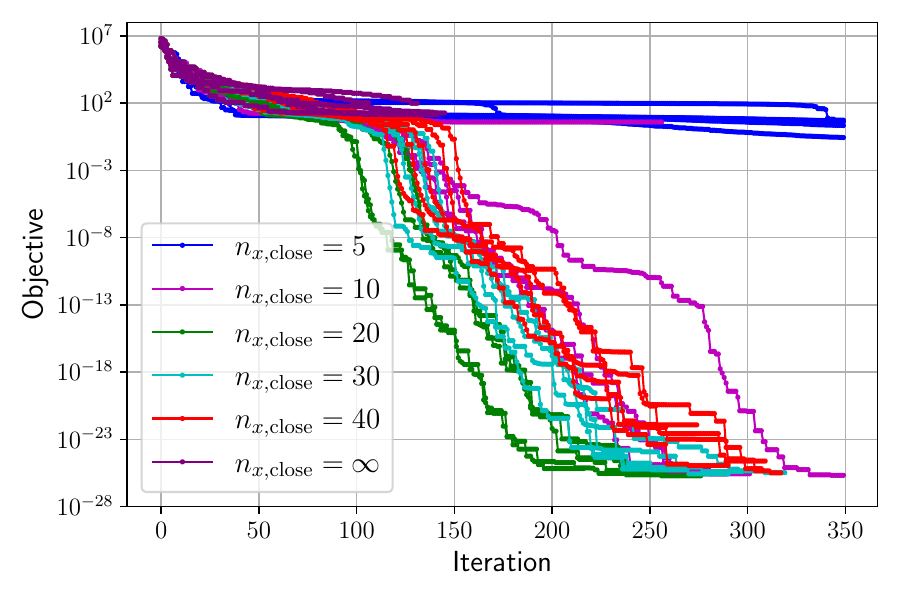}
		\caption{Objective: Rosenbrock $a=100$}
		\label{Fig_Study_nclose_Obj_RosenA100_d20}
	\end{subfigure}	
	\begin{subfigure}[t]{0.328\textwidth}
		\includegraphics[width=\textwidth]{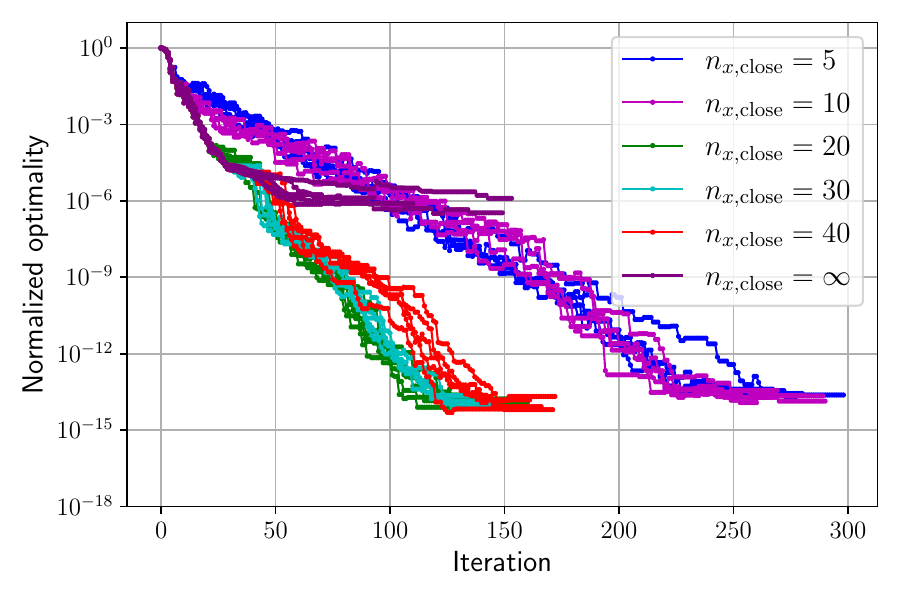}
		\caption{Optimality: quadratic function}
		\label{Fig_Study_nclose_Opt_Quad_d20}
	\end{subfigure}	
	\begin{subfigure}[t]{0.328\textwidth}
		\includegraphics[width=\textwidth]{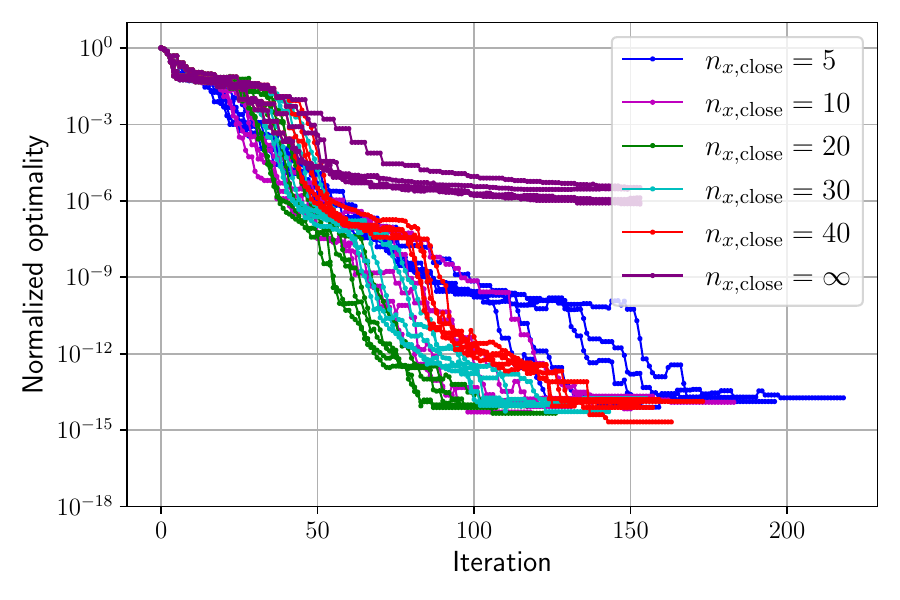}
		\caption{Optimality: bowl function}
		\label{Fig_Study_nclose_Opt_Bowl_d20}
	\end{subfigure}	
	\begin{subfigure}[t]{0.328\textwidth}
		\includegraphics[width=\textwidth]{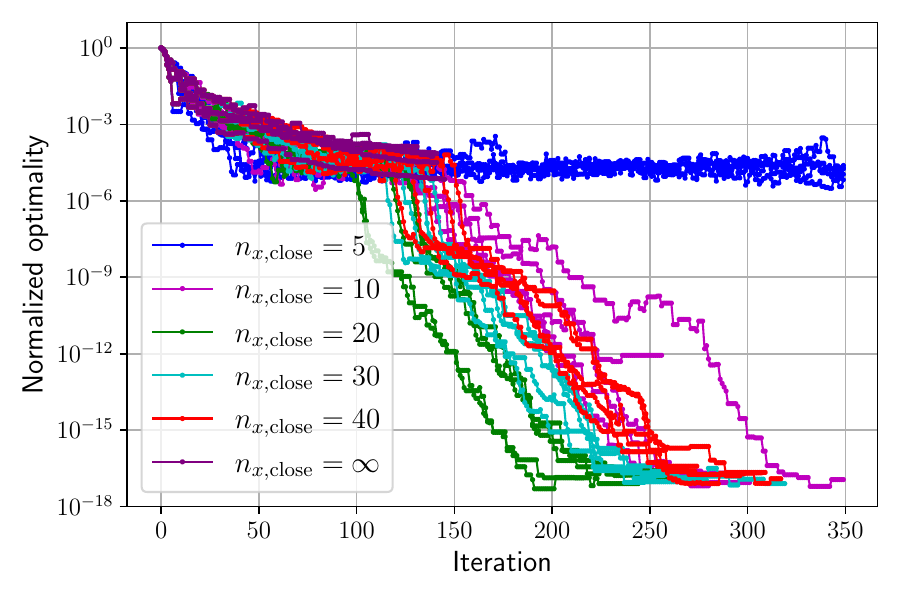}
		\caption{\mbox{Optimality: Rosenbrock $a=100$}}
		\label{Fig_Study_nclose_Opt_RosenA100_d20}
	\end{subfigure}	
	\caption{Unconstrained study with different $\nclose$ for the data region from Algorithm~\ref{Alg_DataRegion} for the Bayesian optimizer. The test cases are the quadratic, bowl, and Rosenbrock functions from \Eqss{Eq_Quadratic_fun}{Eq_Bowl}{Eq_Rosenbrock}, respectively, with $n_d=20$. The case with $\nclose = \infty$ uses all of the function and gradient evaluations to form the surrogate, i.e. no data region is used.}
	\label{Fig_Study_nclose}
\end{figure}

In this section the data region is studied and various values for $\nclose$ from Algorithm~\ref{Alg_DataRegion} are investigated. \Fig{Fig_Study_nclose} shows the minimization of the $n_d = 20$ quadratic, bowl, and Rosenbrock $a=100$ functions from \Sec{Sec_LocalOptz_TestCases}. The subfigures in the top row of \Fig{Fig_Study_nclose} show the objective, while the bottom row displays the normalized optimality. For the case with $\nclose = \infty$, all of the evaluation points are used by the Bayesian optimizer and this is equivalent to not having a data region. It is clear from \Fig{Fig_Study_nclose} that not using a data region results in the Bayesian optimizer having the worst result, \ie the highest final objective and optimality. From 
\Figss{Fig_Study_nclose_Obj_Quad_d20}{Fig_Study_nclose_Obj_Bowl_d20}{Fig_Study_nclose_Obj_RosenA100_d20} we can see that the Bayesian optimizer is not able to reduce the optimality by more than six orders of magnitude when a data region is not used, \ie $\nclose = \infty$ from Algorithm~\ref{Alg_DataRegion}. It is also problematic when there are too few points in the data region. For $\nclose = 5$ it takes significantly more iterations for the Bayesian optimizer to achieve the same optimality reduction as the Bayesian optimizer using a larger finite $\nclose$ for both the quadratic and bowl functions. Moreover, the use of $\nclose = 5$ for the Rosenbrock function with $a=100$ results in an optimality that cannot be reduced further than 6 orders of magnitude, as seen in \Fig{Fig_Study_nclose_Opt_RosenA100_d20}.

\begin{figure}[t!]
	\centering
	\begin{subfigure}[t]{0.328\textwidth}
		\includegraphics[width=\textwidth]{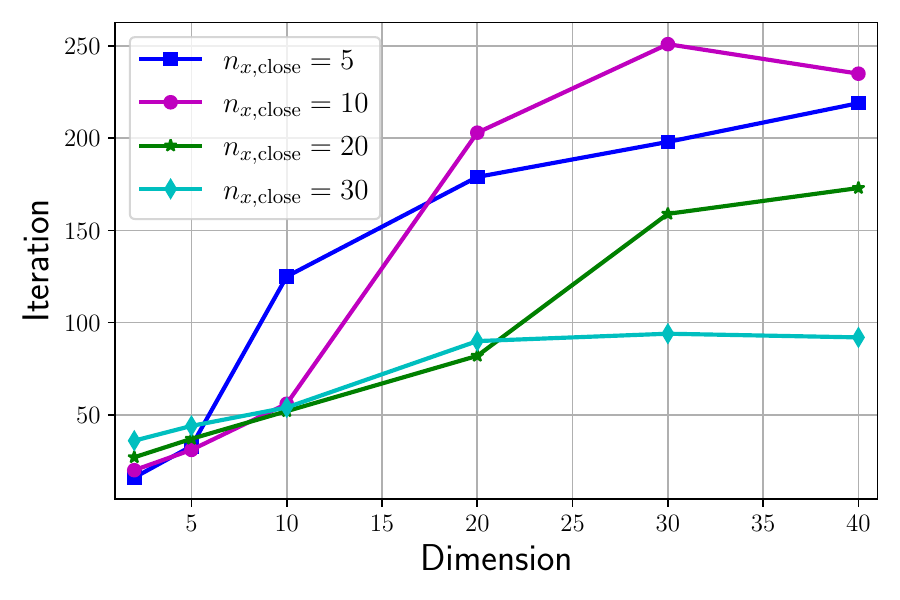}
		\caption{Quadratic function: \Eq{Eq_Quadratic_fun}}
	\end{subfigure}	
	\begin{subfigure}[t]{0.328\textwidth}
		\includegraphics[width=\textwidth]{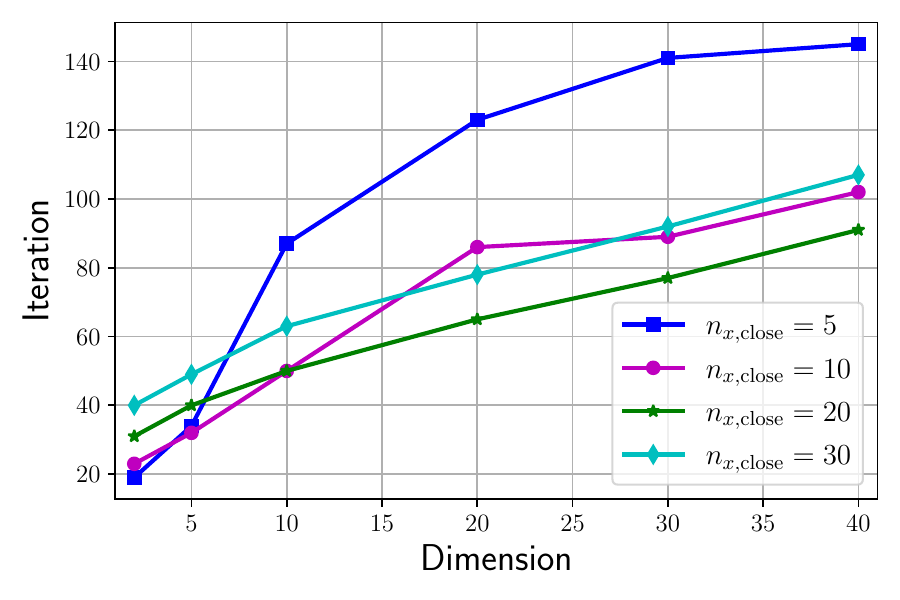}
		\caption{Bowl function: \Eq{Eq_Bowl}}
		\label{Fig_Study_nclose_OptTol_RosenA1}
	\end{subfigure}	
	\begin{subfigure}[t]{0.328\textwidth}
		\includegraphics[width=\textwidth]{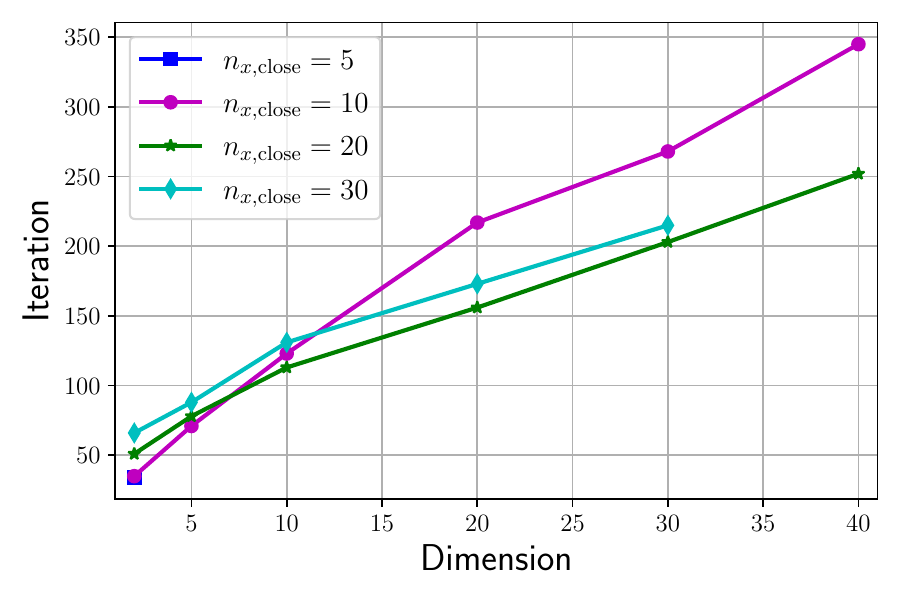}
		\caption{Rosenbrock $a=100$: \Eq{Eq_Rosenbrock}}
		\label{Fig_Study_nclose_OptTol_RosenA100}
	\end{subfigure}	
	\caption{Unconstrained study for the Bayesian optimizer with different $\nclose$ for the data region from Algorithm~\ref{Alg_DataRegion}. The plots indicate the required median number of iterations for the 25 independent optimization runs to reduce the objective below $10^{-5}$ and the optimality by 10 orders of magnitude.}
	\label{Fig_Study_nclose_OptTol}
\end{figure}

\Fig{Fig_Study_nclose_OptTol} plots the median number of iterations required for the Bayesian optimizer using different $\nclose$ to achieve a 10-order reduction in the normalized optimality and an objective below $10^{-5}$ for the quadratic, bowl, and Rosenbrock functions with $n_d \in \{2, 5, 10, 20, 30, 40\}$. The best results for the three test cases are achieved with $\nclose = 20$ or $\nclose = 30$. In contrast, the Bayesian optimizer using smaller or larger values of $\nclose$ either takes more iterations to achieve the desired tolerance or is not able to do so at all. For example, when $\nclose = 5$ the Bayesian optimizer does not achieve the desired tolerance for any of the initial conditions for the $n_d \geq 5$ Rosenbrock test cases. Also, the Bayesian optimizer does not achieve the desired tolerance for the Rosenbrock function with $n_d = 40$ when $\nclose = 30$. The default value for the Bayesian optimizer was selected to be $\nclose = 20$ since it provides good performance and achieves the desired tolerance in all cases.

\subsection{Selecting the unconstrained acquisition function} \label{Sec_LocalOptz_Studies_AcqFun}

\begin{figure}[t!]
	\centering
	\begin{subfigure}[t]{0.328\textwidth}
		\includegraphics[width=\textwidth]{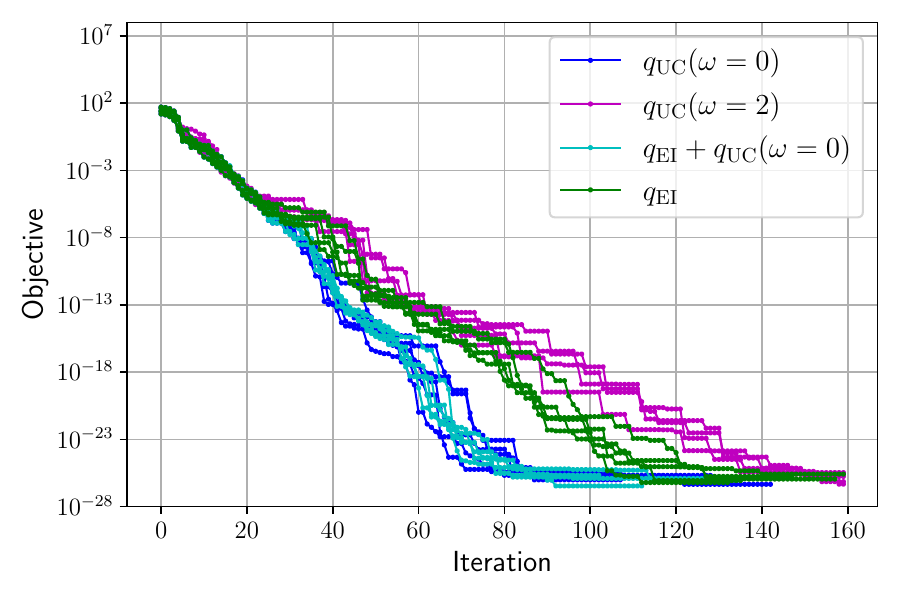}
		\caption{Objective: quadratic function}
		\label{Fig_Study_AcqUncon_Obj_Quad_d20}
	\end{subfigure}	
	\begin{subfigure}[t]{0.328\textwidth}
		\includegraphics[width=\textwidth]{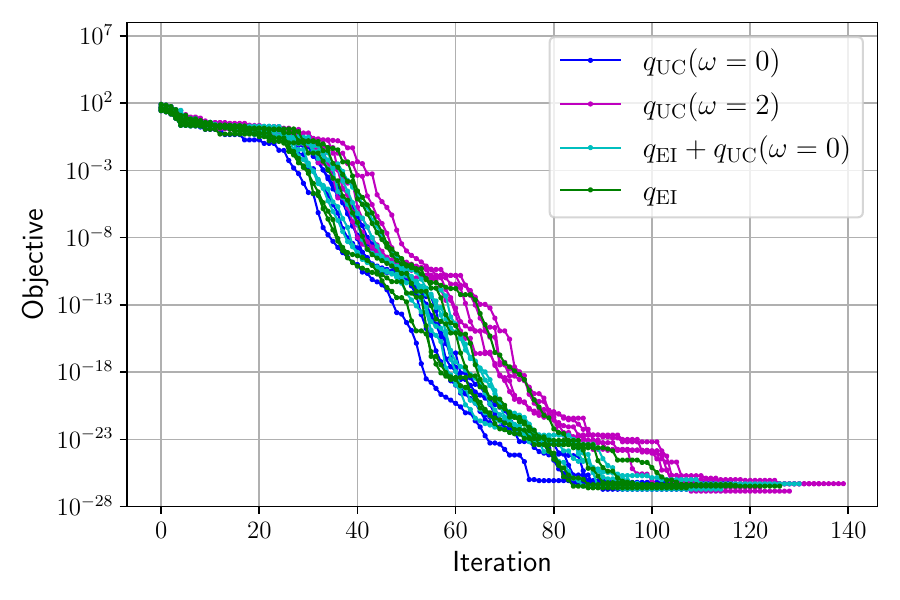}
		\caption{Objective: bowl function}
		\label{Fig_Study_AcqUncon_Obj_Bowl_d20}
	\end{subfigure}	
	\begin{subfigure}[t]{0.328\textwidth}
		\includegraphics[width=\textwidth]{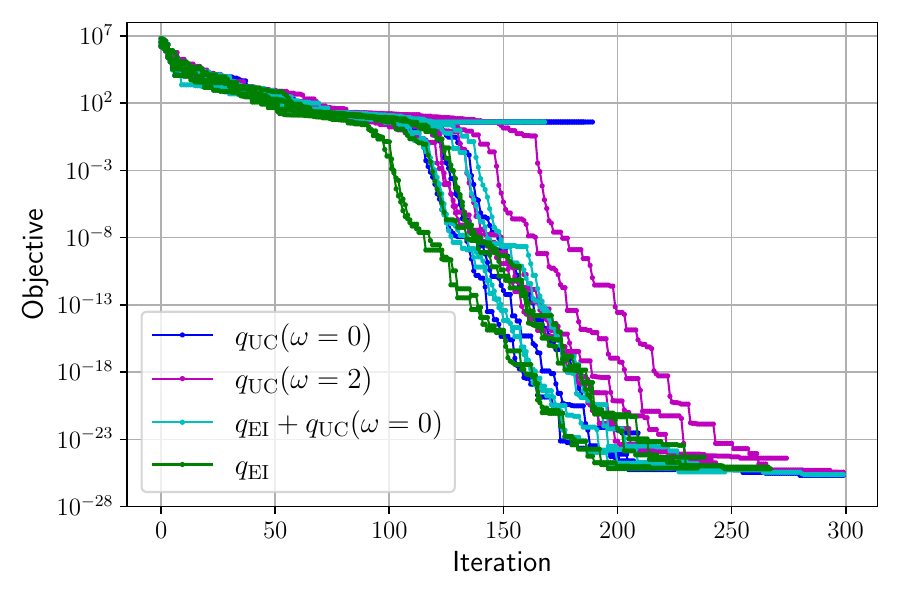}
		\caption{Objective: Rosenbrock $a=100$}
		\label{Fig_Study_AcqUncon_Obj_RosenA100_d20}
	\end{subfigure}	
	\begin{subfigure}[t]{0.328\textwidth}
		\includegraphics[width=\textwidth]{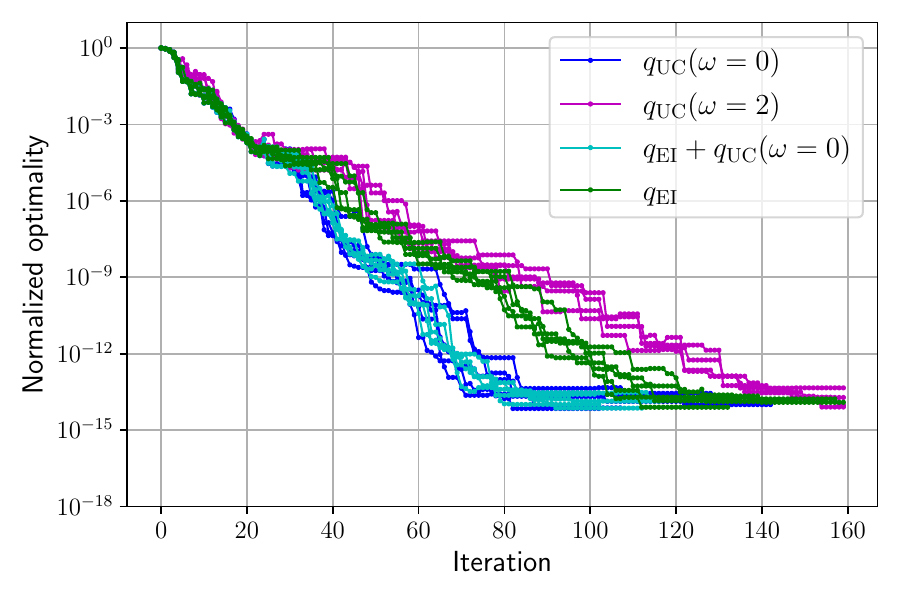}
		\caption{Optimality: quadratic function}
		\label{Fig_Study_AcqUncon_Opt_Quad_d20}
	\end{subfigure}	
	\begin{subfigure}[t]{0.328\textwidth}
		\includegraphics[width=\textwidth]{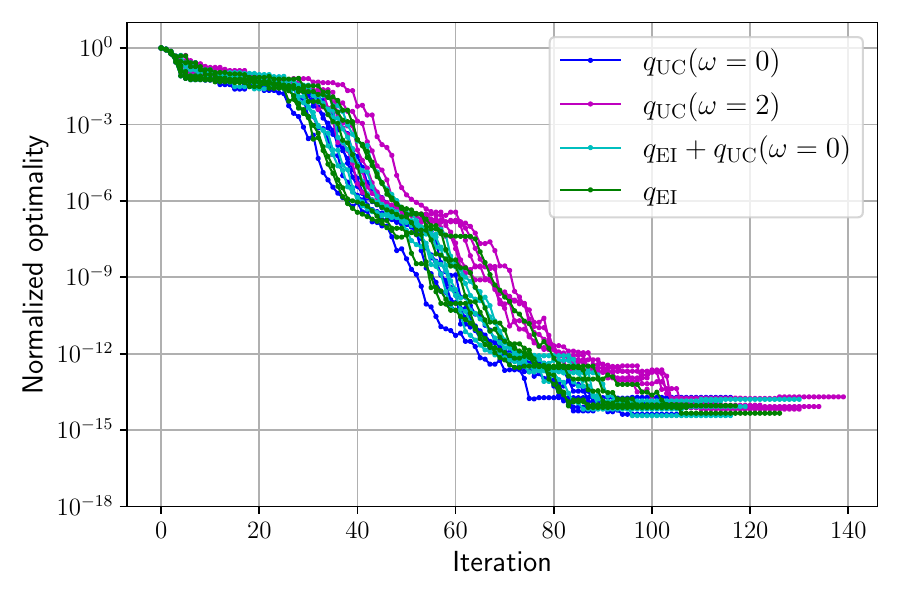}
		\caption{Optimality: bowl function}
		\label{Fig_Study_AcqUncon_Opt_Bowl_d20}
	\end{subfigure}	
	\begin{subfigure}[t]{0.328\textwidth}
		\includegraphics[width=\textwidth]{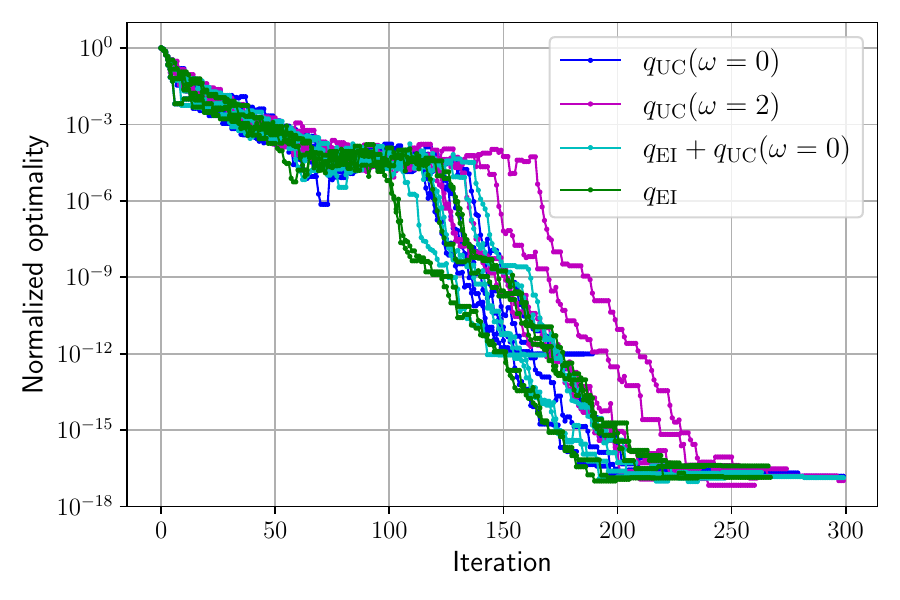}
		\caption{\mbox{Optimality: Rosenbrock $a=100$}}
		\label{Fig_Study_AcqUncon_Opt_RosenA100_d20}
	\end{subfigure}	
	\caption[Unconstrained optimization with different acquisition functions for $n_d=20$.]{Unconstrained Bayesian optimization with the use of different acquisition functions. The test cases are the twenty-dimensional quadratic, bowl, and Rosenbrock functions from \Sec{Sec_LocalOptz_TestCases}.}
	\label{Fig_Study_AcqUncon}
\end{figure}

In this section the use of different acquisition functions for local unconstrained Bayesian optimization is investigated. The upper confidence and expected improvement acquisition functions from \Eqs{Eq_acq_UC}{Eq_acq_EI}, respectively, are considered independently and together in \Fig{Fig_Study_AcqUncon}. The upper confidence acquisition function is used with $\omega = 2$ and with $\omega = 0$; in the latter case, it is simply equal to the negative of the mean of the posterior of the GP, \ie $q_{\text{UC}}(\xvec; \omega = 0) = -\muGP_f(\xvec)$. The test cases are again the $n_d=20$ quadratic, bowl, and $a=100$ Rosenbrock functions from \Sec{Sec_LocalOptz_TestCases}. The top row of subfigures shows the objective evaluations, while the bottom three subfigures are for the normalized optimality. The Bayesian optimization results for the $n_d = 20$ test cases with the use of the upper confidence acquisition function with $\omega = 2$ are noticeably worse than for the other acquisition functions. In contrast, all of the other acquisition functions provide similar results.

\begin{figure}[t!]
	\centering
	\begin{subfigure}[t]{0.328\textwidth}
		\includegraphics[width=\textwidth]{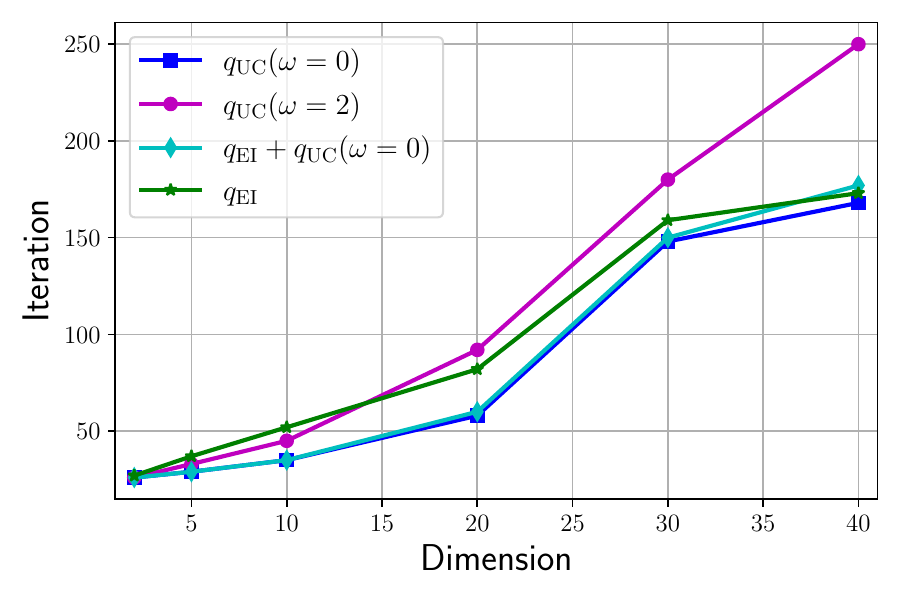}
		\caption{Quadratic function: \Eq{Eq_Quadratic_fun}}
		\label{Fig_Study_AcqUncon_OptTol_Quad}
	\end{subfigure}	
	\begin{subfigure}[t]{0.328\textwidth}
		\includegraphics[width=\textwidth]{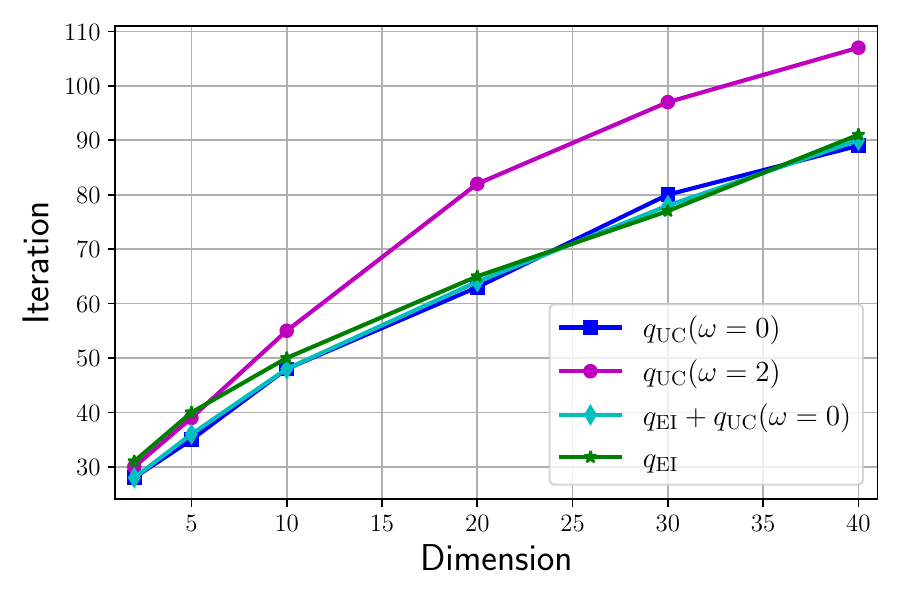}
		\caption{Bowl function: \Eq{Eq_Bowl}}
		\label{Fig_Study_AcqUncon_OptTol_Bowl}
	\end{subfigure}	
	\begin{subfigure}[t]{0.328\textwidth}
		\includegraphics[width=\textwidth]{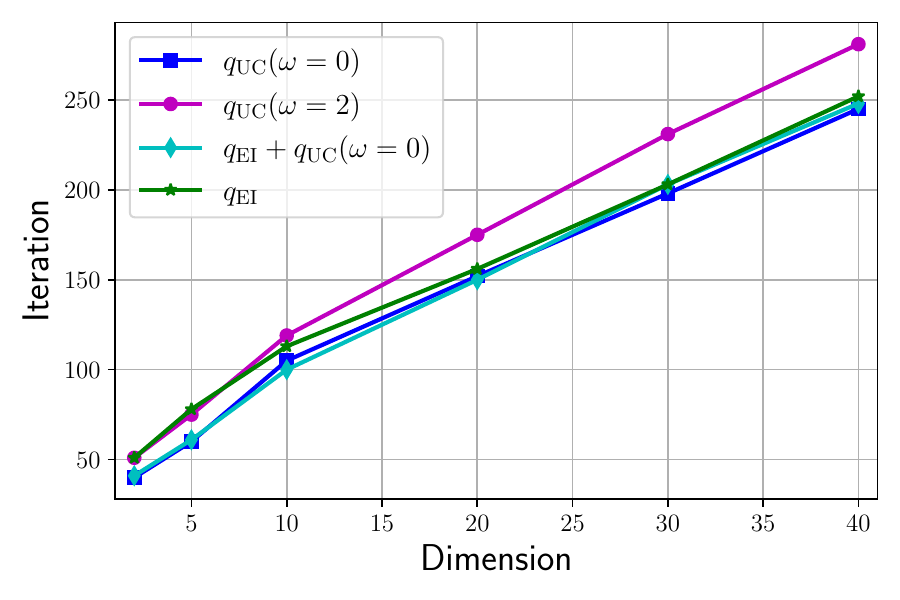}
		\caption{Rosenbrock $a=100$: \Eq{Eq_Rosenbrock}}
		\label{Fig_Study_AcqUncon_OptTol_RosenA100}
	\end{subfigure}	
	\caption{Median number of iterations for the Bayesian optimizer to reduce the objective evaluations below $10^{-5}$ and the optimality by 10 orders of magnitude with different acquisition functions for 25 independent optimization runs. The acquisition functions $q_{\text{UC}}(\omega)$ and $q_{\text{EI}}$ are the upper confidence and expected improvement acquisition functions from \Eqs{Eq_acq_UC}{Eq_acq_EI}, respectively.}
	\label{Fig_Study_AcqUncon_OptTol}
\end{figure}

\Fig{Fig_Study_AcqUncon_OptTol} shows the median number of iterations required for the Bayesian optimizer to reduce the objective evaluations below $10^{-5}$ and the optimality by 10 orders of magnitude for the three test cases using the different acquisition functions. For the three test cases and $2 \leq n_d \leq 40$, all of the acquisition functions have similar performance, except for upper confidence with $\omega = 2$. Having $\omega = 2$ results in the acquisition function promoting exploration in regions that do not result in reductions of the objective evaluation. While the expected improvement acquisition function also promotes exploration, it only does so in regions with significant probabilities of improvement or with higher expected improvements, as the name implies. The default acquisition function for unconstrained optimization was selected to be the expected improvement function $\qEI$ since its use was effective for the three test cases considered.

\subsection{Summary of the studies for the unconstrained Bayesian optimizer} \label{Sec_LocalOptz_Studies_BoOverivew}

The default settings for the Bayesian optimizer are detailed in Algorithm~\ref{Alg_BoDefault}. These are the settings that were used for the Bayesian optimizer with data shown in green in the figures of the studies in the previous subsections. It was clear from \Sec{Sec_LocalOptz_Studies_nclose} that it is beneficial to use a data region such that only the function and gradient evaluations near $\xbest$ are used to select the hyperparameters and evaluate the GP's posterior. The Bayesian optimizer was able to reduce the optimality several additional orders of magnitude by using a data region. 

In \Sec{Sec_LocalOptz_Studies_AcqFun} different acquisition functions are considered and it was found that the expected improvement acquisition function from \Eq{Eq_acq_EI} provides the best results. In Appendix \ref{Sec_AppendixUnconStudy_CondMax} the same test cases were used with the following maximum condition numbers: $\condmax \in \{10^8, 10^{10}, 10^{12}, 10^{14}\}$. The Bayesian optimizer was found to generally be insensitive to different values for $\condmax$. A maximum condition number of $\condmax = 10^{10}$ was selected since it was found to provide good consistent results for the Bayesian optimizer.

The use of the Gaussian, Mat\'ern $\frac{5}{2}$, and rational quadratic kernels from \Eqss{Eq_kern_Gaussian}{Eq_kern_Mat5f2}{Eq_kern_RatQd}, respectively, is investigated in Appendix \ref{Sec_AppendixUnconStudy_Kernel}. The results indicate that there is not a big change in the performance of the Bayesian optimizer with the use of these three kernels for the infinitely differentiable test cases that were considered. While the results may vary for different test cases, the Gaussian kernel was selected since it is simple to use and the Bayesian optimizer was efficient when it was used.

\begin{algorithm}[t!]
	\caption{Local unconstrained optimization framework for the Bayesian optimizer} 
	\label{Alg_BoDefault}
	\begin{algorithmic}[1]
		\Statex{\textbf{Required:} At least one evaluation point with its function and gradient evaluation.}
		\Statex{\textbf{Select:} A desired convergence criterion, \eg a 10-order reduction in the optimality or no reduction of the objective function after 20 iterations.}
		\State{Select the evaluation points for the data region with Algorithm~\ref{Alg_DataRegion}}
		\State{Use the Gaussian kernel from \Eq{Eq_kern_Gaussian} and the preconditioning method from \Sec{Sec_LocalOptz_Framework_Precon} with $\condmax = 10^{10}$}
		\State{Select the hyperparameters of the GP using Algorithm~\ref{Alg_GpHyperparameters}}
		\State{Select the upper bounds for the circular and $\sigma$ trust regions with Algorithms~\ref{Alg_TrustRegion_circle} and \ref{Alg_TrustRegion_sigma}, respectively}
		\State{Use Algorithm~\ref{Alg_AcqSol} with the expected improvement acquisition function to select the next point in the parameter space where the objective and its gradient are evaluated}
		\State{Evaluate the next function and gradient evaluation}
		\State{Check the convergence criterion; if satisfied continue to next step, otherwise return to step 1}
		\State{\textbf{Return:} $\xbest$, \ie the evaluation point with the lowest merit function evaluation along with its objective evaluation and optimality}
	\end{algorithmic}
\end{algorithm}

Previous gradient-enhanced Bayesian optimizers have also been applied to the Rosenbrock function. Shende $\etal$ developed a gradient-enhanced Bayesian optimizer and applied it to the Rosenbrock function with $a=100$ and $n_d= 5$ \cite{shende_systematic_2022}. Their gradient-enhanced Bayesian optimizer was found to outperform the gradient-free Bayesian optimizer but it was not able to converge the objective below $10^{-4}$. In contrast, the gradient-enhanced Bayesian optimizer using Algorithm~\ref{Alg_BoDefault} is able to achieve an objective below $10^{-20}$ for the Rosenbrock function with $a=100$ and $n_d = 5$. Cheng and Zimmermann also developed a gradient-enhanced Bayesian optimizer that they applied to the Rosenbrock function with $a=100$ and $n_d=20$ \cite{cheng_gradient-enhanced_2023}. They started their optimizer with 20 initial evaluation points, ran it for 300 iterations, and repeated this five times with different initial evaluation points each time. The mean of the best function evaluations for the five runs after 320 iterations was 103.81, with a standard deviation of 26.43. In contrast, the gradient-enhanced Bayesian optimizer using Algorithm~\ref{Alg_BoDefault} requires fewer than 200 function evaluations to achieve an objective evaluation below $10^{-20}$ for all five independent runs, as seen in \Fig{Fig_Study_Kernel_Obj_RosenA100_d20} for the Gaussian kernel. In fact, fewer than 130 function evaluations are required to reduce the objective evaluation below the mean of 103.81 achieved by the gradient-enhanced optimizer from Cheng and Zimmermann \cite{cheng_gradient-enhanced_2023}.


\section{Comparing Bayesian, quasi-Newton, and conjugate-gradient optimizers with noise-free gradients} \label{Sec_LocalOptz_UnconQN}

The Bayesian optimizer with its default settings given in Algorithm~\ref{Alg_BoDefault} and the other algorithms it references is compared to the unconstrained optimizers MATLAB fminunc BFGS, SciPy, BFGS, and SciPy CG. The options ``StepTolerance'' and ``OptimalityTolerance'' for the MATLAB fminunc optimizer are reduced from their default values of $10^{-5}$ to $10^{-16}$ to ensure that the optimizer converged as deeply as possible. Similarly, for the SciPy, BFGS and CG optimizers the tolerance ``gtol'' for the gradient norm is set to $10^{-16}$. No other parameters are changed for the MATLAB and SciPy optimizers. The test cases are once again the quadratic, bowl, and $a=100$ Rosenbrock functions from \Sec{Sec_LocalOptz_TestCases}.

\begin{figure}[t!]
	\centering
	\begin{subfigure}[t]{0.328\textwidth}
		\includegraphics[width=\textwidth]{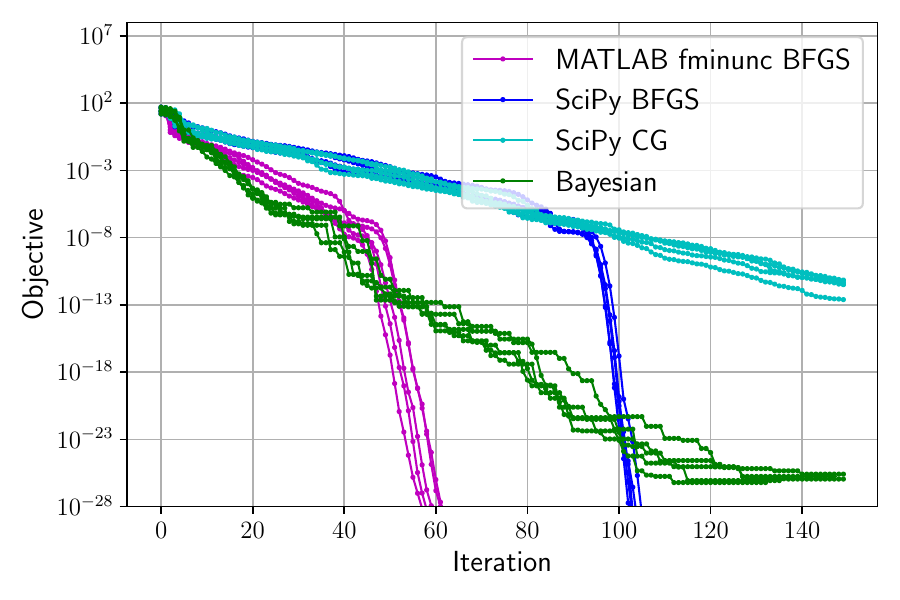}
		\caption{Objective: quadratic function}
		\label{Fig_Baye_qN_NoiseFree_Obj_Quad_d20}
	\end{subfigure}	
	\begin{subfigure}[t]{0.328\textwidth}
		\includegraphics[width=\textwidth]{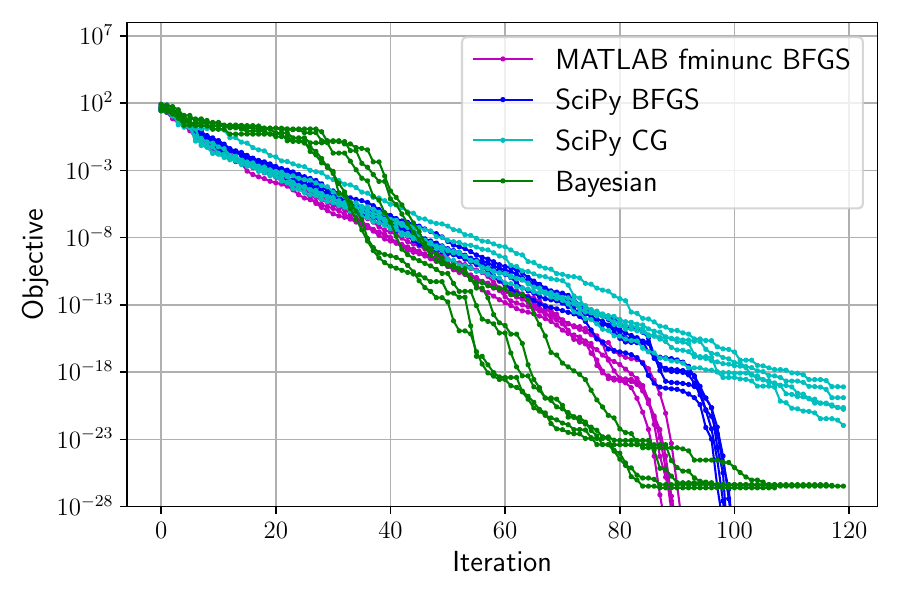}
		\caption{Objective: bowl function}
		\label{Fig_Baye_qN_NoiseFree_Obj_Bowl_d20}
	\end{subfigure}	
	\begin{subfigure}[t]{0.328\textwidth}
		\includegraphics[width=\textwidth]{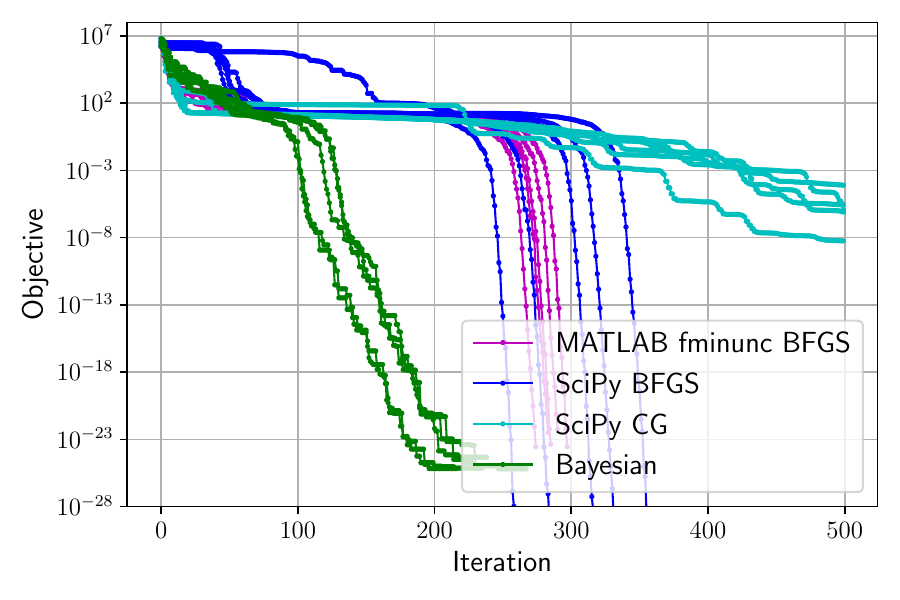}
		\caption{Objective: Rosenbrock $a=100$}
		\label{Fig_Baye_qN_NoiseFree_Obj_RosenA100_d20}
	\end{subfigure}	
	\begin{subfigure}[t]{0.328\textwidth}
		\includegraphics[width=\textwidth]{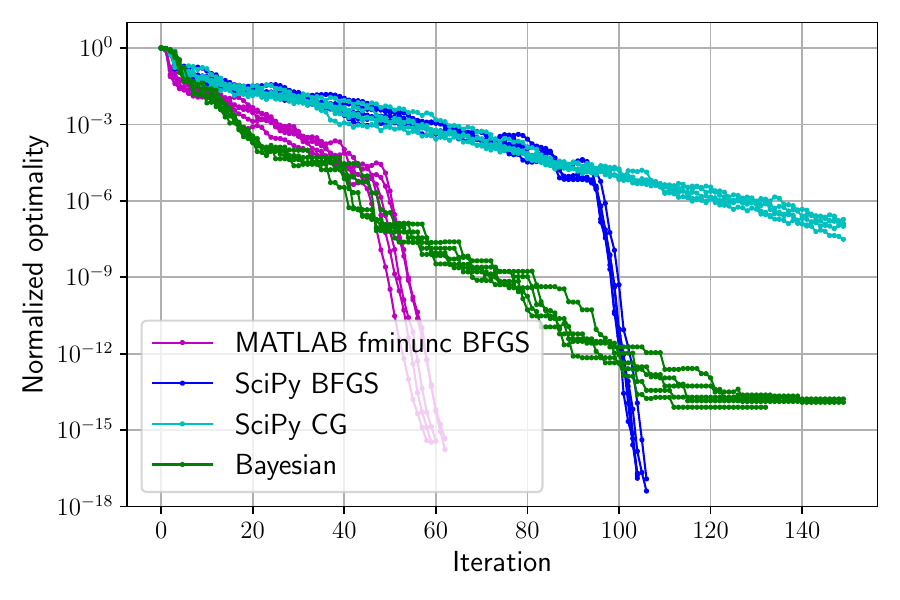}
		\caption{Optimality: quadratic function}
		\label{Fig_Baye_qN_NoiseFree_Opt_Quad_d20}
	\end{subfigure}	
	\begin{subfigure}[t]{0.328\textwidth}
		\includegraphics[width=\textwidth]{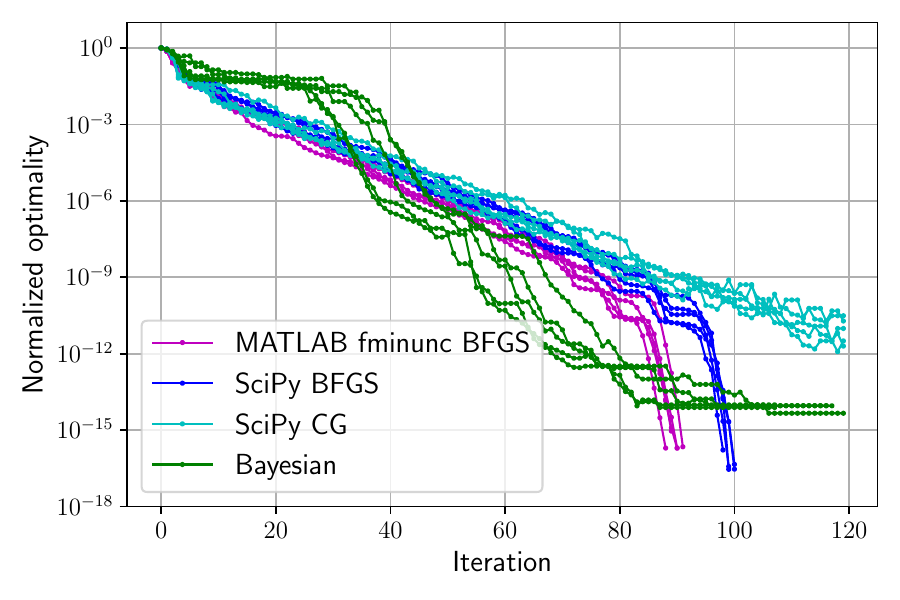}
		\caption{Optimality: bowl function}
		\label{Fig_Baye_qN_NoiseFree_Opt_Bowl_d20}
	\end{subfigure}	
	\begin{subfigure}[t]{0.328\textwidth}
		\includegraphics[width=\textwidth]{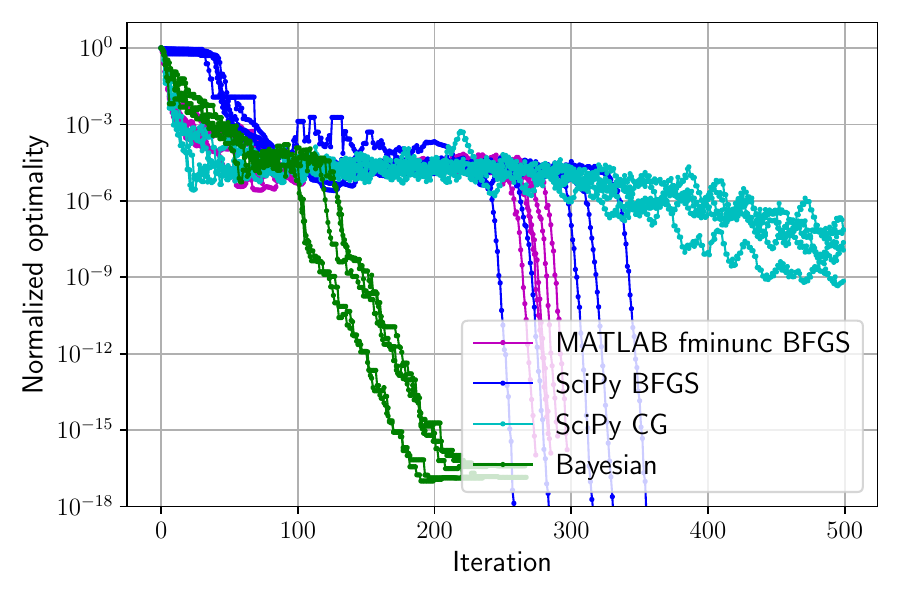}
		\caption{\mbox{Optimality: Rosenbrock $a=100$}}
		\label{Fig_Baye_qN_NoiseFree_Opt_RosenA100_d20}
	\end{subfigure}	
	\caption{Unconstrained optimization comparison of the Bayesian optimizer with the quasi-Newton and conjugate-gradient optimizers from SciPy and MATLAB. The test cases have $n_d=20$ and they are the quadratic, bowl, and Rosenbrock $a=100$ functions from \Eqss{Eq_Quadratic_fun}{Eq_Bowl}{Eq_Rosenbrock}, respectively.}
	\label{Fig_Baye_qN_NoiseFree}
\end{figure}

The comparison of the optimizers for the minimization of these functions can be seen in \Fig{Fig_Baye_qN_NoiseFree} for $n_d = 20$. In \Figs{Fig_Baye_qN_NoiseFree_Obj_Quad_d20}{Fig_Baye_qN_NoiseFree_Opt_Quad_d20} the minimization results for the quadratic test case are shown. As mentioned at the start of \Sec{Sec_LocalOptz_Studies} when the quadratic test case was introduced, quasi-Newton optimizers are well suited to minimize this test case since they approximate the objective function with a quadratic function of their own. Nonetheless, the Bayesian optimizer proves to be competitive with the quasi-Newton optimizers from SciPy and MATLAB for the quadratic test case. It is only once these quasi-Newton optimizers are near the minimum that they overtake the Bayesian optimizer. The SciPy CG optimizer is the slowest of the optimizers to converge for all three test cases.

For the bowl function, the SciPy and MATLAB optimizers initially make quicker progress in minimizing the objective than the Bayesian optimizer. However, the Bayesian and quasi-Newton optimizers achieve roughly the same minimum optimality in approximately the same number of iterations. Finally, for the Rosenbrock function, the Bayesian optimizer achieves the same final optimality as the quasi-Newton optimizers and does so in significantly fewer iterations.

\begin{figure}[t!]
	\centering
	\begin{subfigure}[t]{0.328\textwidth}
		\includegraphics[width=\textwidth]{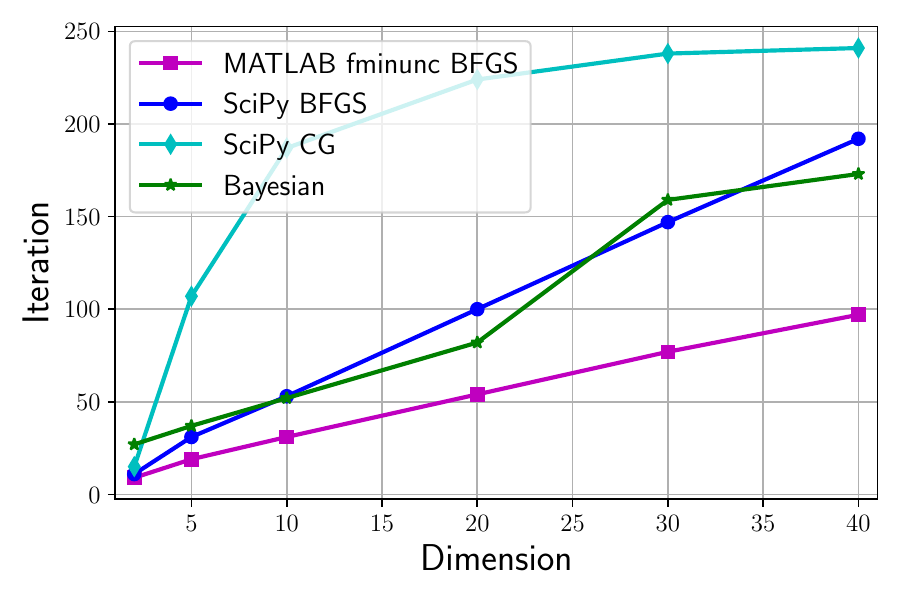}
		\caption{Quadratic function: \Eq{Eq_Quadratic_fun}}
		\label{Fig_Baye_qN_NoiseFree_OptTol_Quad}
	\end{subfigure}	
	\begin{subfigure}[t]{0.328\textwidth}
		\includegraphics[width=\textwidth]{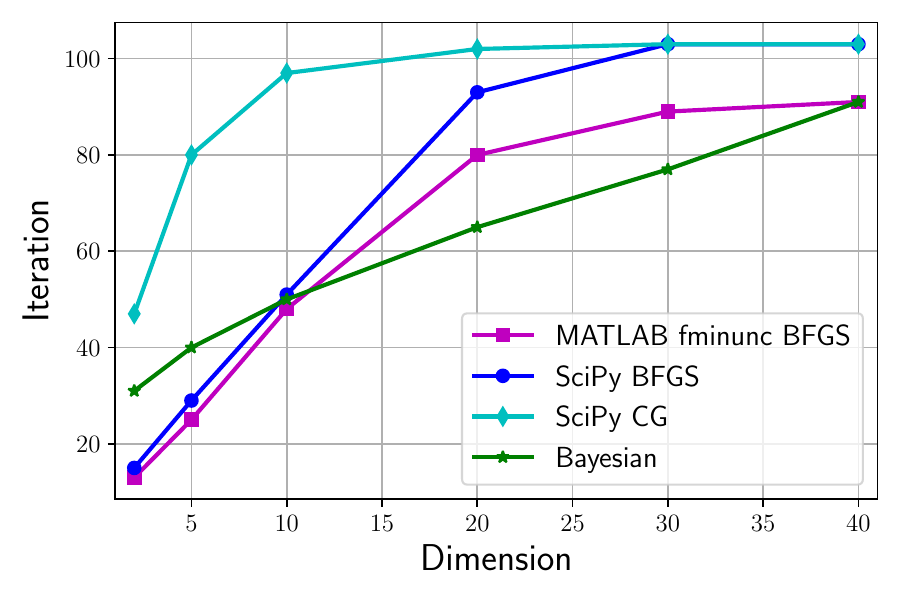}
		\caption{Bowl function: \Eq{Eq_Bowl}}
		\label{Fig_Baye_qN_NoiseFree_OptTol_Bowl}
	\end{subfigure}	
	\begin{subfigure}[t]{0.328\textwidth}
		\includegraphics[width=\textwidth]{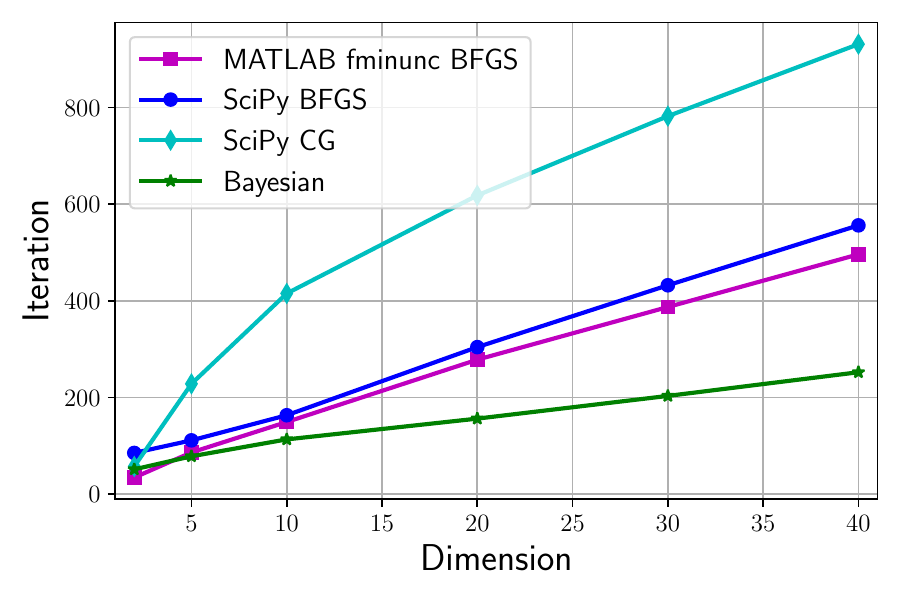}
		\caption{Rosenbrock $a=100$: \Eq{Eq_Rosenbrock}}
		\label{Fig_Baye_qN_NoiseFree_OptTol_RosenA100}
	\end{subfigure}	
	\caption{Median number of iterations for the optimizers to reduce the optimality by 10 orders of magnitude and the objective function below $10^{-5}$ for 25 independent optimization runs.}
	\label{Fig_Baye_qN_NoiseFree_OptTol}
\end{figure}

The median number of iterations required for the Bayesian and quasi-Newton optimizers to achieve a 10-order reduction in the optimality and an objective function below $10^{-5}$ can be seen in \Fig{Fig_Baye_qN_NoiseFree_OptTol}. For the quadratic test case, the performance of the Bayesian optimizer is similar to SciPy BFGS and is only outperformed by MATLAB fminunc BFGS. However, if a higher tolerance was used for the optimality, the difference in results between the Bayesian and quasi-Newton optimizers would be smaller. For the bowl function, which is similar to the quadratic function, as seen in \Fig{Fig_BoTestCases}, the Bayesian optimizer is competitive with both quasi-Newton optimizers for $n_d \geq 10$. Finally, for the Rosenbrock function shown in \Fig{Fig_Baye_qN_NoiseFree_OptTol_RosenA100}, we see that the Bayesian optimizer significantly outperforms all of the optimizers. The advantage of using the Bayesian optimizer for this test case increases significantly as the dimensionality of the problem increases. For $n_d = 40$ the Bayesian optimizer requires approximately half as many iterations as the SciPy BFGS and MATLAB fminunc optimizers to achieve the desired tolerance.

The MATLAB, and SciPy optimizers are able to reduce the objective below $10^{-5}$ and reduce the optimality 10 orders of magnitude for all values of $n_d$ that were considered for the quadratic and bowl test cases. This is also the case for Rosenbrock function with $a=1$, but not for $a=100$. The number of times that the optimizers were able to achieve the convergence criteria for the $a=100$ Rosenbrock function is indicated in \Table{Table_optz_tol_Rosenbrock_a100}. The SciPy BFGS and CG optimizers achieved the desired tolerance the most often out of all of the optimizers, around 90\% of the time. The Bayesian and MATLAB optimizers both did so around 80\% of the time. In general, all of the optimizers achieved the 10-order reduction in the optimality in all cases, but occasionally did so without converging to the minimum, \ie achieving an objective value below $10^{-5}$. In these cases, the optimizers found the valley of the Rosenbrock function, which can be seen in \Eq{Fig_BoTestCases_Rosen_a100_zoomT}, but were not able to traverse this valley to the minimum. The valley of the $a=100$ Rosenbrock function has steep walls with a nearly flat bottom, which makes minimizing this function challenging and causes optimizers to stall.

\begin{table}[t]
	\centering
	\begin{tabular}{l cccccc | c} 
		Optimizers & $n_d = 2$ & $n_d = 5$ & $n_d = 10$ & $n_d = 20$ & $n_d = 30$ & $n_d = 40$ & All cases \\
		\hline 
		MATLAB fminunc 	& 100\% & 80\% & 80\%	& 100\%	& 60\% & 60\% & 80\%\\
		SciPy BFGS 		& 100\% & 76\% & 96\% 	& 88\% 	& 96\% & 88\% & 91\% \\
		SciPy CG		& 100\% & 88\% & 84\% 	& 96\% 	& 92\% & 76\% & 89\% \\
		Bayesian 		& 100\% & 60\% & 80\%	& 72\% 	& 88\% & 72\% & 79\% \\
		\hline
	\end{tabular}
	\caption{The percentage of optimization runs that achieved an objective evaluation below $10^{-5}$ and an optimality reduction of at least 10 orders of magnitude for the Rosenbrock test case with $a=100$.}
	\label{Table_optz_tol_Rosenbrock_a100}
\end{table}


The results in this section demonstrate that quasi-Newton optimizers require fewer iterations to achieve deep convergence of the optimality for the quadratic test case, which is a problem that they are ideally suited to solve. For test cases that are bowl shaped and resemble quadratic functions, the Bayesian optimizer is competitive with quasi-Newton optimizers. For more complicated functions, such as the Rosenbrock function with $a=100$, the Bayesian optimizer was shown to be as robust as the quasi-Newton optimizers at finding the minimum while being able to do so with significantly fewer iterations and hence function evaluations, especially compared to the SciPy CG optimizer.

\section{Optimization with noisy gradients} \label{Sec_LocalOptz_NoisyGrad}

The impact of using inaccurate gradients on the performance of the Bayesian optimizer along with the quasi-Newton and conjugate-gradient
optimizers from SciPy and MATLAB is investigated in this section. Several factors can impact the accuracy of the gradients, such as using approximations, neglecting certain terms, not solving the required equations to a sufficiently small tolerance, or having noise introduced into the gradient. In this section we consider noisy gradients that are provided by
\begin{equation} \label{Eq_noisy_grad}
	\left( \gradfNoisy \right)_i = \left( \nabla f \right)_i + \mathcal{N} \left(0, \sigma^2_{\nabla f} \right) \quad \forall \, i \in \{1, \ldots, n_d \},
\end{equation}
where $\nabla f$ is the noise-free gradient and the noise on each entry of the gradient is independent, zero mean, and normally distributed with variance $\sigma^2_{\nabla f}$. The objective evaluations could also be noisy, but this is not considered in this section.

\subsection{Bayesian optimization with noisy gradients} \label{Sec_LocalOptz_NoisyGrad_BoAllStdGrad}

\begin{figure}[t!]
	\centering
	\begin{subfigure}[t]{0.328\textwidth}
		\includegraphics[width=\textwidth]{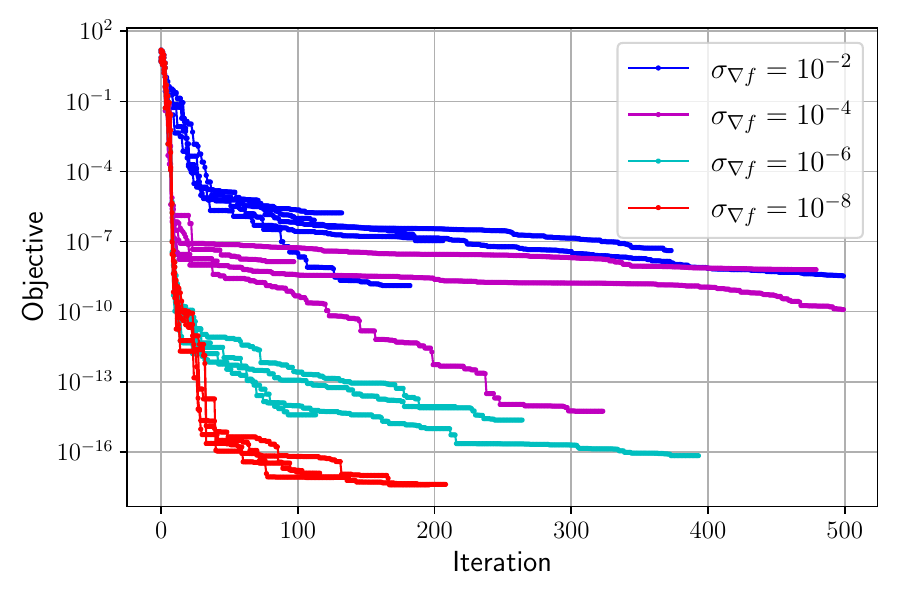}
		\caption{Objective: quadratic function}
		\label{Fig_Bo_NoisyGrad_d5_Obj_Quad}
	\end{subfigure}	
	\begin{subfigure}[t]{0.328\textwidth}
		\includegraphics[width=\textwidth]{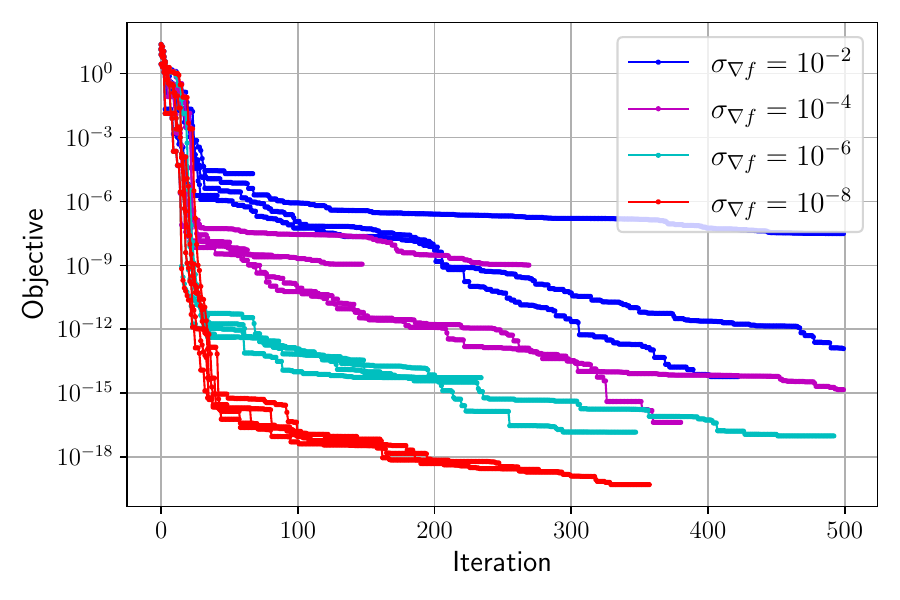}
		\caption{Objective: bowl function}
		\label{Fig_Bo_NoisyGrad_d5_Obj_Bowl}
	\end{subfigure}	
	\begin{subfigure}[t]{0.328\textwidth}
		\includegraphics[width=\textwidth]{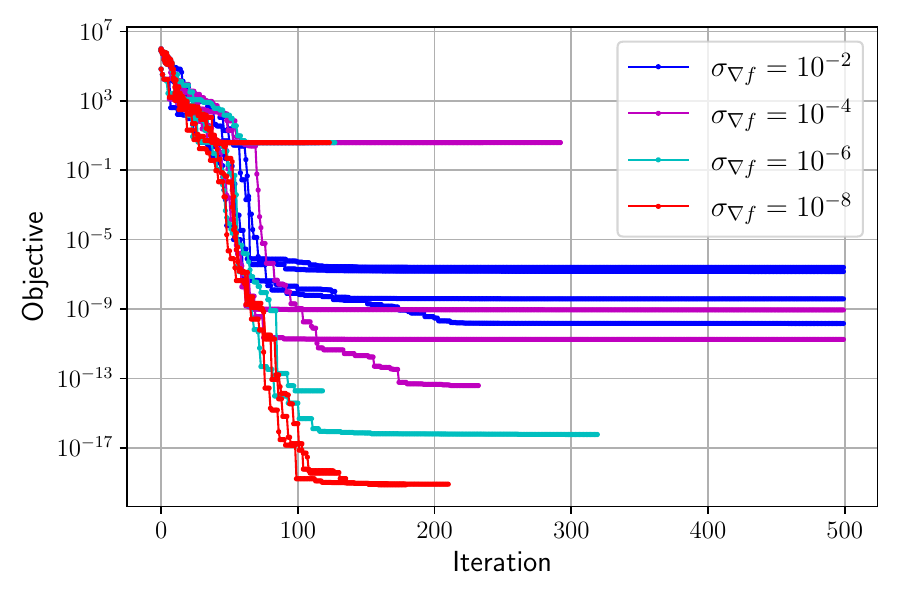}
		\caption{Objective: Rosenbrock $a=100$}
		\label{Fig_Bo_NoisyGrad_d5_Obj_RosenA100}
	\end{subfigure}	
	\begin{subfigure}[t]{0.328\textwidth}
		\includegraphics[width=\textwidth]{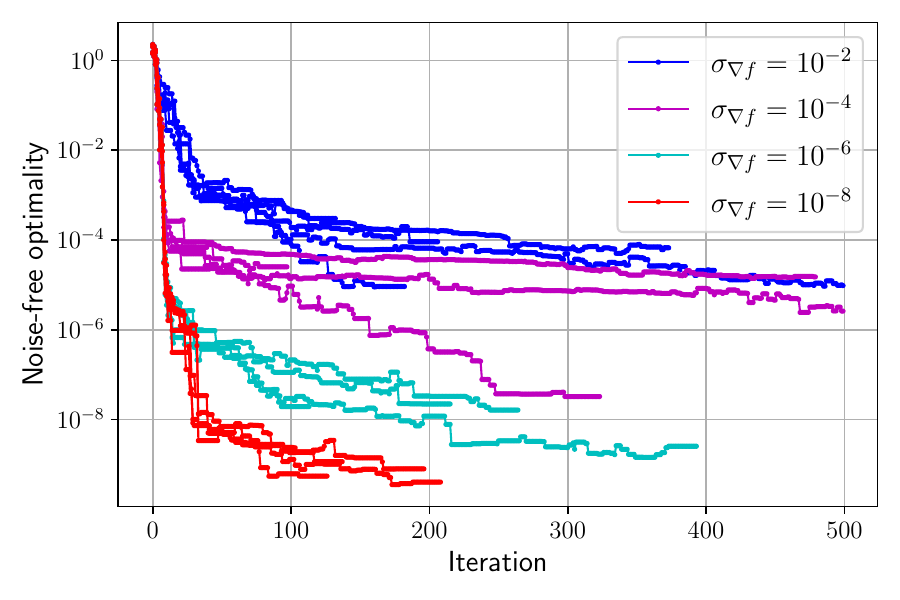}
		\caption{Optimality: quadratic function}
		\label{Fig_Bo_NoisyGrad_d5_Opt_Quad}
	\end{subfigure}	
	\begin{subfigure}[t]{0.328\textwidth}
		\includegraphics[width=\textwidth]{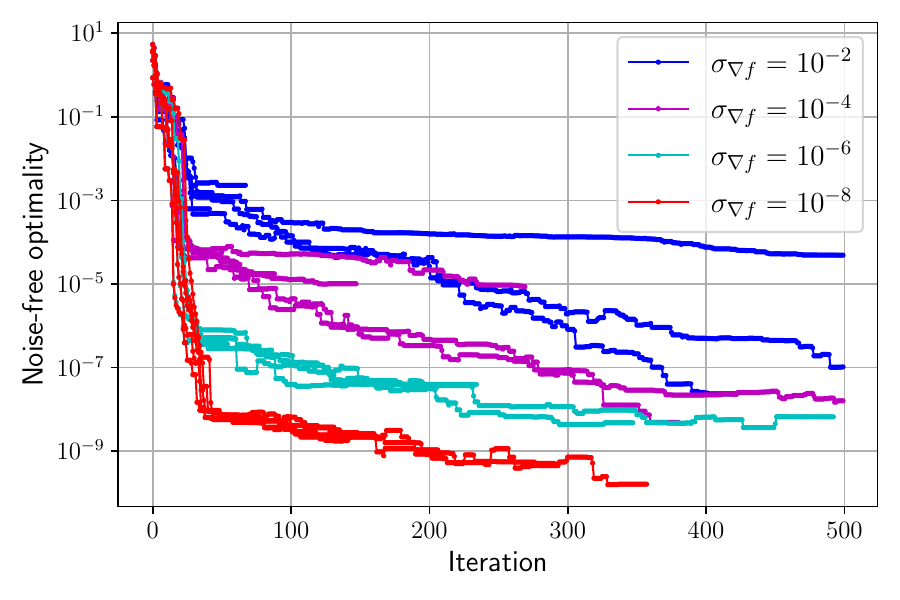}
		\caption{Optimality: bowl function}
		\label{Fig_Bo_NoisyGrad_d5_Opt_Bowl}
	\end{subfigure}	
	\begin{subfigure}[t]{0.328\textwidth}
		\includegraphics[width=\textwidth]{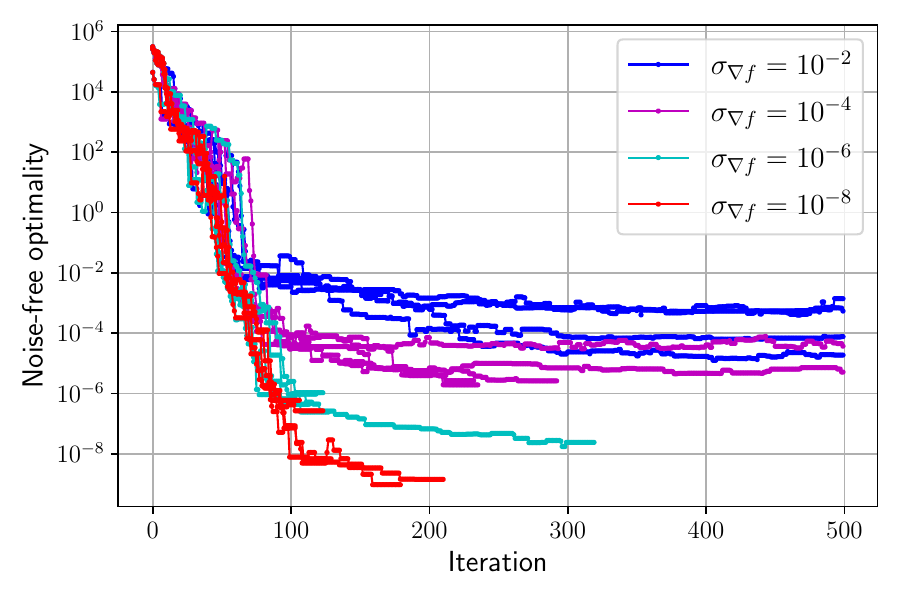}
		\caption{\mbox{Optimality: Rosenbrock $a=100$}}
		\label{Fig_Bo_NoisyGrad_d5_Opt_RosenA100}
	\end{subfigure}	
	\caption[Bayesian optimization using noisy gradients with different amounts of noise.]{Bayesian optimization using noisy gradients of the form from \Eq{Eq_noisy_grad} with different amounts of noise. The test cases have $n_d=5$ and they are the quadratic, bowl, and Rosenbrock $a=100$ functions from \Sec{Sec_LocalOptz_TestCases}. The optimality is calculated with the noise-free gradient $\nabla f$, while the optimizer only has access to the noisy gradient $\gradfNoisy$ from \Eq{Eq_noisy_grad}.}
	\label{Fig_Bo_NoisyGrad_d5}
\end{figure}

The Bayesian optimizer is used to minimize the quadratic, bowl, and Rosenbrock $a=100$ functions from \Sec{Sec_LocalOptz_TestCases}. The Bayesian optimizer is provided with noisy gradients from \Eq{Eq_noisy_grad} with $\sigma_{\nabla f} \in \{10^{-2}, 10^{-4}, 10^{-6}, 10^{-8} \}$. A flag is changed in Bayesian optimization code so that the algorithm does not assume that the gradient evaluations are accurate, as it was doing in \Secs{Sec_LocalOptz_Studies}{Sec_LocalOptz_UnconQN}. The Bayesian optimizer does not know the variance of the noise for the gradient entries $\sigma_{\nabla f}^2$, and instead estimates it with the hyperparameter $\hpstdfgrad^2$. This hyperparameter is selected by maximizing the marginal log likelihood from \Eq{Eq_ln_lkd_noisy}. \Fig{Fig_Bo_NoisyGrad_d5} shows the Bayesian optimization results using noisy gradients for the five-dimensional test cases. The top row of subfigures shows the objective functions and the bottom row plots the noise-free optimality, \ie $\| \nabla f \|_2$. We see from \Figss{Fig_Bo_NoisyGrad_d5_Obj_Quad}{Fig_Bo_NoisyGrad_d5_Obj_Bowl}{Fig_Bo_NoisyGrad_d5_Obj_RosenA100} that having noisy gradients significantly impacts how much the Bayesian optimizer is able to converge the objective function for all three test cases. \Figss{Fig_Bo_NoisyGrad_d5_Opt_Quad}{Fig_Bo_NoisyGrad_d5_Opt_Bowl}{Fig_Bo_NoisyGrad_d5_Opt_RosenA100} show that the Bayesian optimizer is able to converge the optimality only two to three orders of magnitude farther than the magnitude of $\sigma_{\nabla f}$ for these five-dimensional test cases.

\begin{figure}[t!]
	\centering
	\begin{subfigure}[t]{0.328\textwidth}
		\includegraphics[width=\textwidth]{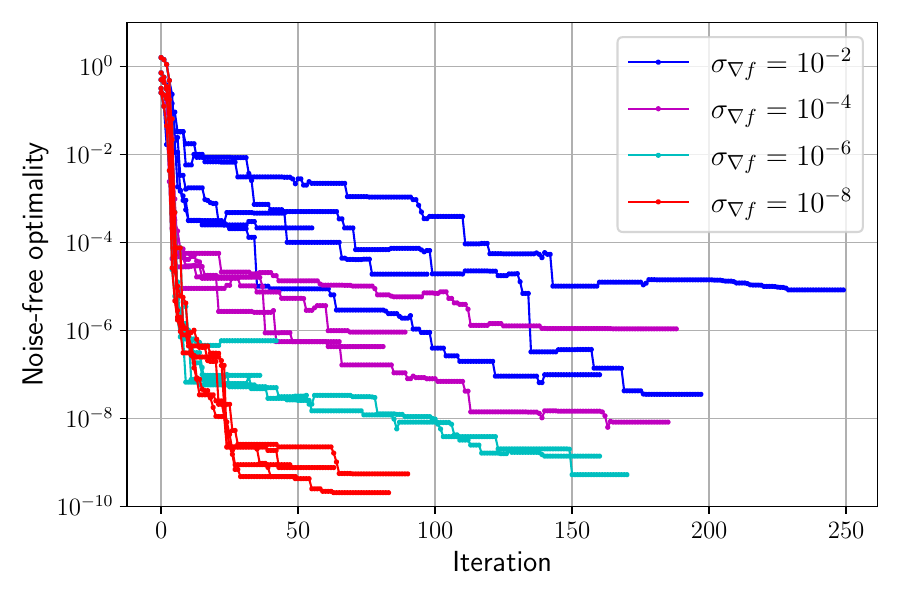}
		\caption{$n_d=2$}
		\label{Fig_Bo_NoisyGrad_ConvexOpt_d2}
	\end{subfigure}	
	\begin{subfigure}[t]{0.328\textwidth}
		\includegraphics[width=\textwidth]{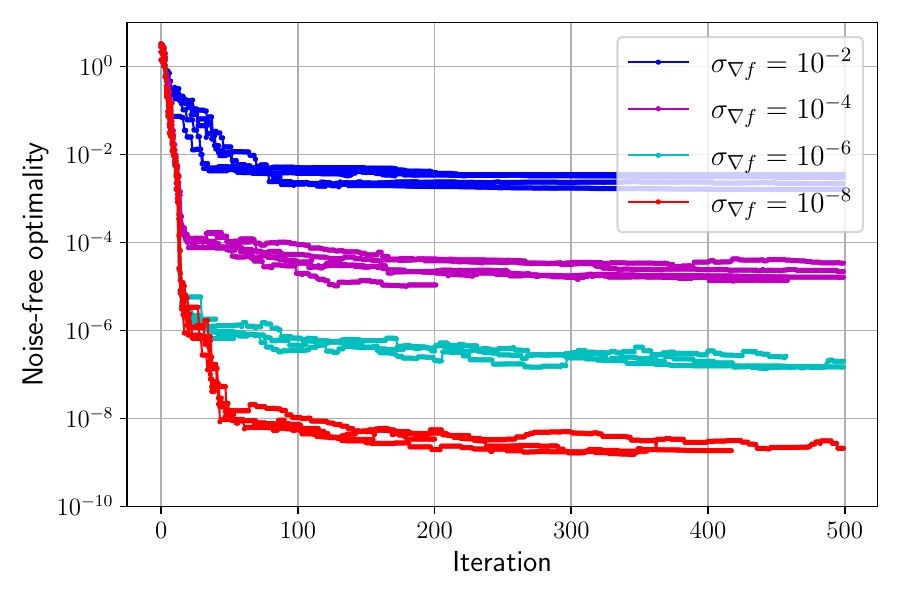}
		\caption{$n_d=10$}
		\label{Fig_Bo_NoisyGrad_ConvexOpt_d10}
	\end{subfigure}	
	\begin{subfigure}[t]{0.328\textwidth}
		\includegraphics[width=\textwidth]{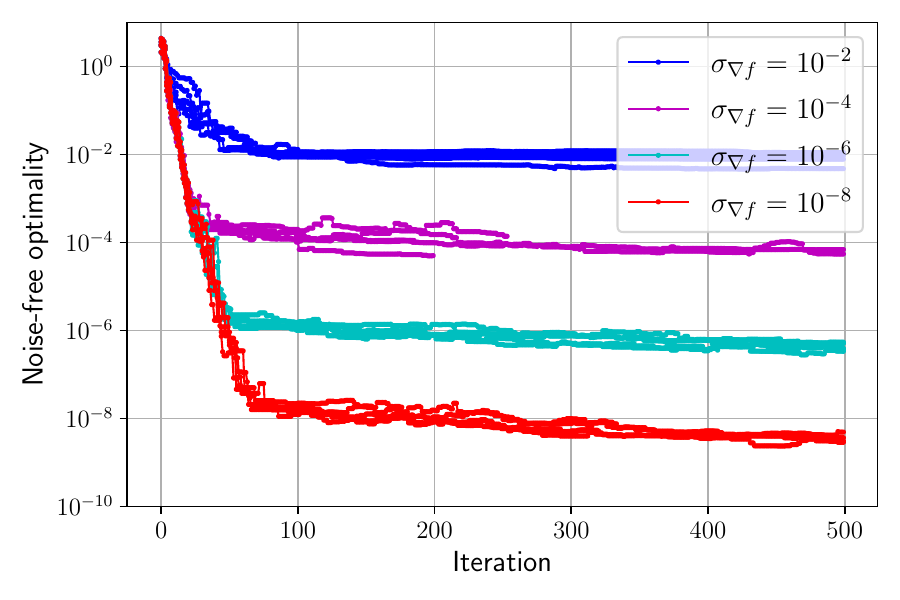}
		\caption{$n_d=20$}
		\label{Fig_Bo_NoisyGrad_ConvexOpt_d20}
	\end{subfigure}	
	\caption[Optimality comparison for the Bayesian optimization using noisy gradients with different amounts of noise.]{Comparing the noise-free optimality for the Bayesian optimizer applied to the quadratic function from \Eq{Eq_Quadratic_fun} with different amounts of noise on the gradient.}
	\label{Fig_Bo_NoisyGrad_ConvexOpt}
\end{figure}

\Fig{Fig_Bo_NoisyGrad_ConvexOpt} shows the optimality for the quadratic function from \Eq{Eq_Quadratic_fun} for $n_d \in \{ 2, 10, 20 \}$ with different amounts of noise added to the entries of the gradient. In \Fig{Fig_Bo_NoisyGrad_ConvexOpt_d2} the Bayesian optimizer is able to reduce the optimality for the $n_d=2$ quadratic test case several orders of magnitude below the magnitude of $\sigma_{\nabla f}$. In contrast, \Figs{Fig_Bo_NoisyGrad_ConvexOpt_d10}{Fig_Bo_NoisyGrad_ConvexOpt_d20} show that for $n_d = 10$ and $n_d=20$, respectively, the Bayesian optimizer is only able to reduce the optimality to approximately the same order as $\sigma_{\nabla f}$. For low-dimensional problems, the Bayesian optimizer has enough information from the noise-free objective to compensate for the noise in the gradient evaluations. However, it becomes increasingly reliant on the gradient as the dimensionality of the test cases increases. Consequently, the noise that is added to the gradients becomes a more significant limitation on how deeply the Bayesian optimizer can converge the optimality. The same results are also observed for the bowl and Rosenbrock test cases.

\begin{figure}[t!]
	\centering
	\begin{subfigure}[t]{0.328\textwidth}
		\includegraphics[width=\textwidth]{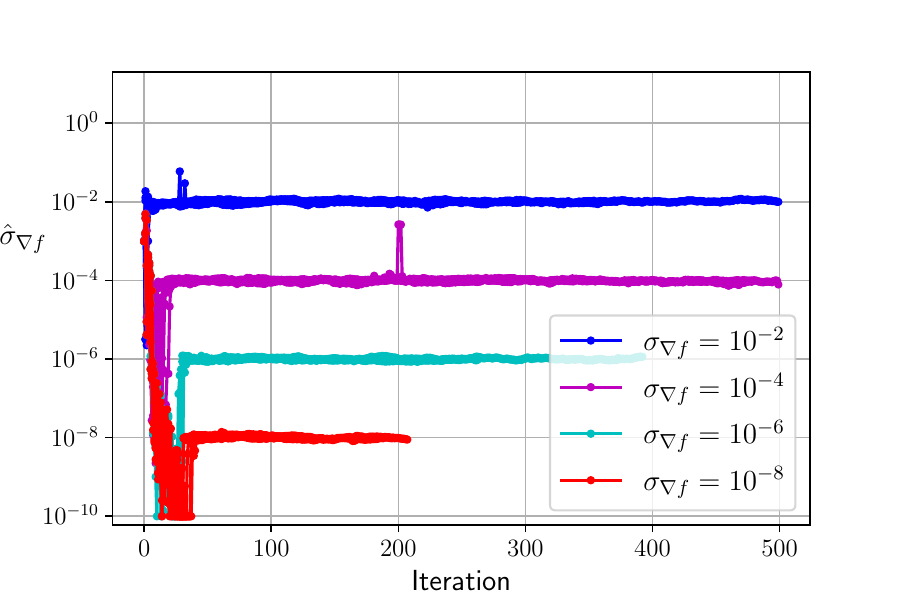}
		\caption{\Eq{Eq_Quadratic_fun}: Quadratic function}
		\label{Fig_Bo_NoisyGrad_StdGradEst_d5_Convex}
	\end{subfigure}	
	\begin{subfigure}[t]{0.328\textwidth}
		\includegraphics[width=\textwidth]{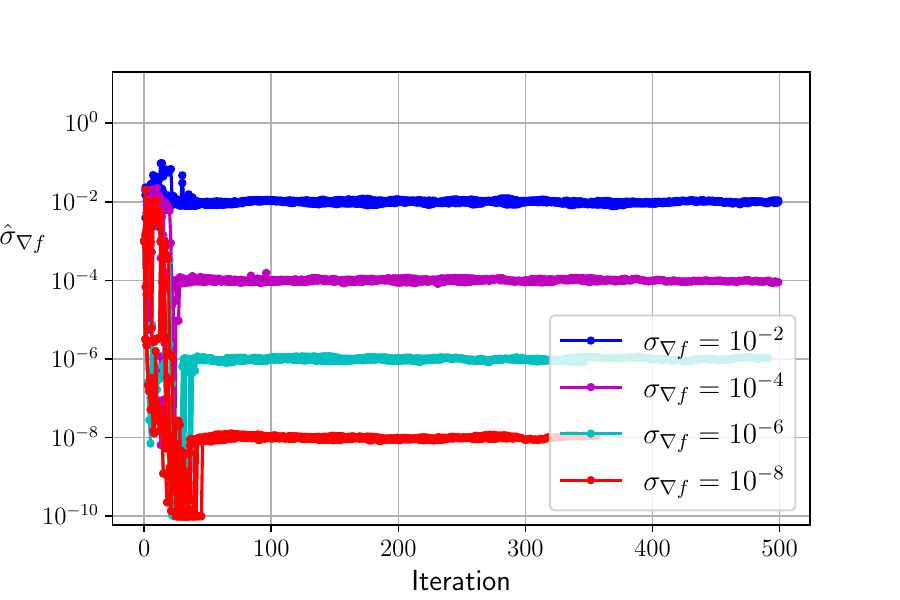}
		\caption{\Eq{Eq_Bowl}: Bowl function}
		\label{Fig_Bo_NoisyGrad_StdGradEst_d5_Bowl}
	\end{subfigure}	
	\begin{subfigure}[t]{0.328\textwidth}
		\includegraphics[width=\textwidth]{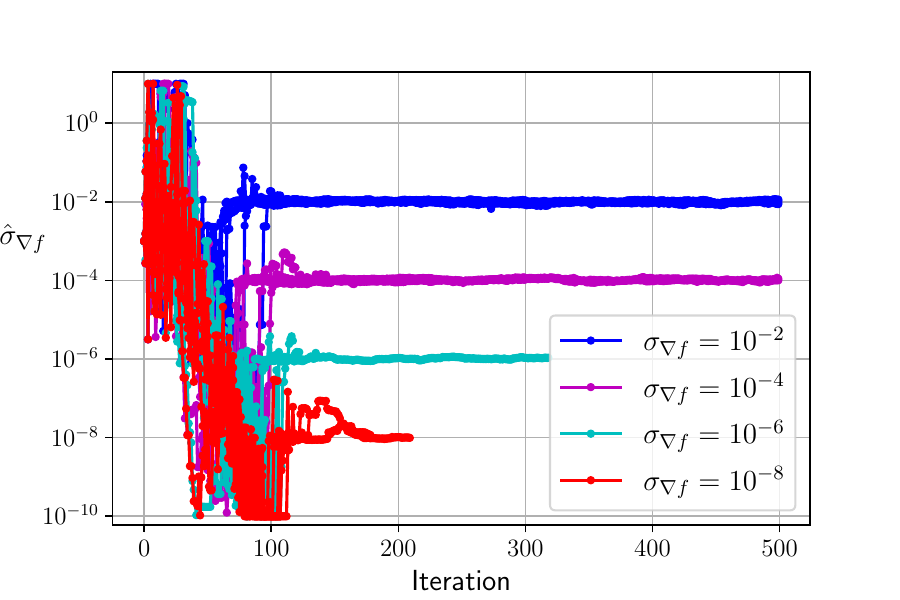}
		\caption{\Eq{Eq_Rosenbrock}: Rosenbrock $a=100$}
		\label{Fig_Bo_NoisyGrad_StdGradEst_d5_RosenA100}
	\end{subfigure}	
	\caption[Value of the GP's hyperparameter that estimates the variance of the noise for the gradient evaluations]{Value of the GP's hyperparameter $\hpstdfgrad^2$ that estimates the variance of the noise $\sigma_{\nabla f}^2$ for the gradient evaluations for three test cases with $n_d=5$.}
	\label{Fig_Bo_NoisyGrad_StdGradEst_d5}
\end{figure}

\Fig{Fig_Bo_NoisyGrad_StdGradEst_d5} shows the hyperparameter $\hpstdfgrad^2$ from the GP that estimates the variance of the noise $\sigma_{\nabla f}^2$ on the gradient evaluations at each iteration for the three unconstrained test cases with $n_d=5$. For all three test cases, the hyperparameter $\hpstdfgrad^2$ initially varies significantly between iterations. However, $\hpstdfgrad^2$ converges to $\sigma_{\nabla f}^2$ when the Bayesian optimizer gets close to the minimum. This indicates that near the minimum, the Bayesian optimizer accurately quantifies the accuracy of the gradient evaluations it is provided. This is a capability that the Bayesian optimizer has since it utilizes a probabilistic surrogate, unlike conjugate-gradient and quasi-Newton optimizers.

\subsection{Comparing Bayesian, quasi-Newton, and conjugate-gradient optimizers when gradients are noisy} \label{Sec_LocalOptz_NoisyGrad_BoVsQN}

Conjugate-gradient and quasi-Newton optimizers have several provable convergence properties but these require the function and gradient evaluations to be noise-free \cite{nocedal_numerical_2006}. Modifications to quasi-Newton optimizers have been explored that enable them to maintain their provable convergence properties for certain classes of problems, such as strongly convex functions, even when function and gradient evaluations are noisy \cite{shi_noise-tolerant_2022}.
In this subsection, the Bayesian optimizer is compared to the same previously considered unconstrained optimizers from SciPy and MATLAB for the minimization of unconstrained test cases with noisy gradients.

\begin{figure}[t!]
	\centering
	\begin{subfigure}[t]{0.328\textwidth}
		\includegraphics[width=\textwidth]{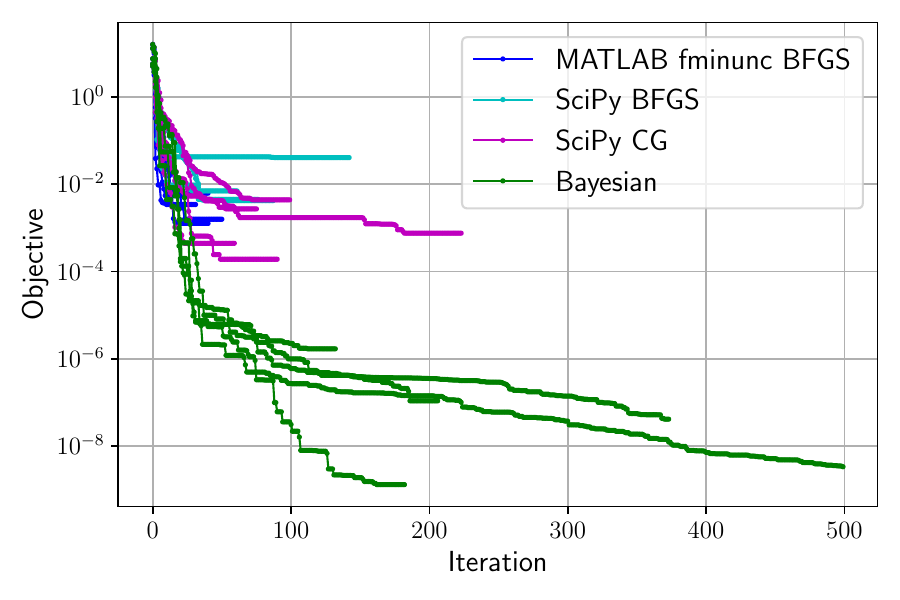}
		\caption{Objective: quadratic function}
		\label{Fig_BoVsQN_NoisyGrad_d5_Obj_Quad}
	\end{subfigure}	
	\begin{subfigure}[t]{0.328\textwidth}
		\includegraphics[width=\textwidth]{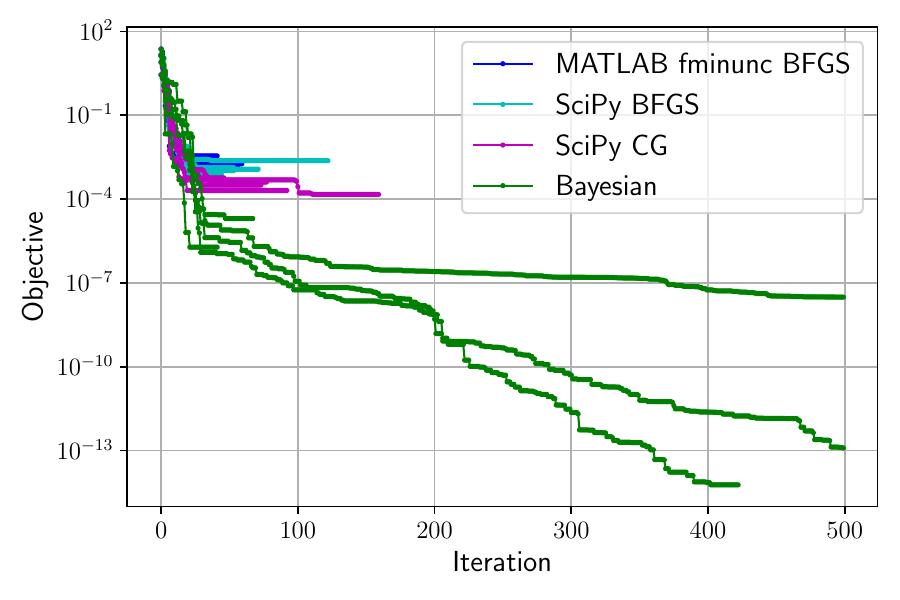}
		\caption{Objective: bowl function}
		\label{Fig_BoVsQN_NoisyGrad_d5_Obj_Bowl}
	\end{subfigure}	
	\begin{subfigure}[t]{0.328\textwidth}
		\includegraphics[width=\textwidth]{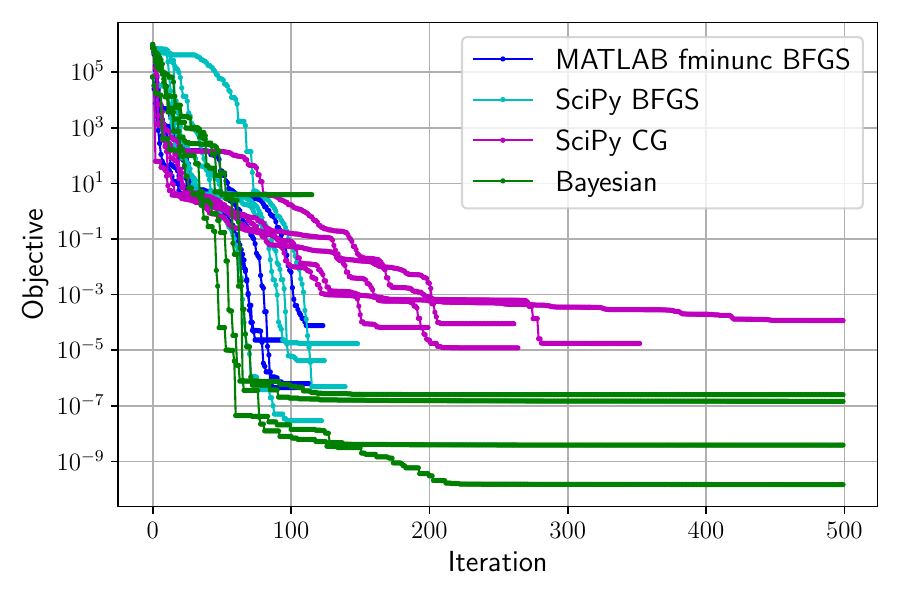}
		\caption{Objective: Rosenbrock $a=100$}
		\label{Fig_BoVsQN_NoisyGrad_d5_RosenA100}
	\end{subfigure}	
	\begin{subfigure}[t]{0.328\textwidth}
		\includegraphics[width=\textwidth]{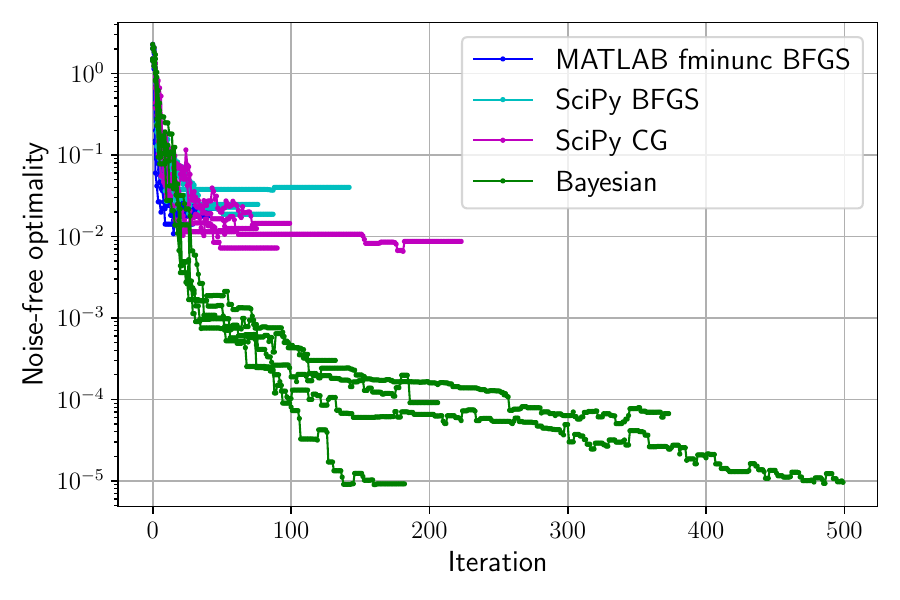}
		\caption{Optimality: quadratic function}
		\label{Fig_BoVsQN_NoisyGrad_d5_Opt_Quad}
	\end{subfigure}	
	\begin{subfigure}[t]{0.328\textwidth}
		\includegraphics[width=\textwidth]{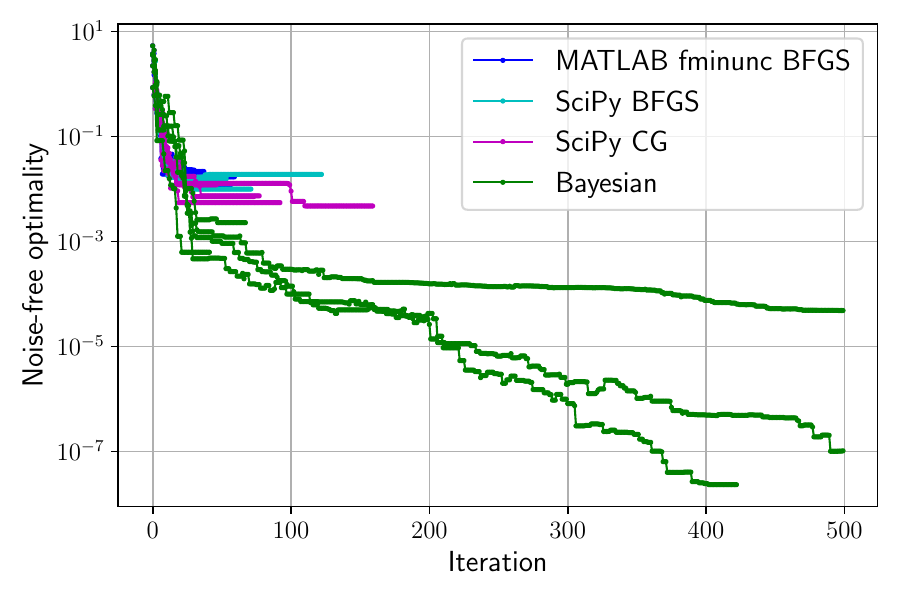}
		\caption{Optimality: bowl function}
		\label{Fig_BoVsQN_NoisyGrad_d5_Opt_Bowl}
	\end{subfigure}	
	\begin{subfigure}[t]{0.328\textwidth}
		\includegraphics[width=\textwidth]{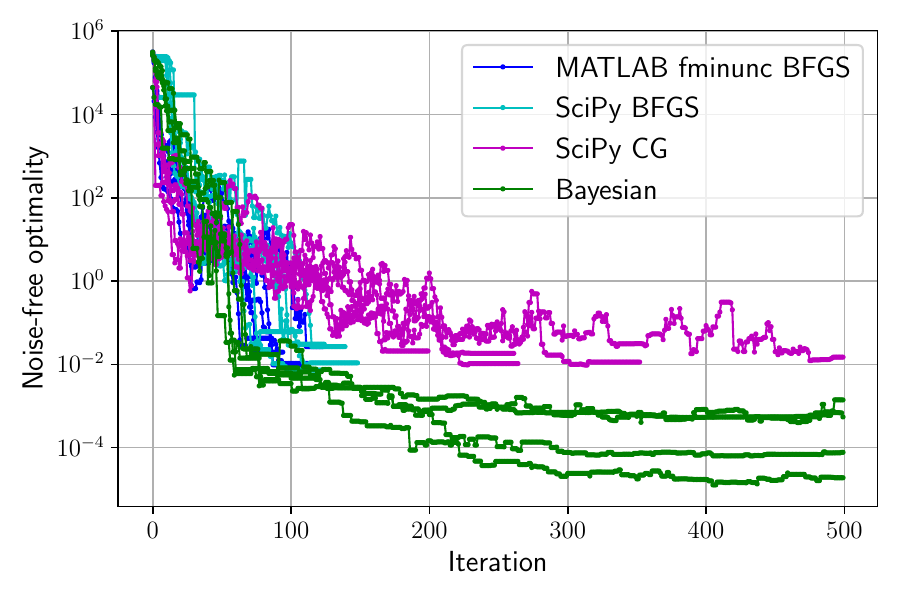}
		\caption{\mbox{Optimality: Rosenbrock $a=100$}}
		\label{Fig_BoVsQN_NoisyGrad_d5_Opt_RosenA100}
	\end{subfigure}	
	\caption{Comparison of the Bayesian optimizer to quasi-Newton and conjugate-gradient optimizers from SciPy and MATLAB for the minimization of $n_d=5$ test cases with noisy gradients. The test cases are the quadratic, bowl, and Rosenbrock $a=100$ functions from \Eqss{Eq_Quadratic_fun}{Eq_Bowl}{Eq_Rosenbrock}, respectively, and the standard deviation for the noise on the entries of the gradient is $\sigma_{\nabla f} = 10^{-2}$.}
	\label{Fig_BoVsQN_NoisyGrad_d5}
\end{figure}

\Fig{Fig_BoVsQN_NoisyGrad_d5} shows the optimization results for the $n_d=5$ quadratic, bowl, and Rosenbrock $a=100$ functions from \Eqss{Eq_Quadratic_fun}{Eq_Bowl}{Eq_Rosenbrock}, respectively, with $\hpstdfgrad = 10^{-2}$. It is clear from the top row of subfigures that the Bayesian optimizer converges the objective function for the three test cases farther than all of the MATLAB and SciPy optimizers. The bottom row of subfigures in \Fig{Fig_BoVsQN_NoisyGrad_d5} shows that the Bayesian optimizer also converges the optimality several additional orders of magnitude relative to the quasi-Newton and conjugate-gradient optimizers for all three test cases. 

In \Sec{Sec_LocalOptz_UnconQN}, the quasi-Newton
optimizers reached the desired optimization tolerance for the quadratic function in fewer iterations than the Bayesian optimizer. However, the results in \Figs{Fig_BoVsQN_NoisyGrad_d5_Obj_Quad}{Fig_BoVsQN_NoisyGrad_d5_Opt_Quad} demonstrate that when the gradients are noisy, the Bayesian optimizer is a more effective optimizer for the quadratic function than the quasi-Newton optimizers.

\begin{figure}[t!]
	\centering
	\begin{subfigure}[t]{0.328\textwidth}
		\includegraphics[width=\textwidth]{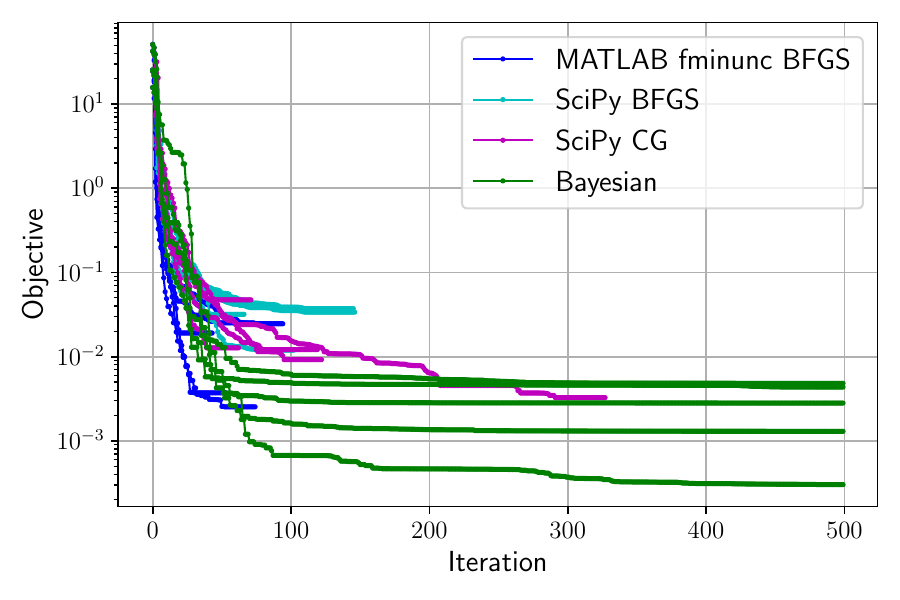}
		\caption{Objective: quadratic function}
		\label{Fig_BoVsQN_NoisyGrad_d20_Obj_Quad}
	\end{subfigure}	
	\begin{subfigure}[t]{0.328\textwidth}
		\includegraphics[width=\textwidth]{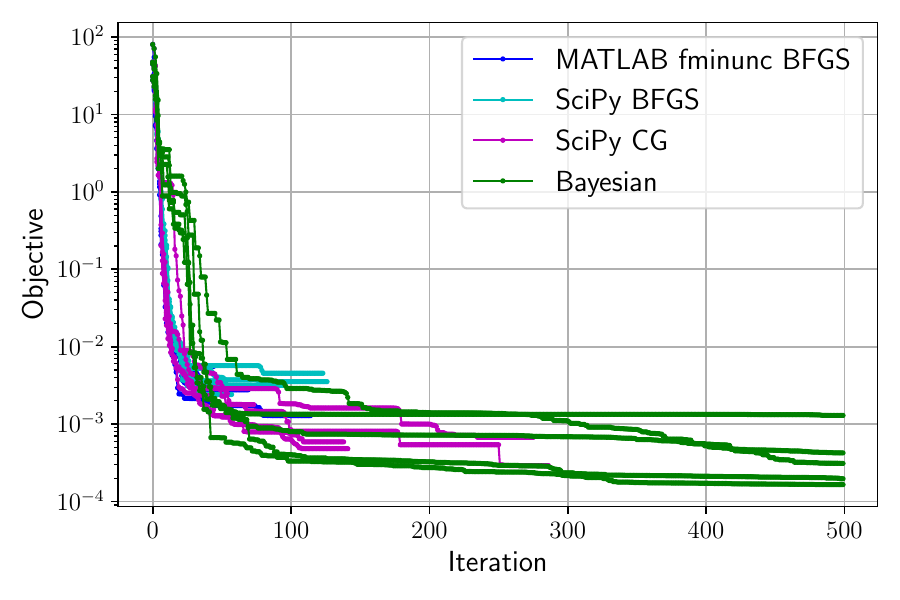}
		\caption{Objective: bowl function}
		\label{Fig_BoVsQN_NoisyGrad_d20_Obj_Bowl}
	\end{subfigure}	
	\begin{subfigure}[t]{0.328\textwidth}
		\includegraphics[width=\textwidth]{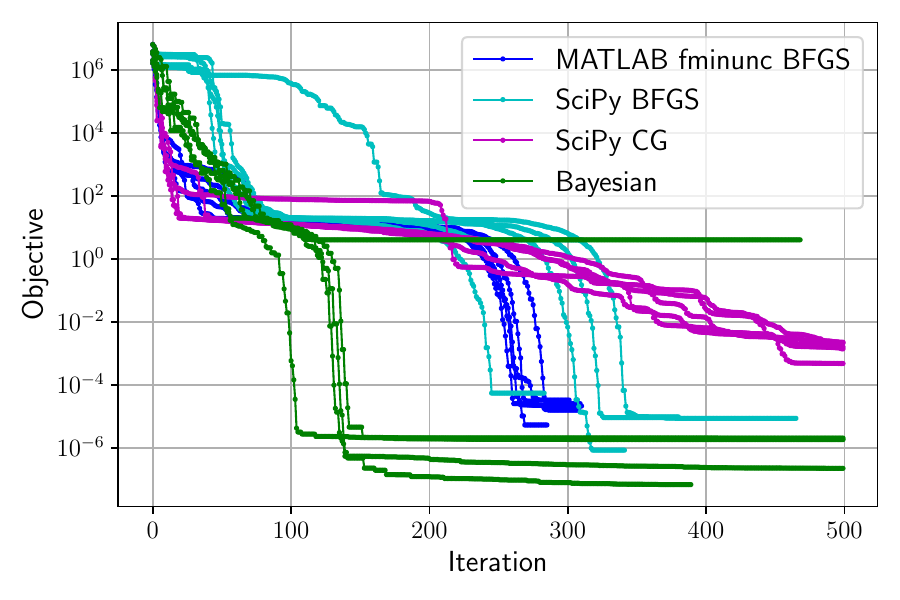}
		\caption{Objective: Rosenbrock $a=100$}
		\label{Fig_BoVsQN_NoisyGrad_d20_RosenA100}
	\end{subfigure}	
	\begin{subfigure}[t]{0.328\textwidth}
		\includegraphics[width=\textwidth]{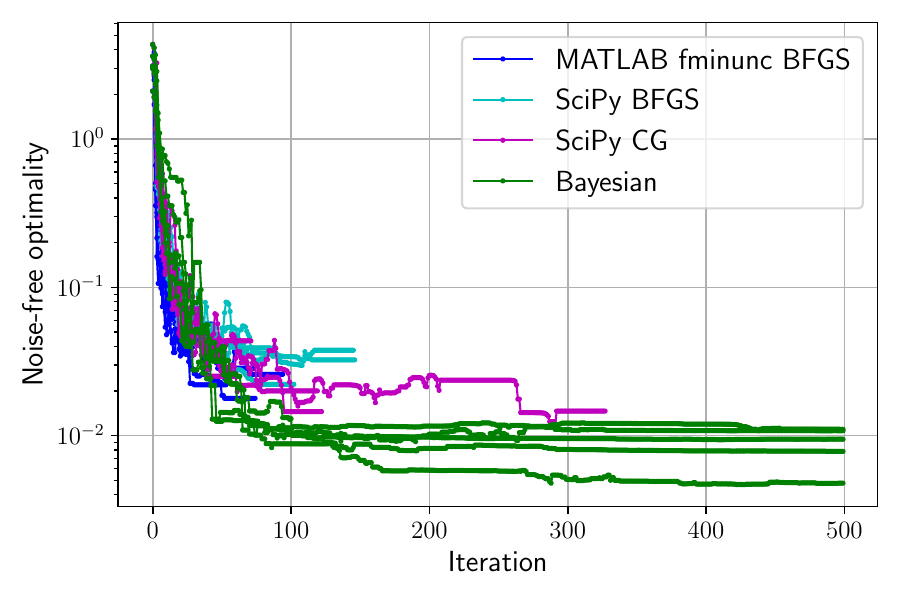}
		\caption{Optimality: quadratic function}
		\label{Fig_BoVsQN_NoisyGrad_d20_Opt_Quad}
	\end{subfigure}	
	\begin{subfigure}[t]{0.328\textwidth}
		\includegraphics[width=\textwidth]{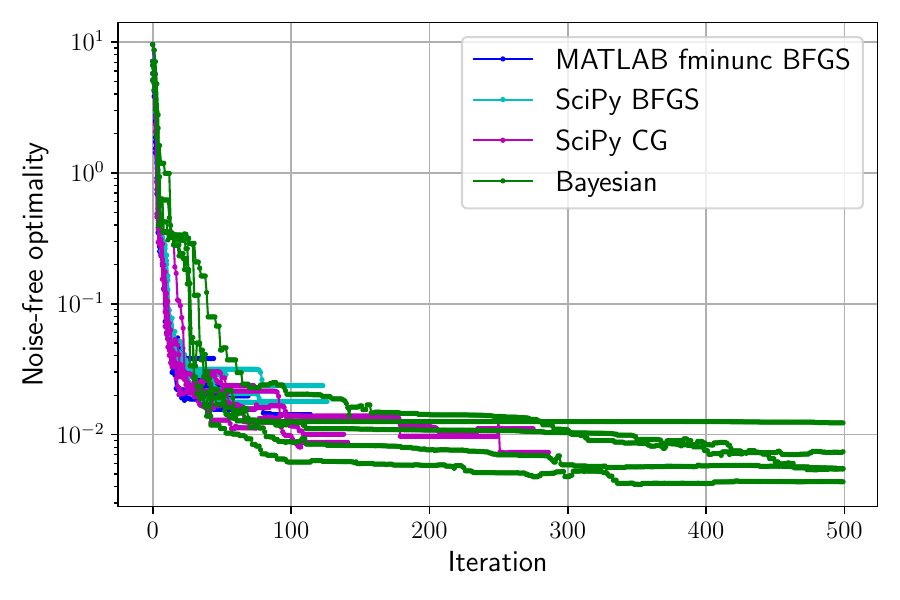}
		\caption{Optimality: bowl function}
		\label{Fig_BoVsQN_NoisyGrad_d20_Opt_Bowl}
	\end{subfigure}	
	\begin{subfigure}[t]{0.328\textwidth}
		\includegraphics[width=\textwidth]{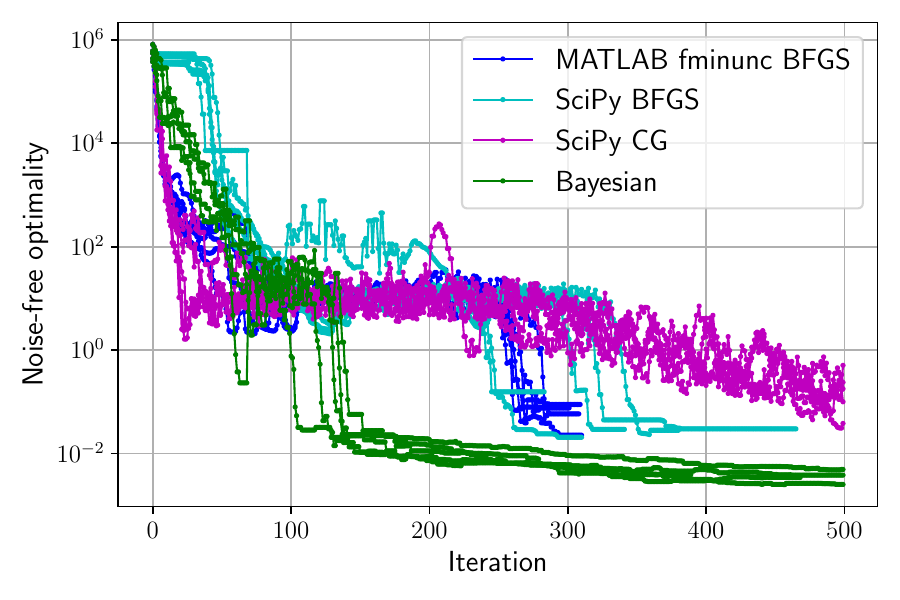}
		\caption{\mbox{Optimality: Rosenbrock $a=100$}}
		\label{Fig_BoVsQN_NoisyGrad_d20_Opt_RosenA100}
	\end{subfigure}	
	\caption{Application of the Bayesian optimizer, conjugate-gradient, and quasi-Newton optimizers from SciPy and MATLAB for unconstrained minimization using noisy gradients with $\sigma_{\nabla f} = 10^{-2}$. The test cases are the twenty-dimensional quadratic, bowl, and Rosenbrock $a=100$ functions from \Sec{Sec_LocalOptz_TestCases}.}
	\label{Fig_BoVsQN_NoisyGrad_d20}
\end{figure}

\Fig{Fig_BoVsQN_NoisyGrad_d20} shows the optimization results for the Bayesian optimizer and optimizers from SciPy and MATLAB for the three unconstrained test cases with $n_d=20$ and $\sigma_{\nabla f} = 10^{-2}$ for the noisy gradients. Once again, we can see from \Figss{Fig_BoVsQN_NoisyGrad_d20_Obj_Quad}{Fig_BoVsQN_NoisyGrad_d20_Obj_Bowl}{Fig_BoVsQN_NoisyGrad_d20_RosenA100} that the Bayesian optimizer is able to converge the objective function farther than the quasi-Newton and conjugate-gradient optimizers using the noisy gradient evaluations for all three test cases. The bottom row of subfigures in \Fig{Fig_BoVsQN_NoisyGrad_d20} shows that the final noise-free optimality for all of the optimizers is on the same order as $\sigma_{\nabla f}$. The Bayesian optimizer is able to converge the optimality further than the MATLAB and SciPy optimizers, but not as much as it was able to for the $n_d=5$ test cases that are plotted in \Fig{Fig_BoVsQN_NoisyGrad_d5}. The primary advantage of using a Bayesian optimizer for optimizations with noisy or inaccurate gradients is that they are less prone to stalling than the optimizers from SciPy and MATLAB. For optimizers that perform a line search, which includes many quasi-Newton and conjugate-gradient optimizers, they will not be able to make additional progress if the latter receives a gradient that provides it with a direction of ascent rather than descent. In contrast, Bayesian optimizers utilize a probabilistic surrogate and thus, one gradient pointing in a direction of ascent will not result in the optimizer stalling. This is demonstrated in the following section when the chaotic Lorenz 63 model is considered.

\section{Optimization of the chaotic Lorenz 63 model} \label{Sec_LocalOptz_LorenzOptz}

In this section, a linearly constrained objective function from a chaotic system is minimized. Since the test case is constrained, the SciPy trust-constr optimizer is used to compared with the Bayesian optimizer. A defining feature of chaotic systems is their extreme sensitivities to perturbations in their initial conditions \cite{strogatz_nonlinear_2018}. Furthermore, conventional analytical methods of calculating sensitivities, such as the adjoint method, break down for chaotic systems \cite{lea_sensitivity_2000,wang_forward_2013}. This presents a significant hurdle to performing optimization of systems that are both high-dimensional, which would typically necessitate the use of gradients, and computationally intensive to evaluate, which rules out using finite differences to approximate the gradients. One such example is aerodynamic shape optimization when the flow is chaotic, which can result for example from the use of large-eddy simulations \cite{blonigan_toward_2017,blonigan_least-squares_2018}. Alternative methods to calculate sensitivities for chaotic systems have been developed but these are generally substantially more expensive and less accurate when applied to chaotic systems relative to conventional sensitivity methods applied to non-chaotic systems \cite{lea_sensitivity_2000,ni_sensitivity_2017}. 

The energy method from Ashley \etal \cite{ashley_towards_2019} is used in this section since its computational cost is significantly lower than alternative methods such as least squares shadowing \cite{blonigan_toward_2017,blonigan_least-squares_2018}, and it has also been found to be more accurate than methods such as the ensemble method \cite{lea_sensitivity_2000,chandramoorthy_analysis_2017}. The energy method works by ensuring that the norm of the sensitivities remains bounded. Consider the following semi-discrete equation, which is either the formulation for an ordinary differential equation or it is a spatially discretized partial differential equation:
\begin{equation} \label{Eq_semi_discrete}
	\frac{d \uvec(t; \svec)}{dt} + \rvec_x(\uvec(t); \svec) = 0 \quad \forall \, t \in \{0, t_{\text{final}} \},
\end{equation}
where $\uvec$ is the solution to \Eq{Eq_semi_discrete}, $\rvec_x$ is the spatial residual, $\svec$ is a vector of variables, and the initial solution is $\uvec(t=0) = \uvec^{0}$. The tangent equation is now derived by differentiating \Eq{Eq_semi_discrete} with respect to a design variable $s$:
\begin{align}
	\frac{d}{dt} \frac{du}{dv} + \frac{d \rvec_x}{d \uvec} \frac{d \uvec}{ds} + \frac{d \rvec_x}{ds} 
		&= 0 \nonumber \\
	\frac{d \vvec}{dt} + \frac{d \rvec_x}{d \uvec}\vvec + \frac{d \rvec_x}{ds} 
		&= 0, \label{Eq_dvdt}
\end{align}
where $\vvec = \frac{d \uvec}{ds}$ is the tangent solution. If \Eq{Eq_dvdt} is solved for a chaotic system, then the growth in the norm $\| \vvec(t) \|$ is exponential \cite{lea_sensitivity_2000}. As explained in Appendix \ref{Sec_AppendixEnergyMethod}, the norm $\| \vvec(t) \|$ can remain bounded for all $t >0$ by clipping positive eigenvalues to zero for the Jacobian $\frac{d \rvec_x}{d \uvec}$ \cite{ashley_towards_2019}. While the tangent solution is considered here, the same concept applies to the adjoint solution as well. 

None of the current methods of approximating sensitivities for chaotic systems, including the energy method, provide machine precision for the gradients of chaotic systems using a finite computational budget \cite{lea_sensitivity_2000,wang_forward_2013,ashley_aerodynamic_2019}. As such, in order to be able to perform gradient-based optimization of chaotic systems, a numerical optimizer that can utilize inexact gradients is required. Optimizers that perform line searches, such as certain quasi-Newton optimizers, can stall if they are provided a gradient that does not point in a descent direction. Meanwhile, Bayesian optimizers can naturally utilize noisy gradients and will not stall if they are provided with a gradient that does not point in a descent direction. This will be demonstrated in this section with the gradient-based optimization of the Lorenz 1963 model:
 \begin{align*}\textit{}
 	\frac{dx}{dt} &= \sigma(y - x) \\
 	\frac{dy}{dt} &= x (\rho - z) - y \yesnumber \label{Eq_Lorenz63} \\
 	\frac{dz}{dt} &= xy - \beta z,
 \end{align*}
 where $\sigma = 10$, $\rho = 28$, and $\beta = \frac{8}{3}$ are the parameter values studied by Lorenz and using them results in the system being chaotic \cite{lorenz_deterministic_1963}. The objective function is given by
 \begin{equation} \label{Eq_Lorenz63_2dObj}
 	J(\rho, \beta) = \frac{1}{n_t} \sum_{i=1}^{n_t} \left( z_i(\rho, \beta) - 35 \right)^2 + \frac{20}{\beta},
 \end{equation}
 where $z_i$ is the third component of the solution for the Lorenz 63 model. The objective function from \Eq{Eq_Lorenz63_2dObj} was selected since it has a local minimum \cite{ashley_aerodynamic_2019}. The parameters $\rho$ and $\beta$ for the Lorenz 63 model are used as the optimization variables while $\sigma = 10$ is used for the Lorenz model. The trapezoidal time-marching method is used with the time step $\Delta t = 0.01$ and the objective is only evaluated for $t > 20$ to ensure that the solution is on the strange attractor of the system.

\begin{figure}[t!]
	\centering
	\begin{subfigure}[t]{0.49\textwidth}
		\includegraphics[width=\textwidth]{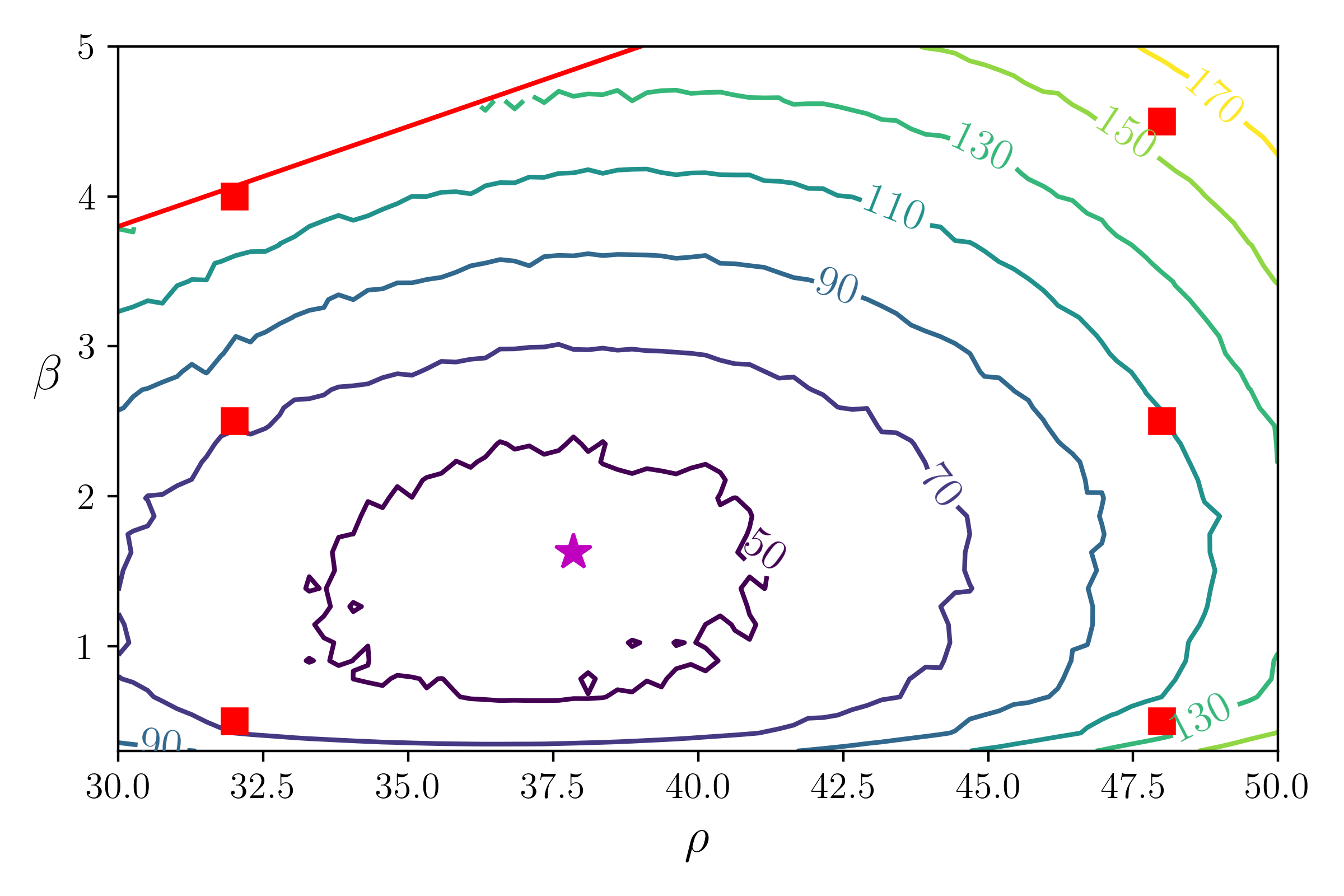}
		\caption{Contour of the objective with the red squares are the starting points of the optimizer.}
		\label{Fig_Lorenz2dOptz_Countour}
	\end{subfigure}	
	\begin{subfigure}[t]{0.49\textwidth}
		\includegraphics[width=\textwidth]{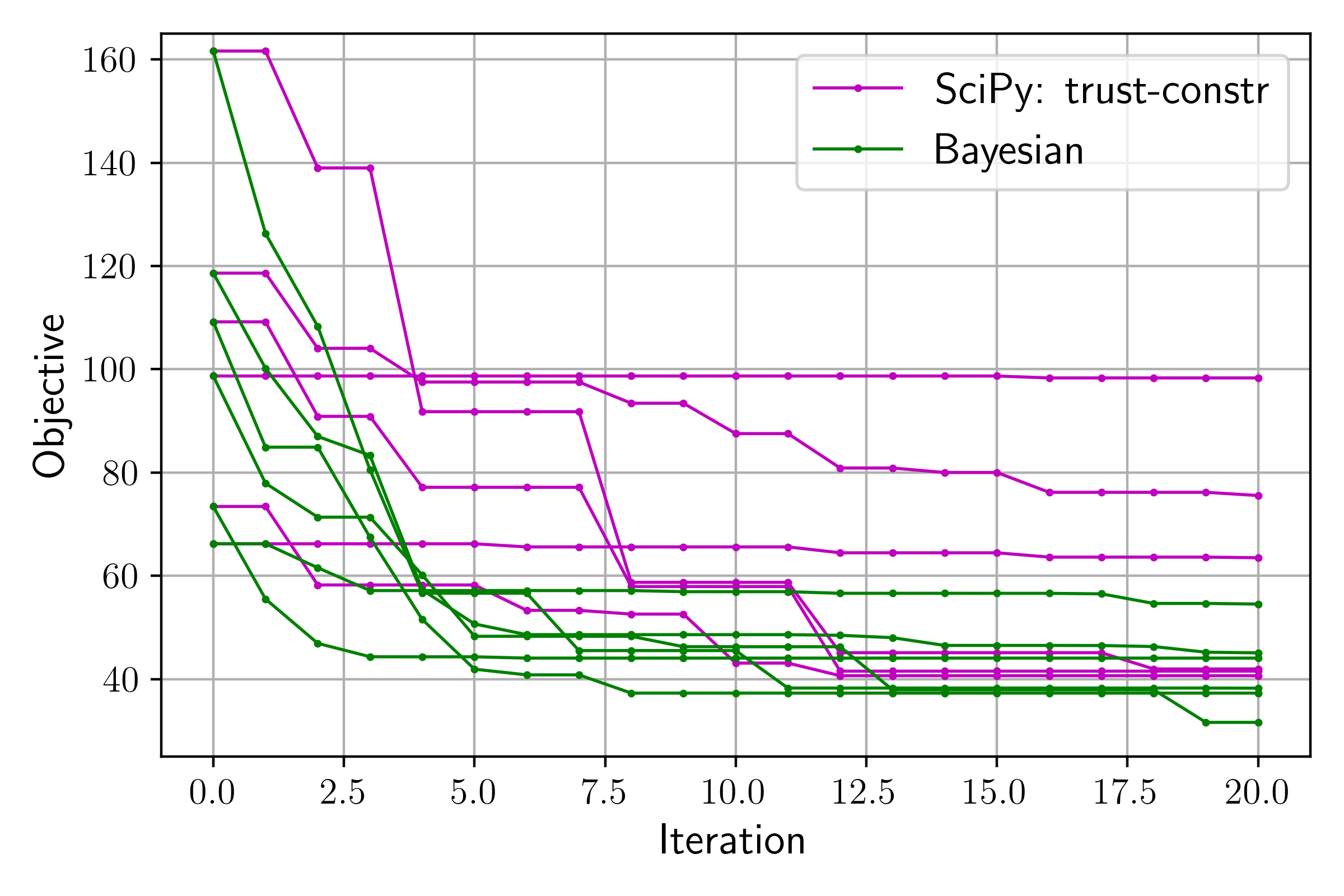}
		\caption{Optimization history using $t_J = 10$.}
		\label{Fig_Lorenz2dOptz_obj}
	\end{subfigure}	
	\caption{Gradient-based optimization of the chaotic Lorenz 63 model with the objective from \Eq{Eq_Lorenz63_2dObj}. The trapezoidal time-marching method is used with the time steps $\Delta t = 0.01$ and the objective starts being evaluated after $t_0=20$. The gradients are calculated with the energy method where the positive eigenvalues of the Jacobian are clipped to zero \cite{ashley_towards_2019}. A linear constraint, which is the red line in (a), is included such that the design space has a unique minimum.}
	\label{Fig_Lorenz2dOptz}
\end{figure}

\Fig{Fig_Lorenz2dOptz_Countour} shows a contour of the objective from \Eq{Eq_Lorenz63_2dObj}. The six red squares indicate the starting points for the optimizers and the red line is a linear constraint that was included such that there is a unique local minimum. The minimization of the objective is shown in \Fig{Fig_Lorenz2dOptz_obj} with the objective function evaluated over a period of $t_J = 10$. A maximum of 20 iterations was set for both optimizers to simulate a finite computational budget for an expensive problem. Additional iterations improve the performance of the Bayesian optimizer, which continues to make progress, while the SciPy optimizer remains stalled. The Bayesian optimizer reduces the objective function more quickly than the SciPy optimizer and achieves a lower final objective value. The SciPy optimizer often stalls at objective values significantly higher than the Bayesian optimizer. Five of the six runs for the Bayesian optimizer finished with an objective evaluation smaller than 50, which is the lowest valued contour shown in \Fig{Fig_Lorenz2dOptz_Countour}. In contrast, only half of the optimization runs for the SciPy achieved this tolerance. These results indicate that the Bayesian optimizer is able to make more effective use of the inexact gradients than the SciPy optimizer. This is despite the distribution of the error for the inexact gradients calculated with the energy method being non-Gaussian, while GPs assume that the error is Gaussian. These results indicate that the Bayesian optimizer can make effective use of inexact gradients even if the distribution of their errors is not Gaussian.

\section{Conclusions} \label{Sec_LocalOptz_Summary}

In this paper a framework was developed that enables a gradient-enhanced Bayesian optimizer to perform efficient unconstrained local optimization with both accurate and inaccurate gradient information. The settings for the Bayesian optimizer were selected through various studies that involved solving three unimodal unconstrained optimization test cases with two to forty dimensions. 

Using a subset of the function and gradient evaluations, \ie a data region, enabled the Bayesian optimizer to converge the optimality more deeply and with fewer function evaluations relative to when all the evaluation points were used. The expected improvement acquisition function was found to be effective for local minimization when combined with two trust regions. The first is a circular trust region and the second is a probabilistic trust region that limits the exploration to locations in the parameter space where the uncertainty from the probabilistic surrogate is below a set threshold. To address the severe ill-conditioning of the gradient-enhanced covariance matrix that is commonly encountered, a preconditioning method was used \cite{marchildon_solution_2024}. Finally, for smooth problems, such as the infinitely differentiable test cases considered in this paper, the use of the Gaussian, Mat\'ern $\frac{5}{2}$, and rational quadratic kernels were all found to provide similar results. The same default settings were found to be effective for the Bayesian optimizer for the test cases with accurate and inaccurate gradients.

The gradient-enhanced local Bayesian optimizer was compared with quasi-Newton and conjugate-gradient optimizers from SciPy and MATLAB. For the test cases with accurate gradients, the optimizers were compared by looking at the total number of function evaluations needed to reduce the optimality by 10 orders of magnitude. This metric was used to investigate which optimizer would be most efficient at achieving a challenging tolerance for problems with expensive function evaluations. The quasi-Newton, and conjugate-gradient, and Bayesian optimizers all reached the same final optimality. The Bayesian optimizer required a comparable number of function evaluations to reach the same tolerance as the quasi-Newton optimizers for the minimization of a quadratic test case, and significantly fewer iterations than the conjugate-gradient optimizer. However, for the Rosenbrock function the Bayesian optimizer was found to be significantly more effective. For example, the Bayesian optimizer required half as many iterations as the SciPy BFGS and MATLAB fminunc optimizers to reduce the optimality 10 orders of magnitude for the $n_d=40$ Rosenbrock function. In general, if a less stringent tolerance were used, the benefit of using the Bayesian optimizer would increase since the Bayesian optimizer is slower than the quasi-Newton optimizers at reducing the optimality near the minimum.

The second set of comparisons involved the same test cases but with zero-mean normally distributed noise added to the entries of the gradients. The Bayesian optimizer was able to converge the optimality several additional orders of magnitude relative to the quasi-Newton and conjugate-gradient optimizers, which are not able to quantify uncertainties in their inputs. Finally, the optimizers were compared for the minimization of an objective function involving the chaotic Lorenz 63 model. Conventional sensitivity methods, such as the adjoint method, break down for chaotic systems. An alternative method was used but it provided gradients with significant errors that are not normally distributed. For this test case, the Bayesian optimizer was able to consistently achieve a significantly larger reduction in the objective function relative to the quasi-Newton optimizer. 

The local Bayesian optimization framework presented in this paper enables Bayesian optimization to achieve the same deep convergence criteria that quasi-Newton and conjugate-gradient optimizers achieve and do so with the same or fewer function evaluations. The advantage of the Bayesian optimizer using a probabilistic surrogate was particularly noticeable for the test cases when the gradients were not accurate. The local optimization framework for Bayesian optimizers from this paper can also be expanded to handle nonlinear constraints \cite{marchildon_framework_2025}.

\bmsection*{Acknowledgments}

The authors are thankful for the financial support provided by the Natural Sciences and Engineering Research Council of Canada and the Ontario Graduate Scholarship Program while this research was being undertaken. The authors would also like to recognize the helpful feedback provided by professors Masayuki Yano and Prasanth Nair at the University of Toronto.


\bmsection*{Conflict of interest}

The authors declare no potential conflict of interests.

\bmsection*{ORCID}

Andr\'e Marchildon: \url{https://orcid.org/0000-0001-6407-3987}

\bibliography{MyLibrary.bib}

\begin{thebibliography}{10}
\providecommand \doibase [0]{http://dx.doi.org/}%

\bibitem{nocedal_numerical_2006}
Nocedal J, Wright SJ. {\it Numerical {Optimization}}.
\newblock Springer series in operation research and financial engineering, New
  York, NY: Springer.
\newblock second~ed., 2006.

\bibitem{jameson_optimum_1998}
Jameson A, Martinelli L, Pierce N. Optimum {Aerodynamic} {Design} {Using} the
  {Navier}-{Stokes} {Equations}. {\it Theoretical and Computational Fluid
  Dynamics.} 1998\string;10(1-4)\string:213--237.
\newblock \href {\doibase 10.1007/s001620050060} {doi: 10.1007/s001620050060}

\bibitem{rasmussen_gaussian_2006}
Rasmussen CE, Williams CKI. {\it Gaussian {Processes} for {Machine}
  {Learning}}.
\newblock Adaptive computation and machine learning, Cambridge, Mass: MIT
  Press, 2006.

\bibitem{shahriari_taking_2016}
Shahriari B, Swersky K, Wang Z, Adams RP, Freitas dN. Taking the {Human} {Out}
  of the {Loop}: {A} {Review} of {Bayesian} {Optimization}. {\it Proceedings of
  the IEEE.} 2016\string;104(1)\string:148--175.
\newblock \href {\doibase 10.1109/JPROC.2015.2494218} {doi:
  10.1109/JPROC.2015.2494218}

\bibitem{paul-dubois-taine_sensitivity-based_2013}
Paul-Dubois-Taine A, Nadarajah S. Sensitivity-{Based} {Sequential} {Sampling}
  of {Cokriging} {Response} {Surfaces} for {Aerodynamic} {Data}. In: American
  Institute of Aeronautics and Astronautics 2013; San Diego, CA

\bibitem{brochu_tutorial_2010}
Brochu E, Cora VM, Freitas dN. A {Tutorial} on {Bayesian} {Optimization} of
  {Expensive} {Cost} {Functions}, with {Application} to {Active} {User}
  {Modeling} and {Hierarchical} {Reinforcement} {Learning}.  2010\string:49.

\bibitem{zhan_expected_2020}
Zhan D, Xing H. Expected improvement for expensive optimization: a review. {\it
  Journal of Global Optimization.} 2020\string;78(3)\string:507--544.
\newblock \href {\doibase 10.1007/s10898-020-00923-x} {doi:
  10.1007/s10898-020-00923-x}

\bibitem{davidon_variable_1991}
Davidon WC. Variable {Metric} {Method} for {Minimization}. {\it SIAM Journal on
  optimization.} 1991\string;1(1)\string:1--17.

\bibitem{morris_bayesian_1993}
Morris MD, Mitchell TJ, Ylvisaker D. Bayesian {Design} and {Analysis} of
  {Computer} {Experiments}: {Use} of {Derivatives} in {Surface} {Prediction}.
  {\it Technometrics.} 1993\string;35(3)\string:243--255.
\newblock \href {\doibase 10.1080/00401706.1993.10485320} {doi:
  10.1080/00401706.1993.10485320}

\bibitem{wu_exploiting_2018}
Wu A, Aoi MC, Pillow JW. Exploiting gradients and {Hessians} in {Bayesian}
  optimization and {Bayesian} quadrature. {\it arXiv:1704.00060 [stat].} 2018.

\bibitem{cheng_gradient-enhanced_2023}
Cheng K, Zimmermann R. Gradient-{Enhanced} {Kriging} for {High}-{Dimensional}
  {Bayesian} {Optimization} with {Linear} {Embedding}. {\it AIAA Journal.}
  2023\string;61(11)\string:4946--4959.
\newblock \href {\doibase 10.2514/1.J062592} {doi: 10.2514/1.J062592}

\bibitem{marchildon_gradient-enhanced_2024}
Marchildon AL, Zingg DW. Gradient-{Enhanced} {Bayesian} {Optimization} {With}
  {Application} to {Aerodynamic} {Shape} {Optimization}. In: AIAA 2024-4405,
  American Institute of Aeronautics and Astronautics 2024; Las Vegas, Nevada

\bibitem{dalbey_efficient_2013}
Dalbey K. Efficient and robust gradient enhanced {Kriging} emulators.. Tech.
  Rep. SAND2013-7022, 1096451, Sandia National Laboratories;  2013

\bibitem{marchildon_non-intrusive_2023}
Marchildon AL, Zingg DW. A {Non}-intrusive {Solution} to the
  {Ill}-{Conditioning} {Problem} of the {Gradient}-{Enhanced} {Gaussian}
  {Covariance} {Matrix} for {Gaussian} {Processes}. {\it Journal of Scientific
  Computing.} 2023\string;95(3).
\newblock \href {\doibase 10.1007/s10915-023-02190-w} {doi:
  10.1007/s10915-023-02190-w}

\bibitem{osborne_gaussian_2009}
Osborne MA, Garnett R, Roberts SJ. Gaussian {Processes} for {Global}
  {Optimization}. In: Learning and Intelligent Optimization (LION) 2009;
  Trento, Italy.

\bibitem{march_gradient-based_2011}
March A, Willcox K, Wang Q. Gradient-based multifidelity optimisation for
  aircraft design using {Bayesian} model calibration. {\it The Aeronautical
  Journal.} 2011\string;115(1174)\string:729--738.
\newblock \href {\doibase 10.1017/S0001924000006473} {doi:
  10.1017/S0001924000006473}

\bibitem{shende_systematic_2022}
Shende S, Gillman A, Buskohl P, Vemaganti K. Systematic cost analysis of
  gradient- and anisotropy-enhanced {Bayesian} design optimization. {\it
  Structural and Multidisciplinary Optimization.} 2022\string;65(8)\string:235.
\newblock \href {\doibase 10.1007/s00158-022-03324-8} {doi:
  10.1007/s00158-022-03324-8}

\bibitem{marchildon_solution_2024}
Marchildon AL, Zingg DW. A solution to the ill‐conditioning of
  gradient‐enhanced covariance matrices for {Gaussian} processes. {\it
  International Journal for Numerical Methods in Engineering.} 2024.
\newblock \href {\doibase 10.1002/nme.7498} {doi: 10.1002/nme.7498}

\bibitem{mortished_aircraft_2016}
Mortished C, Ollar J, Toropov V, Sienz J. Aircraft {Wing} {Optimization} based
  on {Computationally} {Efficient} {Gradient}-{Enhanced} {Ordinary} {Kriging}
  {Metamodel} {Building}. In:  2016; San Diego, California, USA

\bibitem{martins_engineering_2021}
Martins JRRA, Ning A. {\it Engineering {Design} {Optimization}}.
\newblock Cambridge University Press.
\newblock 1~ed., 2021.

\bibitem{marchildon_framework_2025}
Marchildon AL, Zingg DW. A {Framework} for {Nonlinearly}-{Constrained}
  {Gradient}-{Enhanced} {Local} {Bayesian} {Optimization} with {Comparisons} to
  {Quasi}-{Newton} {Optimizers}. 2025.
\newblock arXiv:2506.00648 [math]

\bibitem{zhang_exploiting_2005}
Zhang Y, Leithead WE. Exploiting {Hessian} matrix and trust-region algorithm in
  hyperparameters estimation of {Gaussian} process. {\it Appl. Math. Comput..}
  2005.

\bibitem{svensson_marginalizing_2015}
Svensson A, Dahlin J, Schon TB. Marginalizing {Gaussian} process
  hyperparameters using sequential {Monte} {Carlo}. In: IEEE 2015; Cancun,
  Mexico\string:477--480

\bibitem{ollar_gradient_2017}
Ollar J, Mortished C, Jones R, Sienz J, Toropov V. Gradient based
  hyper-parameter optimisation for well conditioned kriging metamodels. {\it
  Structural and Multidisciplinary Optimization.}
  2017\string;55(6)\string:2029--2044.
\newblock \href {\doibase 10.1007/s00158-016-1626-8} {doi:
  10.1007/s00158-016-1626-8}

\bibitem{chen_exploiting_2022}
Chen L, Qiu H, Gao L, Yang Z, Xu D. Exploiting active subspaces of
  hyperparameters for efficient high-dimensional {Kriging} modeling. {\it
  Mechanical Systems and Signal Processing.} 2022\string;169.
\newblock \href {\doibase 10.1016/j.ymssp.2021.108643} {doi:
  10.1016/j.ymssp.2021.108643}

\bibitem{chen_optimization_2020}
Chen L, Qiu H, Gao L, Jiang C, Yang Z. Optimization of expensive black-box
  problems via {Gradient}-enhanced {Kriging}. {\it Computer Methods in Applied
  Mechanics and Engineering.} 2020\string;362.
\newblock \href {\doibase 10.1016/j.cma.2020.112861} {doi:
  10.1016/j.cma.2020.112861}

\bibitem{han_improving_2013}
Han ZH, Görtz S, Zimmermann R. Improving variable-fidelity surrogate modeling
  via gradient-enhanced kriging and a generalized hybrid bridge function. {\it
  Aerospace Science and Technology.} 2013\string;25(1)\string:177--189.
\newblock \href {\doibase 10.1016/j.ast.2012.01.006} {doi:
  10.1016/j.ast.2012.01.006}

\bibitem{wu_bayesian_2017}
Wu J, Poloczek M, Wilson AG, Frazier P. Bayesian {Optimization} with
  {Gradients}. In:  2017; Long Beach, CA, USA\string:5273--5284.

\bibitem{zimmermann_maximum_2013}
Zimmermann R. On the {Maximum} {Likelihood} {Training} of {Gradient}-{Enhanced}
  {Spatial} {Gaussian} {Processes}. {\it SIAM Journal on Scientific Computing.}
  2013\string;35(6)\string:A2554--A2574.
\newblock \href {\doibase 10.1137/13092229X} {doi: 10.1137/13092229X}

\bibitem{laurent_overview_2019}
Laurent L, Le~Riche R, Soulier B, Boucard PA. An {Overview} of
  {Gradient}-{Enhanced} {Metamodels} with {Applications}. {\it Archives of
  Computational Methods in Engineering.} 2019\string;26(1)\string:61--106.
\newblock \href {\doibase 10.1007/s11831-017-9226-3} {doi:
  10.1007/s11831-017-9226-3}

\bibitem{de_roos_high-dimensional_2021}
De~Roos F, Gessner A, Hennig P. High-{Dimensional} {Gaussian} {Process}
  {Inference} with {Derivatives}. In:  2021\string:2535--2545.

\bibitem{toal_kriging_2008}
Toal DJJ, Bressloff NW, Keane AJ. Kriging {Hyperparameter} {Tuning}
  {Strategies}. {\it AIAA Journal.} 2008\string;46(5)\string:1240--1252.
\newblock \href {\doibase 10.2514/1.34822} {doi: 10.2514/1.34822}

\bibitem{toal_development_2011}
Toal DJ, Bressloff NW, Keane AJ, Holden CM. The development of a hybridized
  particle swarm for kriging hyperparameter tuning. {\it Engineering
  Optimization.} 2011\string;43(6)\string:675--699.
\newblock \href {\doibase 10.1080/0305215X.2010.508524} {doi:
  10.1080/0305215X.2010.508524}

\bibitem{snoek_practical_2012}
Snoek J, Larochelle H, Adams RP. Practical {Bayesian} {Optimization} of
  {Machine} {Learning} {Algorithms}. {\it International Conference on Neural
  Information Processing Systems.} 2012\string;2\string:2951--2959.

\bibitem{bouhlel_improved_2016}
Bouhlel MA, Bartoli N, Otsmane A, Morlier J. An {Improved} {Approach} for
  {Estimating} the {Hyperparameters} of the {Kriging} {Model} for
  {High}-{Dimensional} {Problems} through the {Partial} {Least} {Squares}
  {Method}. {\it Mathematical Problems in Engineering.} 2016\string;2016.
\newblock \href {\doibase 10.1155/2016/6723410} {doi: 10.1155/2016/6723410}

\bibitem{amine_bouhlel_efficient_2018}
Amine~Bouhlel M, Bartoli N, Regis RG, Otsmane A, Morlier J. Efficient global
  optimization for high-dimensional constrained problems by using the {Kriging}
  models combined with the partial least squares method. {\it Engineering
  Optimization.} 2018\string;50(12)\string:2038--2053.
\newblock \href {\doibase 10.1080/0305215X.2017.1419344} {doi:
  10.1080/0305215X.2017.1419344}

\bibitem{zhao_efficient_2020}
Zhao L, Wang P, Song B, Wang X, Dong H. An efficient kriging modeling method
  for high-dimensional design problems based on maximal information
  coefficient. {\it Structural and Multidisciplinary Optimization.}
  2020\string;61(1)\string:39--57.
\newblock \href {\doibase 10.1007/s00158-019-02342-3} {doi:
  10.1007/s00158-019-02342-3}

\bibitem{cheng_sliced_2023}
Cheng K, Zimmermann R. Sliced {Gradient}-{Enhanced} {Kriging} for
  {High}-{Dimensional} {Function} {Approximation}. {\it SIAM Journal on
  Scientific Computing.} 2023\string;45(6)\string:A2858--A2885.
\newblock \href {\doibase 10.1137/22M154315X} {doi: 10.1137/22M154315X}

\bibitem{chung_using_2002}
Chung HS, Alonso J. Using gradients to construct cokriging approximation models
  for high-dimensional design optimization problems. In:  2002; Reno, NV,
  U.S.A.

\bibitem{jones_efficient_1998}
Jones DR, Schonlau M, Welch WJ. Efficient {Global} {Optimization} of
  {Expensive} {Black}-{Box} {Functions}. {\it Journal of Global Optimization.}
  1998\string;13\string:455--492.
\newblock \href {\doibase https://doi.org/10.1023/A:1008306431147} {doi:
  https://doi.org/10.1023/A:1008306431147}

\bibitem{bernardo_optimization_2011}
Gramacy RB, Lee HKH. Optimization {Under} {Unknown} {Constraints}. In:
  Bernardo JM, Bayarri MJ, Berger JO, et al. \kern-2pt, eds. {\it Bayesian
  {Statistics} 9}, , Oxford University Press,  2011\string:229--256

\bibitem{virtanen_scipy_2020}
Virtanen P, Gommers R, Oliphant TE, et al. {SciPy} 1.0: fundamental algorithms
  for scientific computing in {Python}. {\it Nature Methods.}
  2020\string;17(3)\string:261--272.
\newblock \href {\doibase 10.1038/s41592-019-0686-2} {doi:
  10.1038/s41592-019-0686-2}

\bibitem{broyden_convergence_1970}
Broyden CG. The {Convergence} of a {Class} of {Double}-rank {Minimization}
  {Algorithms} 1. {General} {Considerations}. {\it IMA Journal of Applied
  Mathematics.} 1970\string;6(1)\string:76--90.
\newblock \href {\doibase 10.1093/imamat/6.1.76} {doi: 10.1093/imamat/6.1.76}

\bibitem{fletcher_new_1970}
Fletcher R. A new approach to variable metric algorithms. {\it The Computer
  Journal.} 1970\string;13(3)\string:317--322.
\newblock \href {\doibase 10.1093/comjnl/13.3.317} {doi:
  10.1093/comjnl/13.3.317}

\bibitem{goldfarb_family_1970}
Goldfarb D. A {Family} of {Variable}-{Metric} {Methods} {Derived} by
  {Variational} {Means}. {\it Mathematics of computation.}
  1970\string;24(109)\string:23--26.

\bibitem{shanno_conditioning_1970}
Shanno DF. Conditioning of {Quasi}-{Newton} {Methods} for {Function}
  {Minimization}. {\it Mathematics of Computation.}
  1970\string;24(111)\string:647--656.

\bibitem{shi_noise-tolerant_2022}
Shi HJM, Xie Y, Byrd R, Nocedal J. A {Noise}-{Tolerant} {Quasi}-{Newton}
  {Algorithm} for {Unconstrained} {Optimization}. {\it SIAM Journal on
  Optimization.} 2022\string;32(1)\string:29--55.
\newblock \href {\doibase 10.1137/20M1373190} {doi: 10.1137/20M1373190}

\bibitem{strogatz_nonlinear_2018}
Strogatz SH. {\it Nonlinear {Dynamics} and {Chaos}}.
\newblock CRC Press.
\newblock 2~ed., 2018

\bibitem{lea_sensitivity_2000}
Lea DJ, Allen MR, Haine TWN. Sensitivity analysis of the climate of a chaotic
  system. {\it Tellus A: Dynamic Meteorology and Oceanography.}
  2000\string;52(5)\string:523--532.
\newblock \href {\doibase 10.1034/j.1600-0870.2000.01137.x} {doi:
  10.1034/j.1600-0870.2000.01137.x}

\bibitem{wang_forward_2013}
Wang Q. Forward and adjoint sensitivity computation of chaotic dynamical
  systems. {\it Journal of Computational Physics.}
  2013\string;235\string:1--13.
\newblock \href {\doibase 10.1016/j.jcp.2012.09.007} {doi:
  10.1016/j.jcp.2012.09.007}

\bibitem{blonigan_toward_2017}
Blonigan PJ, Fernandez P, Murman SM, Wang Q, Rigas G, Magri L. Toward a chaotic
  adjoint for {LES}. {\it arXiv:1702.06809 [nlin, physics:physics].} 2017.
\newblock arXiv: 1702.06809.

\bibitem{blonigan_least-squares_2018}
Blonigan PJ, Wang Q, Nielsen EJ, Diskin B. Least-{Squares} {Shadowing}
  {Sensitivity} {Analysis} of {Chaotic} {Flow} {Around} a {Two}-{Dimensional}
  {Airfoil}. {\it AIAA Journal.} 2018\string;56(2)\string:658--672.
\newblock \href {\doibase 10.2514/1.J055389} {doi: 10.2514/1.J055389}

\bibitem{ni_sensitivity_2017}
Ni A, Wang Q. Sensitivity analysis on chaotic dynamical systems by
  {Non}-{Intrusive} {Least} {Squares} {Shadowing} ({NILSS}). {\it Journal of
  Computational Physics.} 2017\string;347\string:56--77.
\newblock \href {\doibase 10.1016/j.jcp.2017.06.033} {doi:
  10.1016/j.jcp.2017.06.033}

\bibitem{ashley_towards_2019}
Ashley A, Crean J, Hicken J. Towards {Aerodynamic} {Shape} {Optimization} of
  {Unsteady} {Turbulent} {Flows}. In: AIAA Scitech 2019 Forum, AIAA 2019-0168
  2019; San Diego, California

\bibitem{chandramoorthy_analysis_2017}
Chandramoorthy N, Fernandez P, Talnikar C, Wang Q. An {Analysis} of the
  {Ensemble} {Adjoint} {Approach} to {Sensitivity} {Analysis} in {Chaotic}
  {Systems}. In: 23rd AIAA Computational Fluid Dynamics Conference, AIAA
  2017-3799 2017; Denver, Colorado

\bibitem{ashley_aerodynamic_2019}
Ashley A. {\it Aerodynamic {Shape} {Optimization} of {Unsteady}, {Chaotic}
  {Flows}}. PhD thesis. Rensselaer Polytechnic Institute,  2019.

\bibitem{lorenz_deterministic_1963}
Lorenz EN. Deterministic {Nonperiodic} {Flow}. {\it Journal of the Atmospheric
  Sciences.} 1963\string;20(2)\string:130--141.
\newblock \href {\doibase 10.1175/1520-0469(1963)020<0130:DNF>2.0.CO;2} {doi:
  10.1175/1520-0469(1963)020<0130:DNF>2.0.CO;2}

\bibitem{toal_adjoint_2009}
Toal DJJ, Forrester AIJ, Bressloff NW, Keane AJ, Holden C. An adjoint for
  likelihood maximization. {\it Proceedings of the Royal Society A:
  Mathematical, Physical and Engineering Sciences.}
  2009\string;465(2111)\string:3267--3287.
\newblock \href {\doibase 10.1098/rspa.2009.0096} {doi: 10.1098/rspa.2009.0096}

\bibitem{smith_differentiation_1995}
Smith SP. Differentiation of the {Cholesky} {Algorithm}. {\it Journal of
  Computational and Graphical Statistics.} 1995\string;4(2)\string:134--147.
\newblock \href {\doibase 10.1080/10618600.1995.10474671} {doi:
  10.1080/10618600.1995.10474671}

\end{thebibliography}
\appendix


\section{Noise-free closed form likelihood solution} \label{Apdx_NoiseFreeLkd}

The value of $\beta$ that maximizes the marginal log likelihood from \Eq{Eq_ln_lkd_noisy} is given by
\begin{align}
	\p{\ln(L)}{\betavec} 
	&= \onemod^\top \Sigmag^{-1} \fgrad - \beta \onemod^\top \Sigmag^{-1} \onemod = 0 \nonumber \\
	\betavec
	&= \frac{ \onemod^\top \Sigmag^{-1} \fgrad}{ \onemod^\top \Sigmag^{-1} \onemod}. \label{Eq_lkd_beta}
\end{align}
For the case when there are noisy function or gradient evaluations, \ie $\hpstdfval \neq 0$ or $\hpstdfgrad \neq 0$, there is no closed-form solution for $\sigK^2$ that maximizes $\ln(L)$. However, in the noise-free case we have from \Eq{Eq_Sigmag} $\Sigmag = \sigK^2 \left( \Kg + \etaKg \W \right)$, and the marginal log-likelihood from \Eq{Eq_ln_lkd_noisy} simplifies to
\begin{equation} \label{Eq_ln_lkd_noisefree_init}
	\ln(L) 
	= -\frac{n_x (n_d+1) \ln \left( \sigK^2 \right)}{2} 
	- \frac{\ln \left( \det \left( \Kg + \etaKg \W \right) \right)}{2} 
	- \frac{\left( \fgrad - \onemod \beta \right)^\top \left( \Kg + \etaKg \W \right)^{-1} \left( \fgrad - \onemod \beta \right)}{2 \sigK^2},
\end{equation}
where $\det(\sigK^2 \left( \Kg + \etaKg \W \right)) = \sigK^{2 n_x (n_d+1)} \det \left( \Kg + \etaKg \W \right)$. To find the value of $\sigK^2$ that maximizes the marginal log-likelihood for the noise-free case, the marginal log-likelihood is differentiated with respect to $\sigK^2$, equated to zero, and $\sigK^2$ is then isolated:
\begin{align*}
	\p{\ln(L)}{\sigK^2}
	&= -\frac{n_x (n_d+1)}{2 \sigK^2} + \frac{\left( \fgrad - \onemod \beta \right)^\top \left( \Kg + \etaKg \W \right)^{-1} \left( \fgrad - \onemod \beta \right)}{2 \sigK^4} = 0 \\
	\sigK^2 
	&= \frac{\left( \fgrad - \onemod \beta \right)^\top \left( \Kg + \etaKg \W \right)^{-1} \left( \fgrad - \onemod \beta \right)}{n_x (n_d+1)}. \yesnumber \label{Eq_lkd_sigK2_noisefree}
\end{align*}
Substituting $\sigK^2$ from \Eq{Eq_lkd_sigK2_noisefree} into \Eq{Eq_ln_lkd_noisefree_init} and dropping the constant terms gives
\begin{equation} \label{Eq_ln_lkd_noisefree_final}
	\ln(L) = -\frac{n_x (n_d+1)}{2} \ln \left( \sigK^2 \right)- \frac{1}{2} \ln \left( \det \left( \Kg + \etaKg \W \right) \right).
\end{equation}
We seek to maximize the marginal log-likelihood, which is the equivalent of minimizing $-2 \ln(L)$:
\begin{equation}
	\gammavec^* = \argmin_{\gammavec} \left[ n_x (n_d+1) \ln \left( \sigK^2(\gammavec) \right) + \ln \left( \det \left( \Kg(\gammavec) + \etaKg \W \right)\right) \right].
\end{equation}
While the hyperparameters can be optimized using a gradient-free optimizer \cite{toal_kriging_2008,toal_development_2011}, it is more efficient to use a gradient-based optimizer \cite{zhang_exploiting_2005,ollar_gradient_2017}. An adjoint method has also been developed to calculate the required gradients efficiently \cite{toal_adjoint_2009,smith_differentiation_1995,ollar_gradient_2017}. 


\section{Energy method} \label{Sec_AppendixEnergyMethod}

This appendix provides a short overview of the energy method from Ashley $\etal$ to calculate approximate sensitivities for chaotic systems \cite{ashley_towards_2019}.
The time rate of growth of the tangent solution $\vvec$ from \Eq{Eq_dvdt} is now considered in its homogeneous form, \ie $\p{\rvec_x}{s} = \zero$:
\begin{align*}
	\vvec^\top \left( \frac{d \vvec}{dt} + \p{\rvec_x}{\uvec} \vvec \right)
	&= \frac{1}{2} \frac{d \|\vvec \|_2^2}{dt} + \vvec^\top \p{\rvec_x}{\uvec} \vvec \\
	&= \frac{1}{2} \frac{d \|\vvec \|_2^2}{dt} 
	+ \frac{1}{2} \vvec^\top \left[ \p{\rvec_x}{\uvec} + \left( \p{\rvec_x}{\uvec} \right)^\top \right] \vvec 
	+ \frac{1}{2} \vvec^\top \left[ \p{\rvec_x}{\uvec} - \left( \p{\rvec_x}{\uvec} \right)^\top \right] \vvec \\
	&= \frac{1}{2} \frac{d \|\vvec \|_2^2}{dt} 
	+ \frac{1}{2} \vvec^\top \left[ \p{\rvec_x}{\uvec} + \left( \p{\rvec_x}{\uvec} \right)^\top \right] \vvec = 0, \yesnumber \label{Eq_dvdt_growth}
\end{align*}
where the skew-symmetric portion of $\p{\rvec_x}{\uvec}$ does not lead to any energy growth in $\vvec$. From \Eq{Eq_dvdt_growth} it is clear that ensuring energy stability requires that 
\begin{equation}
	\vvec^\top \left[ \p{\rvec_x}{\uvec} + \left( \p{\rvec_x}{\uvec} \right)^\top \right] \vvec \geq 0.
\end{equation}

The Jacobian $\p{\rvec_x}{\uvec}$ can be modified with the addition of a matrix in order to ensure the following inequality is satisfied:
\begin{equation} \label{Eq_growth_vvec_ineq}
	\vvec^\top \left( \p{\rvec_x}{\uvec} + \A \right) \vvec \geq 0,
\end{equation}
where $\A$ is a symmetric matrix. The matrix $\A$ can be selected to ensure \Eq{Eq_growth_vvec_ineq} is satisfied by first calculating the eigenvalue decomposition
\begin{equation} \label{Eq_vvec_growth_eigdecomp}
	\E \Lambda \E^\top = \p{\rvec_x}{\uvec} + \left( \p{\rvec_x}{\uvec} \right)^\top,
\end{equation}
where the columns of $\E$ are eigenvectors and the diagonal entries of $\Lambda$ are the eigenvalues. To ensure \Eq{Eq_growth_vvec_ineq} is satisfied, the matrix $\A$ is calculated with
\begin{equation} \label{Eq_vvec_growth_A_eig}
	\A = -\E \Lambda^{-} \E^\top,
\end{equation}
where the diagonal matrix $\Lambda^{-}$ holds the non-positive eigenvalues of $\Lambda$, with the positive eigenvalues simply set to zero. 

\section{Supplemental information on the studies for the Bayesian optimizer} \label{Sec_AppendixUnconStudy}

\subsection{Selecting the maximum condition number $\condmax$} \label{Sec_AppendixUnconStudy_CondMax}

\begin{figure}[t!]
	\centering
	\begin{subfigure}[t]{0.328\textwidth}
		\includegraphics[width=\textwidth]{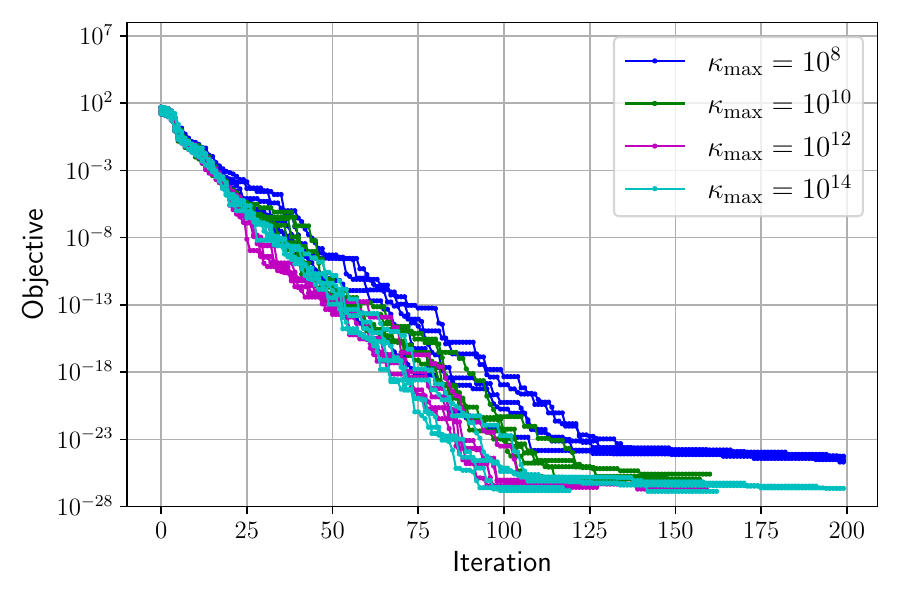}
		\caption{Objective: quadratic function}
		\label{Fig_Study_CondMax_Obj_Quad_d20}
	\end{subfigure}	
	\begin{subfigure}[t]{0.328\textwidth}
		\includegraphics[width=\textwidth]{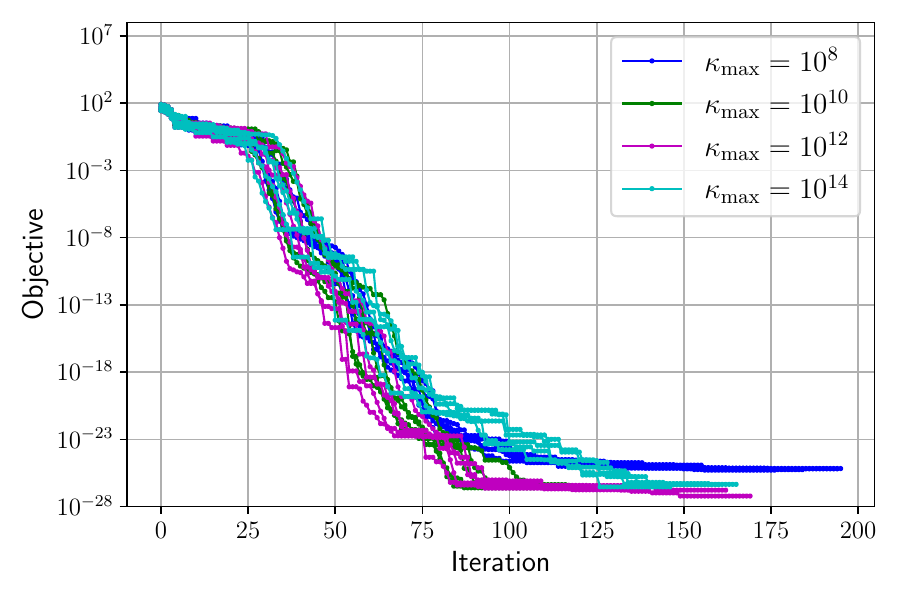}
		\caption{Objective: bowl function}
		\label{Fig_Study_CondMax_Obj_Bowl_d20}
	\end{subfigure}	
	\begin{subfigure}[t]{0.328\textwidth}
		\includegraphics[width=\textwidth]{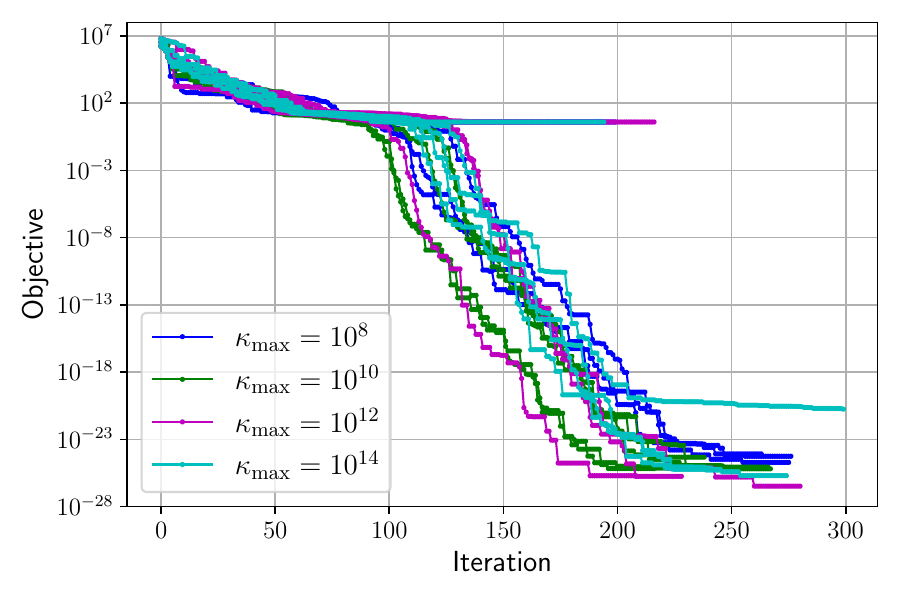}
		\caption{Objective: Rosenbrock $a=100$}
		\label{Fig_Study_CondMax_Obj_RosenA100_d20}
	\end{subfigure}	
	\begin{subfigure}[t]{0.328\textwidth}
		\includegraphics[width=\textwidth]{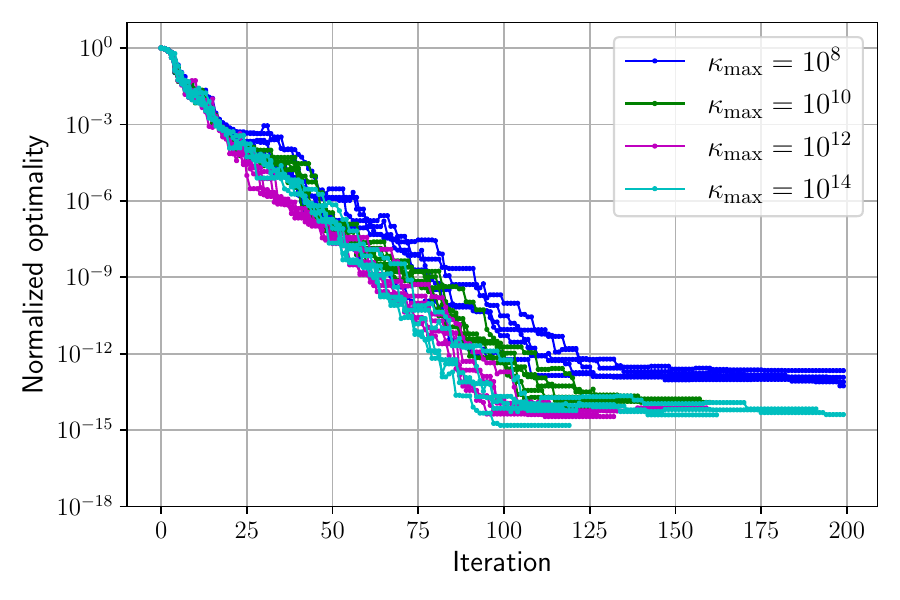}
		\caption{Optimality: quadratic function}
		\label{Fig_Study_CondMax_Opt_Quad_d20}
	\end{subfigure}	
	\begin{subfigure}[t]{0.328\textwidth}
		\includegraphics[width=\textwidth]{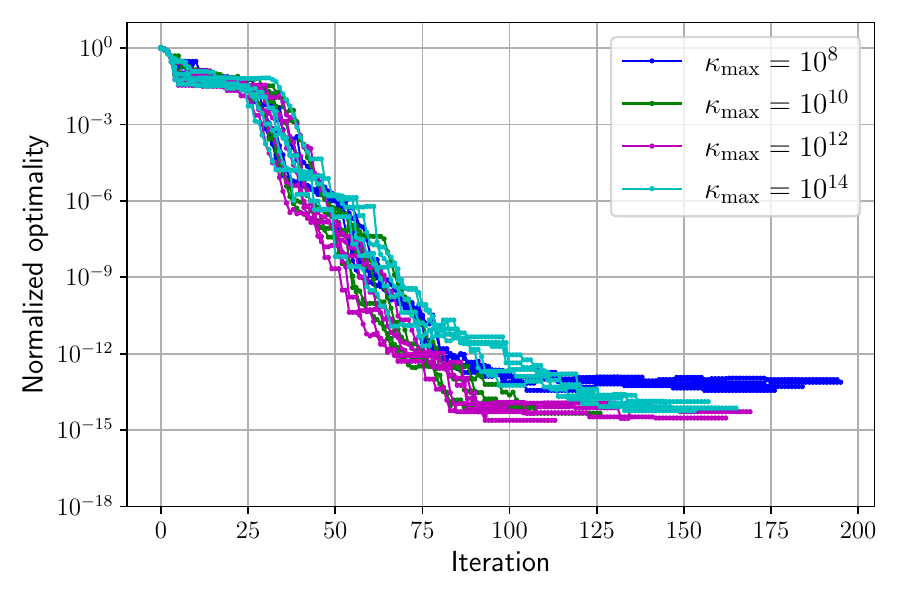}
		\caption{Optimality: bowl function}
		\label{Fig_Study_CondMax_Opt_Bowl_d20}
	\end{subfigure}	
	\begin{subfigure}[t]{0.328\textwidth}
		\includegraphics[width=\textwidth]{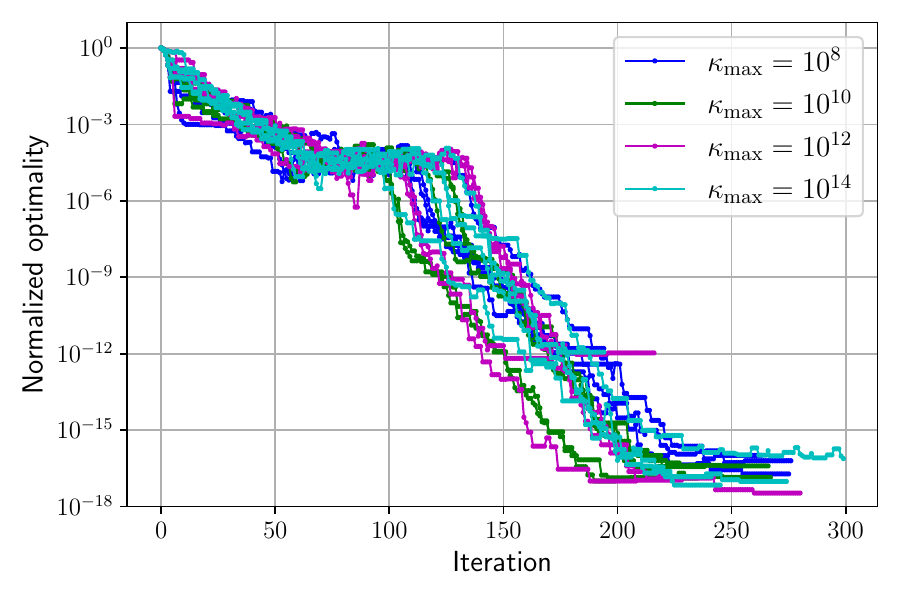}
		\caption{\mbox{Optimality: Rosenbrock $a=100$}}
		\label{Fig_Study_CondMax_Opt_RosenA100_d20}
	\end{subfigure}	
	\caption[Unconstrained optimization for different $\condmax$ with $n_d=20$.]{Unconstrained Bayesian optimization with the preconditioning method and different $\condmax$. The test cases are the $n_d=20$ quadratic, bowl, and $a=100$ Rosenbrock functions from \Eqss{Eq_Quadratic_fun}{Eq_Bowl}{Eq_Rosenbrock}, respectively.}
	\label{Fig_Study_CondMax}
\end{figure}

The Bayesian optimizer with the preconditioning method with different values of $\condmax$ is applied to the twenty-dimensional quadratic, bowl, and $a=100$ Rosenbrock functions from \Eqss{Eq_Quadratic_fun}{Eq_Bowl}{Eq_Rosenbrock}, respectively, and the results are plotted in \Fig{Fig_Study_CondMax}. It is clear that the convergence results for the Bayesian optimizer using $\condmax \in \{10^8, 10^{10}, 10^{12}, 10^{14} \}$ are similar for the three test cases with $n_d = 20$.

\begin{figure}[t!]
	\centering
	\begin{subfigure}[t]{0.328\textwidth}
		\includegraphics[width=\textwidth]{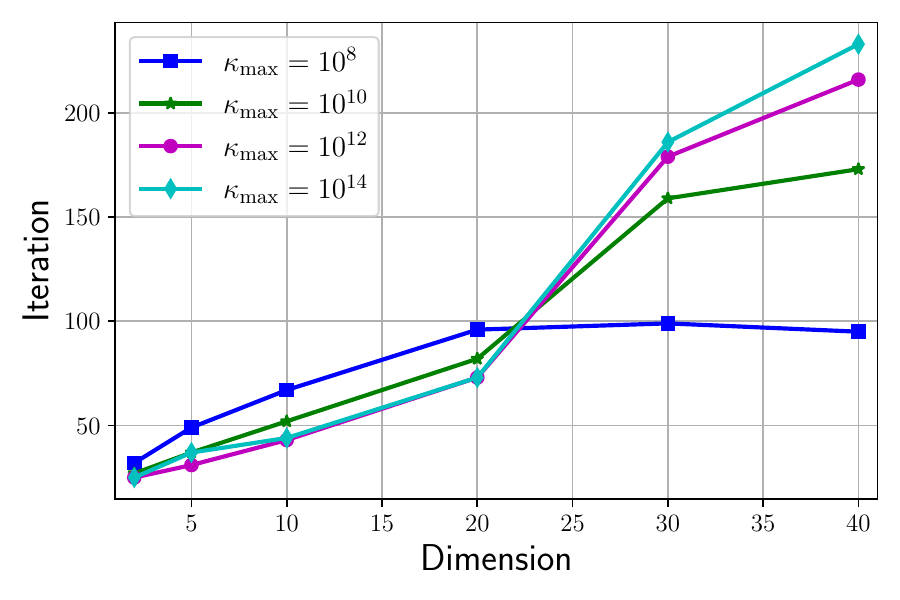}
		\caption{Quadratic function: \Eq{Eq_Quadratic_fun}}
		\label{Fig_Study_CondMax_OptTol_Quad}
	\end{subfigure}	
	\begin{subfigure}[t]{0.328\textwidth}
		\includegraphics[width=\textwidth]{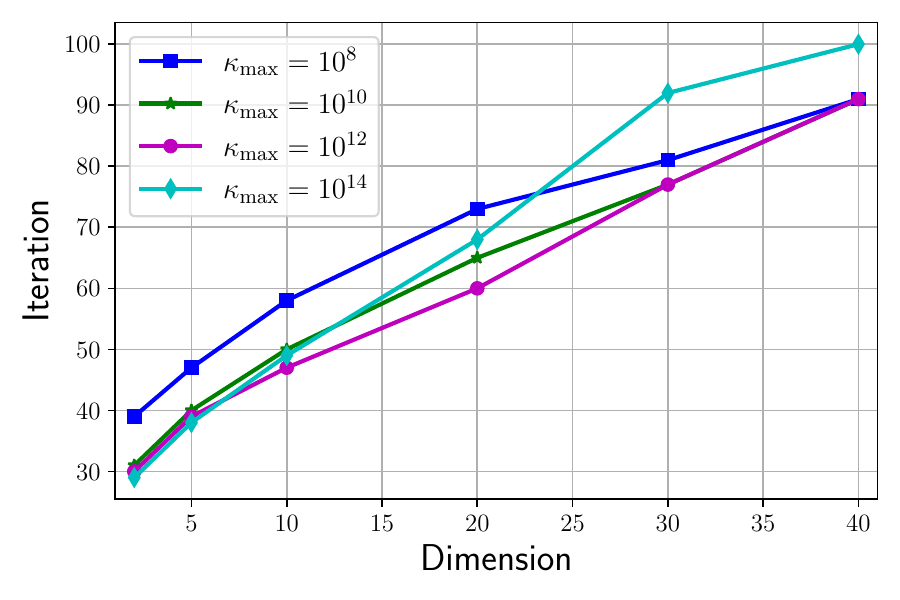}
		\caption{Bowl function: \Eq{Eq_Bowl}}
		\label{Fig_Study_CondMax_OptTol_Bowl}
	\end{subfigure}	
	\begin{subfigure}[t]{0.328\textwidth}
		\includegraphics[width=\textwidth]{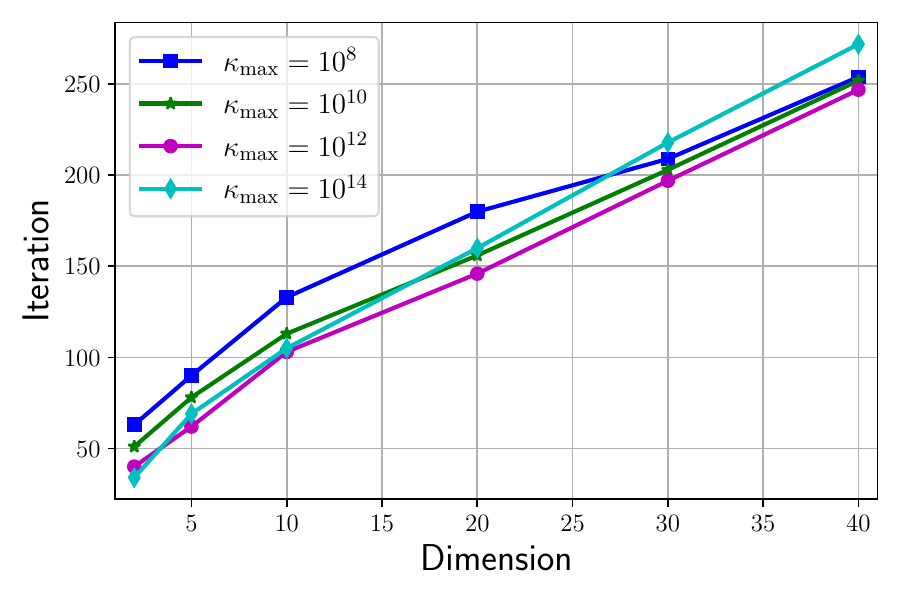}
		\caption{Rosenbrock $a=100$: \Eq{Eq_Rosenbrock}}
		\label{Fig_Study_CondMax_OptTol_RosenA100}
	\end{subfigure}	
	\caption[Median number of iterations using 25 independent optimization runs for the Bayesian optimizer with different $\condmax$ to reduce the optimality by 10 orders of magnitude and the objective evaluations below $10^{-5}$.]{Median number of iterations for the Bayesian optimizer to reduce the optimality by 10 orders of magnitude and the objective evaluations below $10^{-5}$.}
	\label{Fig_Study_CondMax_OptTol}
\end{figure}

\Fig{Fig_Study_CondMax_OptTol} shows the median number of iterations required for the Bayesian optimizer using different values of $\condmax$ to reduce the objective evaluations below $10^{-5}$ and the optimality 10 orders of magnitude for the three test cases. In general, the Bayesian optimizer using $\condmax = 10^{10}$ or $\condmax = 10^{12}$ requires the fewest iterations to achieve the desired tolerance. The use of $\condmax = 10^{8}$ or $\condmax = 10^{14}$ results in inferior results, except for the former on the quadratic function with $n_d = 30$ and $n_d =40$. It other tests it was found that using $\condmax = 10^6$ results in even slower convergence for the Bayesian optimizer while using $\condmax \geq 10^{15}$ results in failures for the Cholesky decomposition. Therefore, $\condmax = 10^{10}$ was selected since it provided consistent results for the test cases considered.

\subsection{Selecting the kernel} \label{Sec_AppendixUnconStudy_Kernel}

\begin{figure}[t!]
	\centering
	\begin{subfigure}[t]{0.328\textwidth}
		\includegraphics[width=\textwidth]{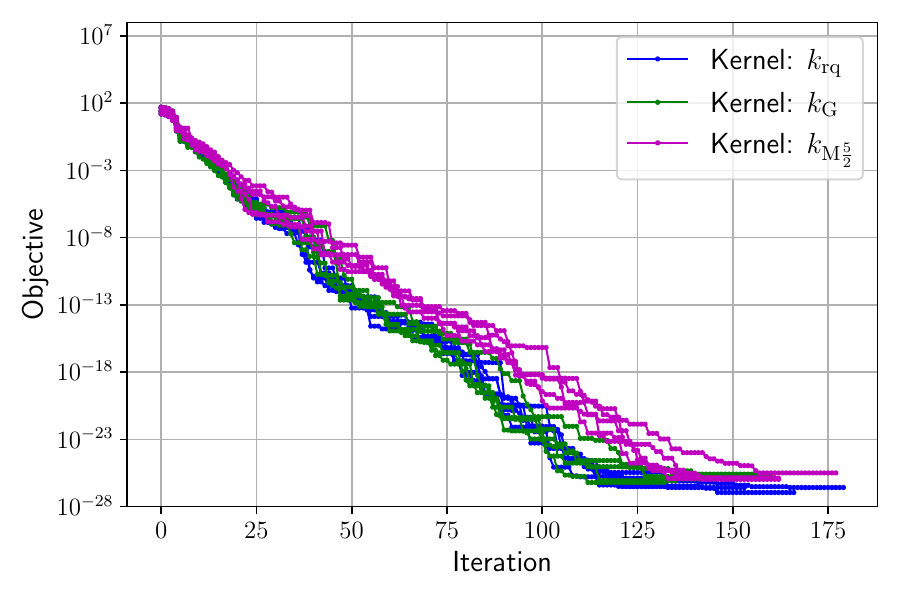}
		\caption{Objective: quadratic function}
		\label{Fig_Study_Kernel_Obj_Quad_d20}
	\end{subfigure}	
	\begin{subfigure}[t]{0.328\textwidth}
		\includegraphics[width=\textwidth]{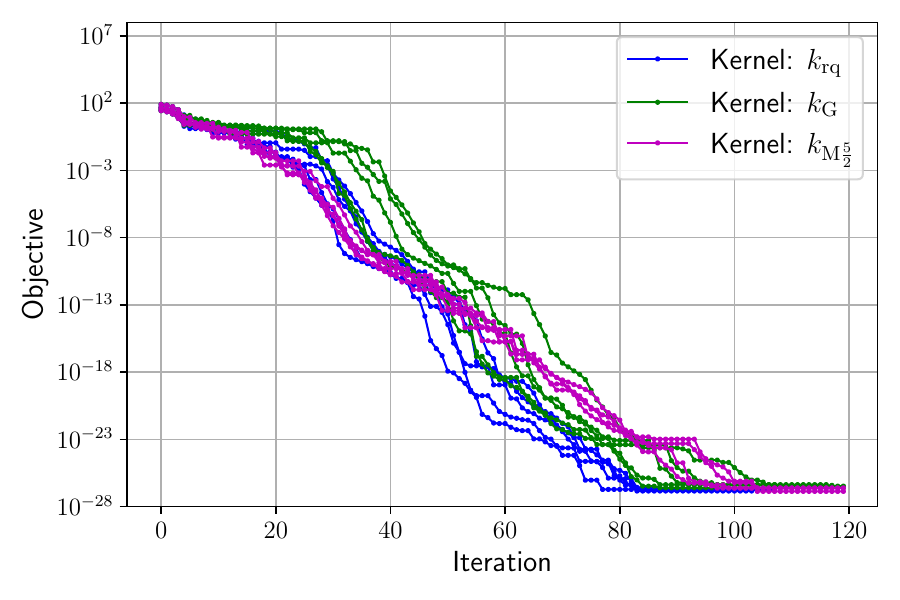}
		\caption{Objective: bowl function}
		\label{Fig_Study_Kernel_Obj_Bowl_d20}
	\end{subfigure}	
	\begin{subfigure}[t]{0.328\textwidth}
		\includegraphics[width=\textwidth]{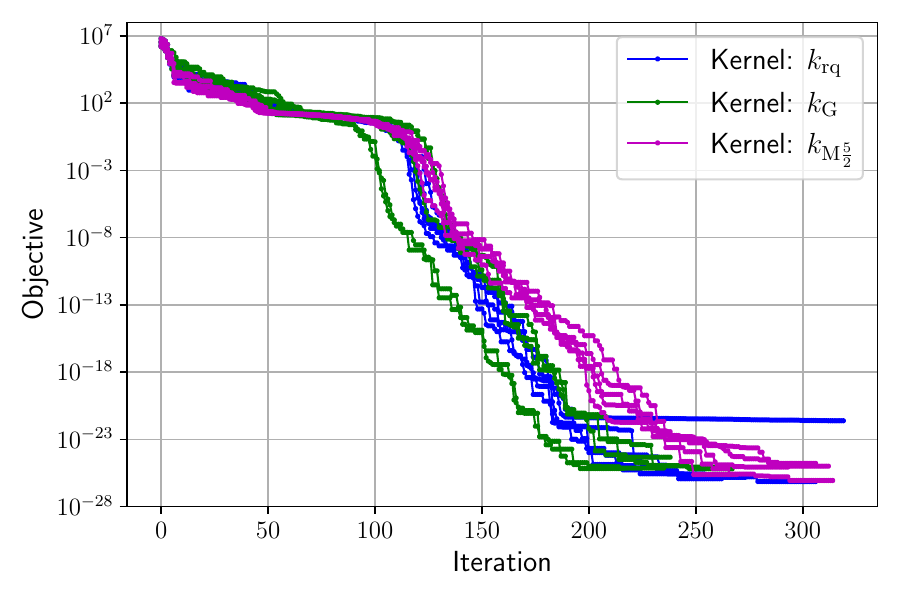}
		\caption{Objective: Rosenbrock $a=100$}
		\label{Fig_Study_Kernel_Obj_RosenA100_d20}
	\end{subfigure}	
	\begin{subfigure}[t]{0.328\textwidth}
		\includegraphics[width=\textwidth]{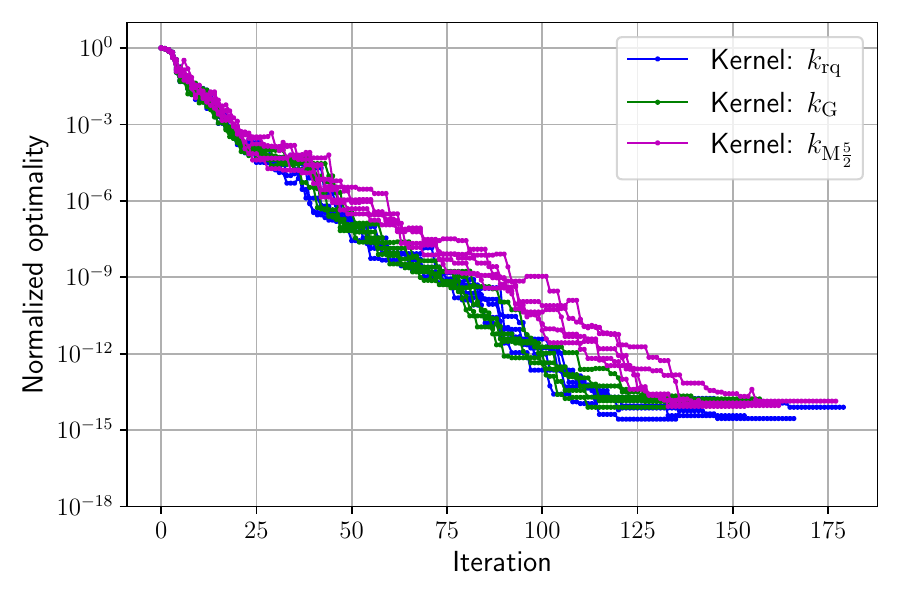}
		\caption{Optimality: quadratic function}
		\label{Fig_Study_Kernel_Opt_Quad_d20}
	\end{subfigure}	
	\begin{subfigure}[t]{0.328\textwidth}
		\includegraphics[width=\textwidth]{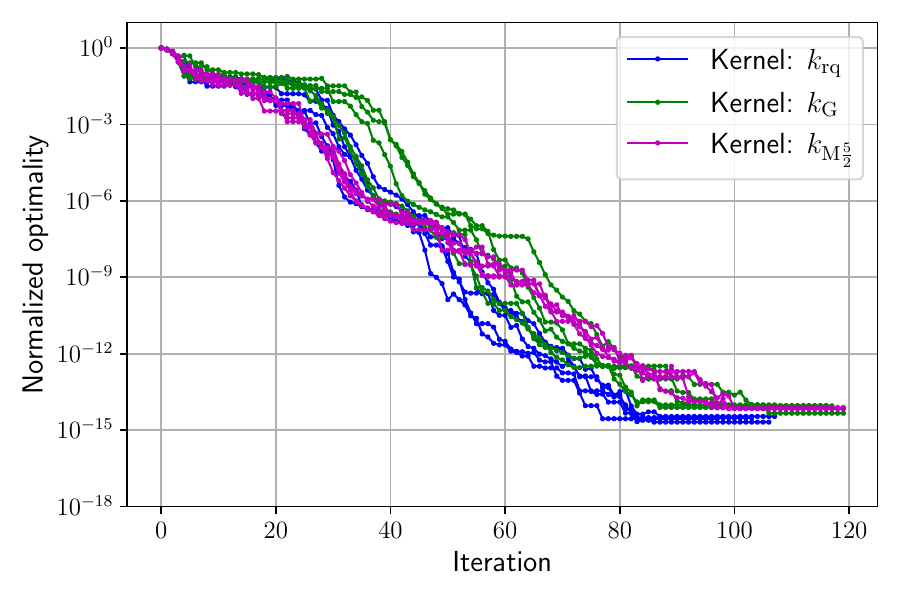}
		\caption{Optimality: bowl function}
		\label{Fig_Study_Kernel_Opt_Bowl_d20}
	\end{subfigure}	
	\begin{subfigure}[t]{0.328\textwidth}
		\includegraphics[width=\textwidth]{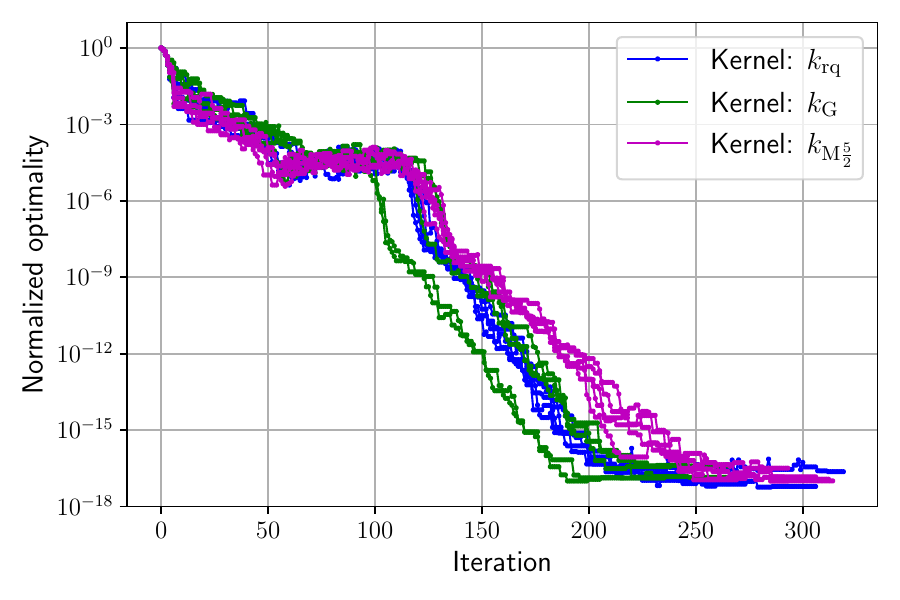}
		\caption{\mbox{Optimality: Rosenbrock $a=100$}}
		\label{Fig_Study_Kernel_Opt_RosenA100_d20}
	\end{subfigure}	
	\caption{Unconstrained study for the Bayesian optimizer with the Gaussian, Mat\'ern $\frac{5}{2}$, and rational quadratic kernels from \Eqss{Eq_kern_Gaussian}{Eq_kern_Mat5f2}{Eq_kern_RatQd}, respectively. The test cases are the quadratic, bowl, and Rosenbrock $a=100$ functions from \Sec{Sec_LocalOptz_TestCases} with $n_d=20$.}
	\label{Fig_Study_Kernel}
\end{figure}

The use of the Gaussian, Mat\'ern $\frac{5}{2}$, and rational quadratic kernels from \Eqss{Eq_kern_Gaussian}{Eq_kern_Mat5f2}{Eq_kern_RatQd}, respectively, is investigated in this subsection. \Fig{Fig_Study_Kernel} shows the Bayesian optimizer with these three kernels minimizing the $n_d=20$ quadratic, bowl, and $a=100$ Rosenbrock functions from \Eqss{Eq_Quadratic_fun}{Eq_Bowl}{Eq_Rosenbrock}, respectively. The results indicate that using a different kernel does not have a substantial impact on the performance of the Bayesian optimizer for the test cases considered.

\begin{figure}[t!]
	\centering
	\begin{subfigure}[t]{0.328\textwidth}
		\includegraphics[width=\textwidth]{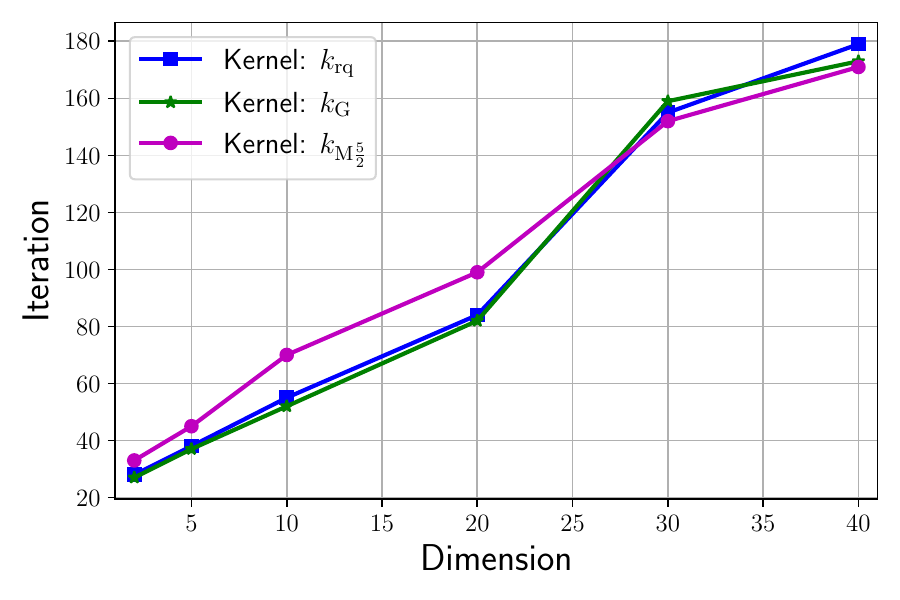}
		\caption{Quadratic function: \Eq{Eq_Quadratic_fun}}
		\label{Fig_Study_Kernel_OptTol_Quad}
	\end{subfigure}	
	\begin{subfigure}[t]{0.328\textwidth}
		\includegraphics[width=\textwidth]{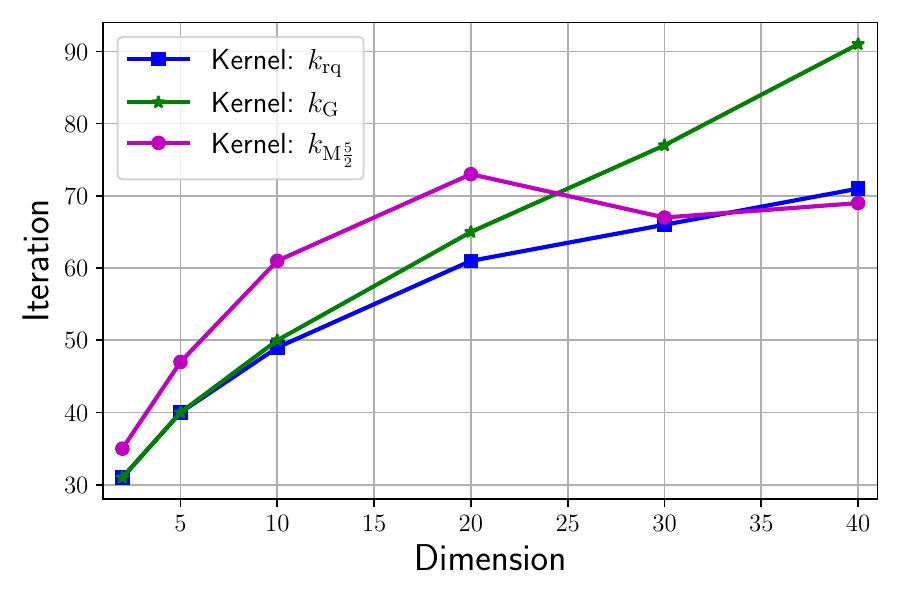}
		\caption{Bowl function: \Eq{Eq_Bowl}}
		\label{Fig_Study_Kernel_OptTol_Bowl}
	\end{subfigure}	
	\begin{subfigure}[t]{0.328\textwidth}
		\includegraphics[width=\textwidth]{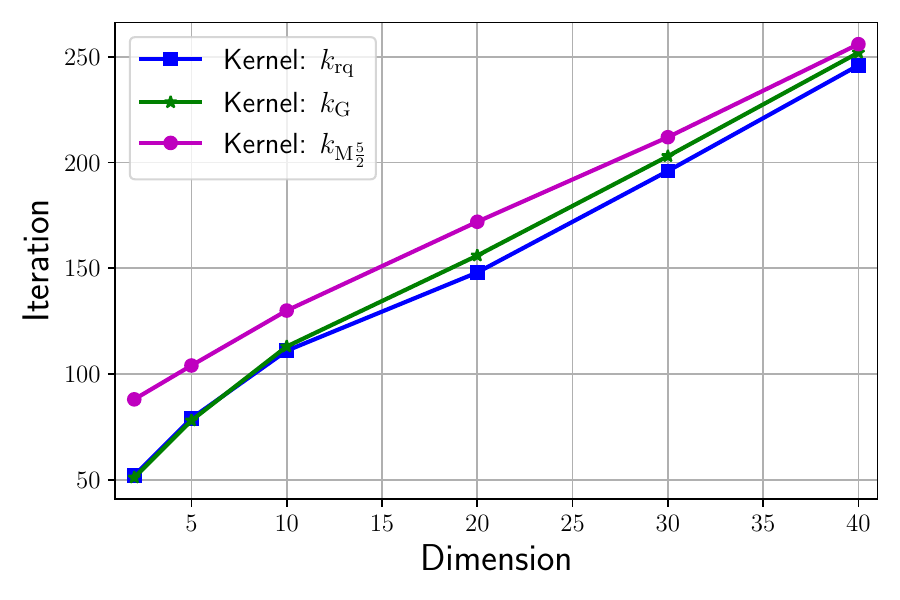}
		\caption{Rosenbrock $a=100$: \Eq{Eq_Rosenbrock}}
		\label{Fig_Study_Kernel_OptTol_RosenA100}
	\end{subfigure}	
	\caption{Median number of iterations to reduce the objective function below $10^{-5}$ and the optimality by 10 orders of magnitude for the Bayesian optimizer using different kernels for 25 independent optimization runs.}
	\label{Fig_Study_Kernel_OptTol}
\end{figure}

For each test case and kernel 25 independent optimization runs were performed and the median number of iterations required to reduce the objective function below $10^{-5}$ and achieve a 10-order reduction in the optimality for the Bayesian optimizer using each of the three kernels is shown in \Fig{Fig_Study_Kernel_OptTol} with the three test cases and $2 \leq n_d \leq 40$. Once again, there is not a significant difference between the performance of the Bayesian optimizer using different kernels. In general, the Bayesian optimizer requires additional iterations to achieve the desired tolerance when it uses the Mat\'ern $\frac{5}{2}$ kernel. All of the test cases considered are infinitely continuously differentiable, just like the Gaussian and rational quadratic kernels. However, as the Mat\'ern $\frac{5}{2}$ kernel is only twice continuously differentiable, it may be advantageous for problems that are not infinitely differentiable.

For these test cases that are infinitely differentiable, the Gaussian kernel was selected as the default kernel for the Bayesian optimizer since it performs well and is simpler to implement than the rational quadratic kernel, which has the additional hyperparameter $\alpha$ that needs to be selected at each iteration.

\end{document}